\documentclass[a4,10pt]{article}
\usepackage{amsmath,amscd,amssymb}

\newcommand{\nbige}{\mathcal{E}}
\newcommand{\nbigf}{\mathcal{F}}
\newcommand{\nbigg}{\mathcal{G}}
\newcommand{\nbigh}{\mathcal{H}}
\newcommand{\nbigi}{\mathcal{I}}

\newcommand{\nbigo}{\mathcal{O}}
\newcommand{\nbigp}{\mathcal{P}}

\newcommand{\nbigv}{\mathcal{V}}

\newcommand{\nbigz}{\mathcal{Z}}

\newcommand{\proj}{\mathbb{P}}
\newcommand{\seisuu}{{\mathbb Z}}

\newcommand{\cnum}{{\mathbb C}}
\newcommand{\real}{{\mathbb R}}

\newcommand{\DD}{\mathbb{D}}

\newcommand{\gminia}{\mathfrak a}

\newcommand{\vecv}{{\boldsymbol v}}
\newcommand{\vecu}{{\boldsymbol u}}
\newcommand{\vecw}{{\boldsymbol w}}

\newcommand{\veczero}{{\boldsymbol 0}}

\newcommand{\veca}{{\boldsymbol a}}
\newcommand{\vecb}{{\boldsymbol b}}

\newcommand{\vecc}{{\boldsymbol c}}

\newcommand{\vecn}{{\boldsymbol n}}
\newcommand{\vecp}{{\boldsymbol p}}

\newcommand{\lrarr}{\longrightarrow}

\newcommand{\pf}{{\bf Proof}\hspace{.1in}}
\newcommand{\qed}{\mbox{\rule{1.2mm}{3mm}}}

\def\Hom{\mathop{\rm Hom}\nolimits}

\def\End{\mathop{\rm End}\nolimits}

\def\Image{\mathop{\rm Im}\nolimits}

\def\Re{\mathop{\rm Re}\nolimits}

\def\rank{\mathop{\rm rank}\nolimits}

\def\tr{\mathop{\rm tr}\nolimits}
\def\Tr{\mathop{\rm Tr}\nolimits}
\def\vol{\mathop{\rm dvol}\nolimits}
\def\dvol{\mathop{\rm dvol}\nolimits}

\def\id{\mathop{\rm id}\nolimits}

\newcommand{\del}{\partial}
\newcommand{\delbar}{\overline{\del}}

\newcommand{\barz}{\overline{z}}
\newcommand{\zbar}{\barz}
\newcommand{\zetabar}{\overline{\zeta}}
\newcommand{\baralpha}{\overline{\alpha}}
\newcommand{\alphabar}{\baralpha}

\newcommand{\etabar}{\overline{\eta}}

\newcommand{\xbar}{\overline{x}}
\newcommand{\sbar}{\overline{s}}
\newcommand{\Hbar}{\overline{H}}

\newcommand{\lefttop}[1]{{}^{#1}\!}

\newcommand{\closedclosed}[2]{[#1,#2]}

\newcommand{\omegatilde}{\widetilde{\omega}}

\newcommand{\Etilde}{\widetilde{E}}

\newcommand{\thetatilde}{\widetilde{\theta}}

\newcommand{\Vtilde}{\widetilde{V}}

\newcommand{\htilde}{\widetilde{h}}

\newcommand{\vtilde}{\widetilde{v}}

\newcommand{\Ftilde}{\widetilde{F}}

\newcommand{\atilde}{\widetilde{a}}

\newcommand{\DDtilde}{\widetilde{\DD}}

\newcommand{\Utilde}{\widetilde{U}}
\newcommand{\Dtilde}{\widetilde{D}}
\newcommand{\Xtilde}{\widetilde{X}}

\newcommand{\Ltilde}{\widetilde{L}}

\def\op{\mathop{\rm op}\nolimits}

\newcommand{\Ptilde}{\widetilde{P}}

\newcommand{\Qtilde}{\widetilde{Q}}

\newcommand{\vecatilde}{\widetilde{\veca}}

\newcommand{\DDD}{\boldsymbol D}

\newcommand{\Gbar}{\overline{G}}

\newcommand{\vecell}{\boldsymbol \ell}
\newcommand{\vecq}{\boldsymbol q}

\newcommand{\etilde}{\widetilde{e}}
\newcommand{\etatilde}{\widetilde{\eta}}
\newcommand{\nbiggbar}{\overline{\nbigg}}

\newtheorem{thm}{Theorem}[section]
\newtheorem{cor}[thm]{Corollary}

\newtheorem{rem}[thm]{Remark}
\newtheorem{lem}[thm]{Lemma}
\newtheorem{prop}[thm]{Proposition}

\newtheorem{assumption}[thm]{Assumption}

\begin{document}

\title{Asymptotic behaviour of
 certain families of harmonic bundles
 on Riemann surfaces}
\author{Takuro Mochizuki}
\date{}
\maketitle

\begin{abstract}
Let $(E,\overline{\partial}_E,\theta)$ be a stable Higgs bundle
of degree $0$ on a compact connected Riemann surface.
Once we fix a flat metric $h_{\det(E)}$
on the determinant of $E$,
we have the harmonic metrics $h_t$ $(t>0)$
for the stable Higgs bundles 
$(E,\overline{\partial}_E,t\theta)$
such that $\det(h_t)=h_{\det(E)}$.
We study the behaviour of $h_t$ when $t$ goes to $\infty$.
First, we  show that the Hitchin equation is
asymptotically decoupled 
under the assumption that the Higgs field
is generically regular semisimple.
We apply it to the study of 
the so called Hitchin WKB-problem.
Second, 
we study the convergence of the sequence
$(E,\overline{\partial}_E,\theta,h_t)$
in the case $\rank E=2$.
We introduce a rule 
to determine the parabolic weights
of a ``limiting configuration'',
and we show the convergence of the sequence
to the limiting configuration in an appropriate sense.
The results can be appropriately generalized
in the context of 
Higgs bundles with a Hermitian-Einstein metric
on curves.

\vspace{.1in}


\end{abstract}

\section{Introduction}

Let $X$ be a compact connected Riemann surface.
Let $(E,\delbar_E,\theta)$ be a Higgs bundle of rank $r$
on $X$ with $\deg(E)=0$.
Let $h$ be a Hermitian metric of $E$.
We have the Chern connection $\nabla_h$
associated to $(E,\delbar_E,h)$.
Let $R(h)$ denote the curvature of $\nabla_h$.
Let $\theta_h^{\dagger}$ denote the adjoint of
$\theta$ with respect to $h$.
Recall the celebrated Hitchin equation 
\cite{Hitchin-self-duality}:
\[
 R(h)+\bigl[\theta,\theta_h^{\dagger}\bigr]=0
\]
If the Hitchin equation is satisfied,
$h$ is called a harmonic metric of $(E,\delbar_E,\theta)$,
and $(E,\delbar_E,\theta,h)$ is called
a harmonic bundle.

\begin{rem}
The Hitchin equation
makes sense for
a Higgs bundle with a Hermitian metric
on any complex curve.
In this introduction, 
$X$ is assumed to be compact
to simplify the explanation.

If $\deg(E)$ is not $0$,
we take a Hermitian metric $h_{\det(E)}$ of
the determinant line bundle $\det(E)$,
and we consider the Hermitian-Einstein condition
$R(h)^{\bot}+[\theta,\theta^{\dagger}_h]=0$,
where $R(h)^{\bot}$ is the trace-free part of $R(h)$.
The condition is also called the Hitchin equation.
In this introduction, we assume $\deg(E)=0$
for simplicity.
\hfill\qed
\end{rem}

If $\rank E=1$,
we always have $[\theta,\theta_h^{\dagger}]=0$.
Hence, the Hitchin equation is reduced to
$R(h)=0$,
i.e., the metric $h$ is flat with respect to
the Chern connection.
By the classical harmonic theory,
we can always find such a harmonic metric
in the rank one case,
which is unique up to the multiplication of
positive constants.

In the higher rank case,
we fix a harmonic metric $h_{\det E}$
of $(\det(E),\delbar_{\det E},\tr\theta)$.
According to Hitchin \cite{Hitchin-self-duality}
and Simpson \cite{Simpson88},
if the Higgs bundle $(E,\delbar_E,\theta)$ is stable,
we have a unique harmonic metric $h$ of
$(E,\delbar_E,\theta)$
such that $\det(h)=h_{\det E}$.

For any non-zero complex number $t$,
the Higgs bundle
$(E,\delbar_E,t\theta)$ is also stable.
We obtain a family of harmonic metrics
$h_{t}$ $(t\in\cnum^{\ast})$
satisfying $\det(h_t)=h_{\det E}$.
It is easy to observe that
$h_{t_1}=h_{t_2}$ if $|t_1|=|t_2|$.
So, it is enough to consider the case
where $t$ are positive numbers.

Simpson studied the behaviour of
$(E,\delbar_E,t\theta,h_{t})$ when $t\to 0$.
(See \cite{Simpson88,Simpson90,s5,s4}, for example.)
He discovered 
the convergence to a polarized variation of Hodge structure,
and he gave various applications of this interesting phenomena.

More recently,
there has been a growing interest to the behaviour 
of $(E,\delbar_E,t\theta,h_{t})$ when $t\to\infty$.
In \cite{KNPS},
Katzarkov, Noll, Pandit and Simpson
proposed ``Hitchin WKB-problem''
on the behaviour of 
the family of harmonic metrics $h_t$
and the monodromy of the associated flat connections
$\nabla_h+\theta+\theta^{\dagger}_h$,
in relation with their magnificent theory of
harmonic maps to buildings.
In \cite{MSWW,MSWW2},
Mazzeo, Swoboda, Weiss and Witt
studied the rank $2$ case
under the assumption that
the zeroes of $\det(\theta)-(\tr\theta)^2/4$
are simple,
i.e.,
the spectral curve of the Higgs field
is smooth irreducible and simply ramified over $X$,
motivated by the study on the structure of the end
of the moduli spaces of Higgs bundles.
They introduced the concept of
``limiting configuration'',
and they proved the convergence 
to the limiting configuration
under the assumption,
inspired by the work of 
Gaiotto, Moore and Neitzke
\cite{GMN10,GMN13}.
In \cite{Collier-Li},
Collier and Li closely studied the issue 
for some Toda-like harmonic bundles
in a rather explicit way,
and they resolved Hitchin WKB-problem
in these cases
for some kind of non-critical paths.

In this paper,
we shall give two results
on the asymptotic behaviour of 
the harmonic bundles
$(E,\delbar_E,t\theta,h_t)$ $(t\to\infty)$.
One is the asymptotic decoupling,
and the other is the convergence to the limiting configuration
for harmonic bundles of rank two.

\subsection{Asymptotic decoupling}

The Hitchin equation is much simplified
if $R(h)=[\theta,\theta_h^{\dagger}]=0$ holds.
The equation $R(h)=0$ implies that 
$(E,\nabla_h,h)$ is a unitary flat bundle.
The additional condition $[\theta,\theta_h^{\dagger}]=0$ implies that,
at least locally, 
we have a flat decomposition
$(E,\nabla,h)=
 \bigoplus_{i=1}^r (E_i,\nabla_i,h_i)$
into flat line bundles
such that $\theta=\bigoplus \phi_i\cdot \id_{E_i}$,
where $\phi_i$ are holomorphic one forms.
In \cite{MSWW},
this kind of simplification seems to be called 
``decoupling'' of the Hitchin equation.

Our first purpose is to show that 
if $t$ is sufficiently large,
the Hitchin equation for 
$(E,\delbar_E,t\theta)$
is almost decoupled in some sense.

\subsubsection{Generically regular semisimple Higgs bundles}
\label{subsection;15.7.22.1}

To state the claim more precisely,
we introduce a condition for Higgs bundles.
Let $(E,\delbar_E,\theta)$ be a Higgs bundle on $X$.
We have the associated coherent sheaf
$M_{E,\theta}$ on the cotangent bundle $T^{\ast}X$.
The support $\Sigma(E,\theta)$ is called
the spectral curve of the Higgs bundle.
The number of the points 
of $\rho(P):=T_P^{\ast}X\cap \Sigma(E,\theta)$
are finite for any $P\in X$.
We say that 
the Higgs bundle $(E,\delbar_E,\theta)$
is generically regular semisimple if the following holds:
\begin{itemize}
\item
We have a discrete subset $D\subset X$
such that  $\rho(P)=\rank E$
for any $P\in X\setminus D$. 
\end{itemize}
Let $D(E,\theta)$ denote the set of 
the points $P\in X$ such that $\rho(P)<\rank E$,
which we call the discriminant of the Higgs bundle.

Suppose that $(E,\delbar_E,\theta,h)$
is generically regular semisimple.
Then, the following holds
for any point $P\in X\setminus D(E,\theta)$
with a small neighbourhood $U_P$ of $P$:
\begin{itemize}
\item
We have holomorphic $1$-forms
$\phi_{P,1},\ldots,\phi_{P,r}$ on $U_P$
and a decomposition of the Higgs bundle
\[
 (E,\delbar_E,\theta)_{|U_P}
=\bigoplus_{i=1}^r
 (E_{P,i},\delbar_{E_{P,i}},\phi_{P,i}\id_{E_{P,i}}),
\]
where we assume that
$\rank E_i=1$ $(i=1,\ldots,r)$,
and that $\phi_{P,i}-\phi_{P,j}$ $(i\neq j)$
have no zero.
\end{itemize}

\subsubsection{Asymptotic decoupling}
\label{subsection;15.7.20.30}

Let $(E,\delbar_E,\theta)$ be 
a stable Higgs bundle of degree $0$ on $X$.
Suppose that it is generically regular semisimple.
We take any K\"ahler metric $g_X$ of $X$.
For any local section $s$ of $\End(E)\otimes\Omega^{p,q}$,
we have the function 
$|s|_{h_t,g_X}:X\lrarr \real$,
where $|s|_{h_t,g_X}(P)$ $(P\in X)$
are the norm of $s_{|P}$ with respect to $h_t$ and $g_X$.
We have the asymptotic decoupling in the following sense.

\begin{thm}[Theorem
 \ref{thm;15.7.12.11}]
\label{thm;15.7.20.1}
We take any neighbourhood $N$ of 
the discriminant $D(E,\theta)$.
Then, there exist positive constants $C_0$ and $\epsilon_0$
such that the following holds on $X\setminus N$:
\[
 \bigl|R(h_t)\bigr|_{h_t,g_X}
=|t|^2\bigl|
 \bigl[\theta,\theta^{\dagger}_{h_t}\bigr]
 \bigr|_{h_t,g_X}
\leq
 C_0\exp\bigl(-\epsilon_0t\bigr).
\]
The constants $C_0$ and $\epsilon_0$
may depend only on $(X,g_X)$, $N$ and $\Sigma(E,\theta)$.
\end{thm}

We also have the family of flat connections
$\DD^1_{h_t}:=\nabla_{h_t}+t\theta+(t\theta)^{\dagger}_{h_t}$,
which are correctly associated to
the harmonic bundles $(E,\delbar_E,t\theta,h_t)$.
Let us describe that we have nice approximations
of these connections.

Let $P$ be any point of $X\setminus D(E,\theta)$.
We take a small neighbourhood $U_P$
of $P$ in $X\setminus D(E,\theta)$.
We have a decomposition of the Higgs bundle
$(E,\delbar_E,\theta)_{|U_P}
=\bigoplus_{i=1}^r(E_{P,i},\delbar_{E_{P,i}},\theta_{P,i})$,
where $\rank E_{P,i}=1$.
Let $h_{t,E_{P,i}}$ be the restriction of $h_t$ to $E_{P,i}$.
By taking the direct sum,
we obtain 
a Hermitian metric $h_{t,P,0}:=\bigoplus_{i=1}^r h_{t,E_{P,i}}$
of $E_{|U_P}$.
Note that 
Theorem \ref{thm;15.7.20.1}
implies the following.
\begin{lem}
We have positive constants
$C_{P,0}'$ and $\epsilon_{P,0}'$
such that the following estimate holds
for any local sections $u_i$ and $u_j$ $(i\neq j)$ 
of $E_{P,i}$ and $E_{P,j}$:
\[
 \bigl|
 h_t(u_i,u_j)
 \bigr|
\leq
 C_{P,0}'\exp(-\epsilon_{P,0}'t)
 |u_i|_{h_t}\,|u_j|_{h_t}
\]
In particular,
we have a constant $K_{P}>1$
such that 
$K_P^{-1}h_{t,P,0}\leq h_{t|U_P}\leq K_Ph_{t,P,0}$
for any $t>1$.
\end{lem}

By varying $P\in X\setminus D(E,\theta)$
and by gluing $h_{t,P,0}$, 
we obtain a family of Hermitian metrics $h_{t,0}$  $(t>1)$
on $E_{|X\setminus D(E,\theta)}$.
We have the Chern connection $\nabla_{t,0}$
of $(E_{|X\setminus D(E,\theta)},h_{t,0})$.
Let $(t\theta)^{\dagger}_{h_{t,0}}$
denote the adjoint of $t\theta_{|X\setminus D(E,\theta)}$
with respect to $h_{t,0}$.
We set 
$\DD^1_{h_{t,0}}:=
 \nabla_{t,0}+
 t\theta+
 (t\theta)^{\dagger}_{h_{t,0}}$.

\begin{thm}[Proposition \ref{prop;15.7.12.12},
 Corollary \ref{cor;15.7.17.40}]
\label{thm;15.7.19.30}
Take any neighbourhood $N$ of $D(E,\theta)$.
Then, we have a constant $K>1$
such that
$K^{-1}h_{t,0|X\setminus N}
\leq
 h_{t|X\setminus N}
\leq
 Kh_{t,0|X\setminus N}$
for any $t>1$.
We also have positive constants
$C_1$ and $\epsilon_1$
such that the following holds on $X\setminus N$:
\[
 \bigl|
  R(h_{t,0})
 \bigr|_{h_{t},g_X}
\leq
 C_1\exp(-\epsilon_1t),
\quad\quad
 \bigl|
 \DD^1_{h_t}
-\DD^1_{h_{t,0}}
 \bigr|_{h_t,g_X}
\leq
 C_1\exp(-\epsilon_1t)
\]
The constants $C_1$ and $\epsilon_1$
may depend only on 
$(X,g_X)$, $N$ and $\Sigma(E,\theta)$.
\end{thm}

We note that we can obtain these estimates
in an elementary way which is standard in the study of 
the asymptotic behaviour of harmonic bundles
around the singularity,
pioneered by Simpson \cite{Simpson90},
and pursued further by the author 
\cite{mochi2,Mochizuki-wild}.
We also emphasize that 
we can obtain these kinds of estimates 
without the assumption that harmonic bundles
are given on a compact Riemann surface.
Indeed, we shall study harmonic bundles
given on discs in \S\ref{section;15.7.16.10},
which is clearly enough for the above estimates
on relatively compact regions.

Finally, we remark that the estimates
can be generalized in the case $\deg(E)\neq 0$,
i.e.,
in the context of Higgs bundles with Hermitian-Einstein metrics
on curves.
(See \S\ref{subsection;16.5.23.30}.)

\subsubsection{Hitchin WKB-problem}

Together with a rather standard argument of singular perturbations,
we can apply Theorem \ref{thm;15.7.19.30}
to the Hitchin WKB-problem in \cite{KNPS}.

We recall a notation in \cite{KNPS}.
Let $V$ be an $r$-dimensional complex vector space.
For Hermitian metrics $h_1,h_2$,
we can take a base $e_1,\ldots,e_r$ of $V$
which is orthogonal with respect to both $h_i$ $(i=1,2)$.
We have the real numbers $\kappa_j$ $(j=1,\ldots,r)$
determined by
$\kappa_j:=\log|e_j|_{h_2}-\log|e_j|_{h_1}$.
We impose 
$\kappa_1\geq \kappa_2\geq\cdots\geq \kappa_r$.
Then, we set
\[
 \vec{d}(h_1,h_2):=
 \bigl(
 \kappa_1,\ldots,\kappa_r
 \bigr)\in\real^r.
\]

Let us return to the study on
the family of harmonic bundles
$(E,\delbar_E,t\theta,h_t)$ $(t>0)$
for a stable Higgs bundle
$(E,\delbar_E,\theta)$ of rank $r$
with $\deg(E)=0$
on a compact Riemann surface $X$,
which is generically regular semisimple.
We take a universal covering
$\pi:Y\lrarr X\setminus D(E,\theta)$.
Then, 
we have the decomposition
of the Higgs bundle
$\pi^{\ast}(E,\delbar_E,\theta)
=\bigoplus_{i=1}^r (E_i,\delbar_{E_i},\phi_i\id_{E_i})$,
where $\phi_i$ are holomorphic $1$-forms.
We have $\rank E_i=1$,
and 
$\phi_i-\phi_j$ $(i\neq j)$ 
have no zeroes.

Let $[0,1]$ denote the closed interval
$\{0\leq s\leq 1\}$.
Let $\gamma:[0,1]\lrarr Y$ be a $C^{\infty}$-path.
We have the expressions
$\gamma^{\ast}(\phi_i)=a_i\,ds$
where $a_i$ are $C^{\infty}$-functions on $[0,1]$.
The path $\gamma$ is called non-critical
if $\Re a_i(s)\neq \Re a_j(s)$ $(i\neq j)$ for any $s$.
In that case, we may assume 
$\Re a_i(s)<\Re a_j(s)$ $(i<j)$.
We set
\[
 \alpha_i:=-\int_0^1\Re(a_i)ds.
\]

We have the families of Hermitian metrics
$h_{t,\gamma(\kappa)}$ $(t>0)$
on the fibers $E_{\gamma(\kappa)}$ $(\kappa=0,1)$,
induced by the harmonic metrics $h_t$.
Let $\Pi_{\gamma,t}:E_{\gamma(0)}\lrarr E_{\gamma(1)}$ 
denote the parallel transport of 
the flat connection $\DD^1_{h_t}$
along $\gamma$.
Let $\Pi_{\gamma,t}^{\ast}h_{t,\gamma(1)}$
denote the family of Hermitian metrics on $E_{\gamma(0)}$
induced by $h_{t,\gamma(1)}$
and $\Pi_{\gamma,t}$.
\begin{thm}[Theorem
 \ref{thm;15.7.20.100}]
If $\gamma$ is non-critical,
there exist positive constants $C_2$ and $\epsilon_2$
such that the following holds:
\[
\Bigl|
 \frac{1}{t}\vec{d}\bigl(h_{t,\gamma(0)},\Pi_{\gamma,t}^{\ast}h_{t,\gamma(1)}\bigr)
-(2\alpha_1,\ldots,2\alpha_r)
\Bigr|
\leq
 C_2\exp(-\epsilon_2t)
\]
The constants $C_2$ and $\epsilon_2$
may depend only on 
$X$, $\phi_1,\ldots,\phi_r$ and $\gamma$.
\end{thm}

The theorem was conjectured in \cite{KNPS},
and a different version of the problem
called the Riemann-Hilbert WKB-problem
was studied in detail.
Some cases were verified in \cite{Collier-Li}.
(See \cite{Collier-Li, KNPS} for the precise statements.)
We emphasize that
we can obtain this kind of estimate 
for more general families of harmonic bundles
given on complex curves
which are not necessarily compact.

\subsection{Convergence to the limiting configuration}

\subsubsection{Limiting configuration}

Let $(E,\delbar_E,\theta)$ 
be a stable Higgs bundle of rank $r$
with $\deg(E)=0$
on a compact Riemann surface $X$,
which is generically regular semisimple.
We have the family of harmonic bundles
$(E,\delbar_E,t\theta,h_t)$ $(t>0)$.
Let $\nabla_t$ denote the Chern connection of
$(E,\delbar_E,h_t)$.
By Theorem \ref{thm;15.7.20.1},
we can take a sub-sequence
$t_i\to\infty$ such that
the sequence of vector bundles
with a Hermitian metric and a unitary connection
$(E,\nabla_{t_i},h_{t_i})
 _{|X\setminus D(E,\theta)}$ 
converges in some sense
to a vector bundle with 
a Hermitian metric and a unitary flat connection
$(E_{\infty},\nabla_{\infty},h_{\infty})$
on $X\setminus D(E,\theta)$.
It is interesting to determine
the flat connection $\nabla_{\infty}$.
Following \cite{MSWW,MSWW2},
such a limit (or the associated parabolic bundle)
is called a limiting configuration of
the Higgs bundle $(E,\delbar_E,\theta)$
in this paper.
Note that it is not clear, in general,
whether $\nabla_{\infty}$ is independent
of the choice of a sub-sequence.
In this paper, we shall study the case $\rank E=2$
by assuming that $\theta$ is generically regular semisimple,
but without assuming that
the zeroes of 
$\det(\theta)-(\tr\theta)^2/4$ are simple.
See \S\ref{subsection;15.10.12.1}
for a remark on the case where the Higgs bundle
is not generically regular semisimple.

\begin{rem}
The higher rank case has not yet been studied in general.
Under some assumptions on the spectral curves,
it looks possible to generalize
the method in {\rm\cite{MSWW,MSWW2}} and 
our method in this paper.
See also {\rm\cite{Collier-Li}} for some Toda-like cases.
\hfill\qed
\end{rem}

\subsubsection{The case where 
the spectral curve is smooth irreducible
and simply ramified over $X$}
\label{subsection;15.7.20.10}

In \cite{MSWW,MSWW2},
Mazzeo, Swoboda, Weiss, and Witt
studied the case
under the assumption that
$\rank E=2$
and that the zeroes of 
$\det(\theta)-(\tr\theta)^2/4$
are simple,
i.e.,
the orders of the zeroes of
the quadratic differential 
$\det(\theta)-(\tr\theta)^2/4$
are at most $1$,
inspired by the work of 
Gaiotto, Moore and Neitzke \cite{GMN13}.
(We remark that in \cite{MSWW, MSWW2},
 the condition $\deg(E)=0$ is not imposed.)
Let us describe the limiting configuration
in this case.
It is enough to consider the case $\tr(\theta)=0$.
Under the assumption,
the spectral curve $\Xtilde:=\Sigma(E,\theta)$
is smooth and connected.
The natural projection
$\pi:\Xtilde\lrarr X$ is the ramified covering of degree $2$.
The ramification index is at most $2$.
The discriminant $D(E,\theta)$ is exactly 
the set of the points on which $\pi$ is ramified.
Let $\iota:\Xtilde\lrarr T^{\ast}X$ denote the inclusion.
We have the line bundle $\Ltilde$ on $\Xtilde$
such that 
$M_{E,\theta}\simeq \iota_{\ast}\Ltilde$.
We have
$E\simeq \pi_{\ast}\Ltilde$.

We have the unique non-trivial involution
$\rho:\Xtilde\lrarr \Xtilde$ over $X$,
i.e., $\pi\circ\rho=\pi$,
$\rho\circ\rho=\id_{\Xtilde}$
and $\rho\neq \id_{\Xtilde}$.
We have the line bundle
$\rho^{\ast}\Ltilde$ on $\Xtilde$.
We have a natural inclusion of 
$\nbigo_{\Xtilde}$-modules
$\pi^{\ast}E\lrarr
 \Ltilde\oplus \rho^{\ast}\Ltilde$,
and the cokernel 
is isomorphic to the structure sheaf of
$\Dtilde(E,\theta):=\pi^{-1}D(E,\theta)$.
We have
$\Ltilde\otimes
 \rho^{\ast}\Ltilde
\simeq
 \pi^{\ast}\det(E)\otimes
 \nbigo_{\Xtilde}(\Dtilde(E,\theta))$.

Because 
$\deg(\Ltilde)-\bigl|\Dtilde(E,\theta)\bigr|/2=0$,
we have a Hermitian metric 
$h^{\lim}_{\Ltilde}$ of 
$\Ltilde_{|\Xtilde\setminus\Dtilde(E,\theta)}$
with the following property.
\begin{itemize}
\item
 The Chern connection of 
$(\Ltilde_{|\Xtilde\setminus\Dtilde(E,\theta)},
 h_{\Ltilde}^{\lim})$
 is flat.
\item
 Let $\Ptilde\in \Dtilde(E,\theta)$.
 Let $e_{\Ptilde}$ be a local frame of 
 $\Ltilde$ around $\Ptilde$.
 Let $(\Utilde_{\Ptilde},w)$ be 
 a holomorphic coordinate system 
 around $\Ptilde$ with $w(\Ptilde)=0$.
 Then, 
 $\Bigl|
 \log\bigl(
 |w|h_{\Ltilde}^{\lim}(e_{\Ptilde},e_{\Ptilde})
 \bigr)\Bigr|$ is bounded on $\Utilde_{\Ptilde}\setminus\Ptilde$.
\item
We have
 $h^{\lim}_{\Ltilde}\otimes
 \rho^{\ast}h^{\lim}_{\Ltilde}
=\pi^{\ast}h_{\det(E)}$
 on 
 $\Xtilde\setminus\Dtilde(E,\theta)$.
\end{itemize}
Such $h^{\lim}_{\Ltilde}$ is uniquely determined.
We have the induced Hermitian metric
$h^{\lim}_{E,\theta}$
of $E_{|X\setminus D(E,\theta)}
=\pi_{\ast}(\Ltilde)_{|X\setminus D(E,\theta)}$.
We have the Chern connection
$\nabla^{\lim}_{E,\theta}$ of 
$(E_{|X\setminus D(E,\theta)},h^{\lim}_{E,\theta})$.
This is the limiting configuration
of $(E,\delbar_E,\theta)$ in this case.
Interestingly,
in \cite{MSWW, MSWW2},
it is proved that the family
$(E,\delbar_E,h_t,\theta)_{|X\setminus D(E,\theta)}$ $(t>0)$
is convergent to
$(E_{|X\setminus D(E,\theta)},\delbar_E,h^{\lim}_{E,\theta},\theta)$
in an appropriate sense.

\subsubsection{Limiting configuration in the general case}
\label{subsection;15.7.24.111}

It is natural to study the case
where $\theta$ is generically regular semisimple,
but the zeroes of $\det(\theta)-(\tr\theta)^2/4$ 
are not necessarily simple.
We may assume $\tr\theta=0$.
It is enough to consider the case
where the spectral curve of the Higgs bundle
$(E,\delbar_E,\theta)$
is reducible
to the two components,
i.e.,
we have holomorphic $1$-forms $\omega\neq 0$
such that
$\Sigma(E,\theta)=\Image(\omega)\cup\Image(-\omega)$.
Indeed, if the spectral curve is irreducible,
we have only to consider the pull back of
the Higgs bundle by a ramified covering map of degree $2$
given by the normalization of the spectral curve.
(See \S\ref{subsection;15.7.24.112}.
 See \S\ref{subsection;16.5.21.20} 
and \S\ref{subsection;16.5.21.30} for more details.)
Let $Z(\omega)$ denote the zero set of $\omega$,
which is equal to the discriminant $D(E,\theta)$.

We have two line bundles
$L_i$ $(i=1,2)$
with an inclusion $\iota_E:E\lrarr L_1\oplus L_2$
with the following property:
\begin{itemize}
\item
 Set $\theta_{L_1\oplus L_2}:=
 \omega\id_{L_1}\oplus 
(-\omega)\id_{L_2}$.
Then, we have
 $\theta_{L_1\oplus L_2}\circ\iota_E
=\iota_E\circ\theta$.
\item
 The restriction $\iota_{E|X\setminus Z(\omega)}$
 is an isomorphism.
\item
 The induced morphisms
 $E\lrarr L_i$ are surjective.
\end{itemize}

Set $d_i:=\deg(L_i)$.
We may assume that $d_1\leq d_2$.
For each $P\in Z(\omega)$,
let $m_P$ denote the order of the zero of $\omega$ at $P$.
Namely,
for a holomorphic coordinate system $(U_P,z)$ around $P$
with $z(P)=0$,
we have 
$\omega_{|U_P}=g_P\cdot z^{m_P}dz$ with $g_P(P)\neq 0$.
The support of the cokernel
$(L_1\oplus L_2)/E$ 
is contained in $Z(\omega)$.
For each $P\in Z(\omega)$,
let $\ell_P$ denote the length
of the stalk of $(L_1\oplus L_2)/E$ at $P$.
Note that we have
$L_1\otimes L_2=
 \det(E)\otimes
 \nbigo\bigl(\sum_{P\in Z(\omega)}\ell_PP\bigr)$.
In particular, we have
\[
 d_1+d_2-\sum_{P\in Z(\omega)}\ell_P=\deg(E)=0.
\]
Let $L_1'$ denote the kernel of
$E\lrarr L_2$.
Similarly,
let $L_2'$ denote the kernel of 
$E\lrarr L_1$.
Proper subbundles of $E$ preserved by $\theta$
are only $L_1'$ and $L_2'$.
Hence, the stability condition
for $(E,\theta)$ is equivalent to
$\deg(L_i')<0$ $(i=1,2)$.
It is equivalent to
$d_i=\deg(L_i)>0$ $(i=1,2)$.

To give the limiting configuration of $(E,\delbar_E,\theta)$,
we would like to give ``parabolic weights''
of $L_i$ $(i=1,2)$ at $Z(\omega)$.
In \S\ref{subsection;15.7.20.10},
all the parabolic weights are $1/2$.
In the general case,
it turns out that 
the parabolic weights are given by the following rule.

For each $P\in Z(\omega)$,
we consider the piecewise linear function
$\chi_P:\real_{\geq 0}\lrarr \real_{\leq 0}$
given by
\[
 \chi_P(a):=
 \left\{
\begin{array}{ll}
 (m_P+1)a-\ell_P/2
 & (a\leq \ell_P/2(m_P+1))
 \\
 0 &(a> \ell_P/2(m_P+1))
\end{array}
 \right.
\]
We set 
$\chi_{E,\theta}(a):=\sum_{P\in Z(\omega)}\chi_P(a)$.
Then, we have a unique 
$0\leq a_{E,\theta}<
\max_{P\in Z(\omega)}\bigl\{\ell_P/2(m_P+1)\bigr\}$
such that
$d_1+\chi_{E,\theta}(a_{E,\theta})=0$.

We determine 
the parabolic weights of $L_1$
at $P\in Z(\omega)$ as $-\chi_P(a_{E,\theta})$,
and 
the parabolic weights of $L_2$
at $P\in Z(\omega)$
as $\chi_P(a_{E,\theta})+\ell_P$.
We have associated singular Hermitian metrics on $L_i$,
i.e., we consider Hermitian metrics
$h^{\lim}_{L_i}$ $(i=1,2)$ 
on $L_{i|X\setminus Z(\omega)}$
satisfying the following conditions.
\begin{itemize}
\item
 The Chern connections of
 $(L_{i|X\setminus Z(\omega)},h^{\lim}_{L_i})$
 are flat.
\item
 Let $P\in Z(\omega)$. 
 Let $v_{P,i}$ be local frames of $L_i$ around $P$.
 Let $(U_P,z)$ be a holomorphic coordinate neighbourhood 
around $P$ with $z(P)=0$.
 Then, the functions
\[
 \log\bigl(
 |z|^{-2\chi_P(a_{E,\theta})}h^{\lim}_{L_1}(v_{P,1},v_{P,1})
 \bigr),
\quad
 \log\bigl(
 |z|^{2\chi_P(a_{E,\theta})+2\ell_P}
 h^{\lim}_{L_2}(v_{P,2},v_{P,2})
 \bigr)
\]
are bounded
on $U_P\setminus P$.
\item
 We have
 $h^{\lim}_{L_1}\otimes h^{\lim}_{L_2}
=h_{\det(E)}$.
\end{itemize}

We obtain the Hermitian metric
$h^{\lim}_{E,\theta}=h^{\lim}_{L_1}\oplus h^{\lim}_{L_2}$
of $E_{|X\setminus Z(\omega)}$,
and the Chern connection $\nabla^{\lim}_{E,\theta}$.

\begin{rem}
\label{rem;15.7.24.110}
We still have the ambiguity for $h^{\lim}_{L_i}$.
Namely, if $h^{\lim}_{L_i}$ $(i=1,2)$
satisfy the above conditions,
$\alpha h^{\lim}_{L_1}$
and $\alpha^{-1}h^{\lim}_{L_2}$
also satisfy the above conditions for any $\alpha>0$.
In {\rm\S\ref{subsection;15.7.20.10}},
we have an additional condition
$\rho^{\ast}h^{\lim}_{L_1}=h^{\lim}_{L_2}$
under the natural isomorphism
$\rho^{\ast}L_1\simeq L_2$,
with which the metrics are uniquely determined.
But, in general,
it seems that we do not have such an extra condition.

The metric $h^{\lim}_{E,\theta}$ is also characterized
as a harmonic metric of the polystable Higgs bundle
$(E,\delbar_{E},\theta)_{|X\setminus Z(\omega)}$
adapted to the parabolic structure
given as above.

Although we have the ambiguity of the metrics,
the connection
$\nabla^{\lim}_{E,\theta}$ is uniquely determined.
\hfill\qed
\end{rem}

\begin{rem}
Let us consider the case where 
$m_P=\ell_P=1$ for any $P\in Z(\omega)$,
and $d_1=d_2=|Z(\omega)|/2$.
Then, we have $a_{E,\theta}=0$,
and hence
$-\chi_{P}(a_{E,\theta})=\chi_P(a_{E,\theta})+\ell_P=1/2$.
These are the parabolic weights appeared
in {\rm\S\ref{subsection;15.7.20.10}}.
\hfill\qed
\end{rem}

\begin{rem}
It might be instructive
to mention that 
if we have $\deg(L_1)=\deg(L_2)$
then 
we have
$a_{E,\theta}=0$,
and $-\chi_P(a_{E,\theta})=\chi_P(a_{E,\theta})+\ell_P=\ell_P/2$.
\hfill\qed
\end{rem}

\subsubsection{Convergence to the limiting configuration}

Suppose that $(E,\delbar_E,\theta)$
is stable of degree $0$ with $\rank E=2$
and that $\theta$ is generically regular semisimple.
For simplicity, we assume $\tr\theta=0$.
Let $h_t$ $(t>0)$ denote the harmonic metric
of the Higgs bundle $(E,\delbar_E,t\theta)$
with $\det(h_t)=h_{\det(E)}$.

For any $\alpha>0$,
let $\Psi_{\alpha}$ be the automorphism
of $L_1\oplus L_2$ given by 
$\Psi_{\alpha}=
 \alpha\id_{L_1}
\oplus 
 \alpha^{-1}\id_{L_2}$.
Let $\Psi_{\alpha}^{\ast}h_t$
be the metric of $E_{|X\setminus D(E,\theta)}$
given by
$h_t(\Psi_{\alpha}s_1,\Psi_{\alpha}s_2)$
for local sections $s_j$ $(j=1,2)$
of $E_{|X\setminus D(E,\theta)}$.
We take any point $Q\in X\setminus D(E,\theta)$
and a frame $v_Q$ of $L_{1|Q}$,
and we put
\[
\gamma(t,Q):=\left(
 \frac{h^{\lim}_{L_1}(v_Q,v_Q)}{h_t(v_Q,v_Q)}
 \right)^{1/2}.
\]
The following theorem is our second main result
in this paper.

\begin{thm}[The degree $0$ case of Theorem 
\ref{thm;15.7.15.10}]
\label{thm;15.7.20.31}
When $t$ goes to $\infty$,
the sequence
$\Psi_{\gamma(t,Q)}^{\ast}h_t$ 
converges to 
$h^{\lim}_{E,\theta}$
in the $C^{\infty}$-sense
on any compact subset 
in $X\setminus Z(\omega)$.
In particular,
the sequence of the connections $\nabla_{t}$
converges to $\nabla^{\lim}_{E,\theta}$.
(See {\rm\S\ref{subsection;15.7.24.111}}
for $h^{\lim}_{E,\theta}$ and $\nabla^{\lim}_{E,\theta}$.)
\end{thm}

For the proof of Theorem \ref{thm;15.7.15.10},
we need to study the global property of $(E,\delbar_E,\theta)$,
in contrast that 
we can obtain the estimates in \S\ref{subsection;15.7.20.30}
locally at any point of $X\setminus D(E,\theta)$.
A key is the construction of a family of 
Hermitian metrics $h^0_{t}$ of $E$
such that
the $L^p$-norms of
$R(h^0_t)
+\bigl[t\theta,(t\theta)^{\dagger}_{h^0_t}\bigr]$
are bounded,
for which we naturally encounter the above parabolic weights
for the limiting configuration.
We have the family of the self-adjoint endomorphisms
$k_t$ of $(E,h^0_t)$
such that
$h_t(u_1,u_2)=h^0_t(k_tu_1,u_2)$
for any local sections $u_i$ $(i=1,2)$.
By applying a variant of the arguments 
used in \cite{mochi5}
with the tools given in \cite{Simpson88},
we shall observe that 
the sequence $\kappa_i k_{t_i}$ converges 
in some sense
for an appropriate sequence $\kappa_i>0$.
Then, the claim of the theorem follows.

We note that for the construction of 
the family of metrics $h^0_t$,
we use the general theory of wild harmonic bundles.
We need a family of harmonic metrics given on a neighbourhood
of the discriminant locus,
for which we apply the Kobayashi-Hitchin correspondence
for wild harmonic bundles on curves
given in \cite{biquard-boalch}
(see \cite{Simpson90} for the tame case).
On the basis of the general results on the asymptotic behaviour
of wild harmonic bundles studied in 
\cite{Simpson90} and \cite{Mochizuki-wild},
we can deduce rather detailed properties of 
the family of the harmonic metrics
as in Proposition \ref{prop;15.7.13.40}
and Proposition \ref{prop;15.7.13.110}.

\vspace{.1in}

At this moment,
it is not clear to the author
whether we could directly use the argument in \cite{MSWW}
to prove Theorem \ref{thm;15.7.20.31} or its variant
in our setting.
But, we should note that we do not study 
the order of the convergence in Theorem \ref{thm;15.7.20.31}.
In contrast, 
when all the zeroes of $\det\theta-(\tr\theta)^2/4$ are simple,
the method in \cite{MSWW} is strong enough
to give the order of the convergence.

\begin{rem}
Let $(E,\delbar_E,\theta)$ be a stable Higgs bundle of rank $2$
on $X$
such that $\theta$ is generically regular semisimple,
but that  $\deg(E)$ is not necessarily $0$.
We fix a Hermitian metric $h_{\det(E)}$ of $\det(E)$.
For each $t>0$,
we have the Hermitian-Einstein metric $h^{HE}_t$
of $(E,\delbar_E,\theta)$ such that 
$\det(h^{HE}_t)=h_{\det E}$,
according to 
{\rm\cite{Hitchin-self-duality}} and {\rm\cite{Simpson88}}.
Here, the Hermitian-Einstein condition means
the trace-free part of
$R(h^{HE}_t)+[t\theta,(t\theta)^{\dagger}_{h^{HE}_t}]$
is $0$.
We can study the behaviour of $h^{HE}_t$ $(t\to\infty)$
by using Theorem {\rm\ref{thm;15.7.20.31}}.
(See Theorem {\rm\ref{thm;15.7.15.10}}.)

Suppose that $\deg(E)=2m$ for an integer $m$.
We take a line bundle $L$ with $\deg(L)=-m$.
We can easily reduce the study on the behaviour of $h^{HE}_t$
to the study on the behaviour of harmonic metrics
for $(E,\delbar_E,\theta)\otimes L$.
Suppose that $\deg(E)$ is odd.
We take any (unramified) covering map $\varphi:X'\lrarr X$ of degree $2$.
Then, 
the degree of $\varphi^{\ast}(E,\delbar_E,\theta)$
is even.
Hence, the study can be reduced
to the degree $0$ case.
(See also {\rm\S\ref{subsection;16.5.19.22}}.)
\hfill\qed
\end{rem}

\subsubsection{Symmetric case and irreducible case}
\label{subsection;15.7.24.112}

Let us explain that we can obtain a stronger result
if the Higgs bundle is equipped with an extra symmetry.
Suppose that $X$ is equipped with a non-trivial involution $\rho$,
i.e., $\rho:X\lrarr X$ is a holomorphic automorphism
such that $\rho\circ\rho=\id_X$ and $\rho\neq\id_X$.
Let $(E,\delbar_E,\theta)$ be a stable Higgs bundle of degree $0$
on $X$ with $\tr\theta=0$
such that the spectral curve is reducible,
i.e.,
$\Sigma(E,\theta)=\Image(\omega)\cup\Image(-\omega)$
for a holomorphic one form $\omega\neq 0$ on $X$.
We impose the following conditions.
\begin{itemize}
\item
 We have $\rho^{\ast}\omega=-\omega$.
\item
 The Higgs bundle $(E,\delbar_E,\theta)$
 is equivariant with respect to the action of 
 $\{1,\rho\}$ on $X$.
\end{itemize}
We have the induced isomorphism
$\rho^{\ast}\det(E)\simeq \det(E)$.
We naturally have 
$\rho^{\ast}h_{\det(E)}=h_{\det(E)}$.

Let $L_i$ $(i=1,2)$
be as in \S\ref{subsection;15.7.24.111}.
We naturally have the isomorphisms
$\rho^{\ast}L_1\simeq L_2$
and 
$\rho^{\ast}L_2\simeq L_1$,
which are compatible with 
the isomorphism
$\rho^{\ast}E\simeq E$.
Because $\deg(L_1)=\deg(L_2)$,
we have $a_{E,\theta}=0$
and $-\chi_P(a_{E,\theta})=\chi_P(a_{E,\theta})+\ell_P
=\ell_P/2$ for any $P\in Z(\omega)$.
We can uniquely determine the metrics
$h^{\lim}_{L_j}$
by imposing the extra condition
$\rho^{\ast}h^{\lim}_{L_1}=h^{\lim}_{L_2}$.

\begin{thm}[The degree $0$ case of Theorem 
\ref{thm;15.7.24.143}]
\label{thm;15.7.24.130}
For $t>0$, let $h_t$ be the harmonic metrics
of the Higgs bundles $(E,\delbar_E,t\theta)$
such that $\det(h_t)=h_{\det(E)}$.
Then, the sequence $h_t$ is convergent
to $h^{\lim}_{E,\theta}=h^{\lim}_{L_1}\oplus h^{\lim}_{L_2}$
on any compact subset in $X\setminus Z(\omega)$.
\end{thm}

Let us remark that we can apply 
Theorem \ref{thm;15.7.24.130}
to the case where the spectral curve is irreducible
and the Higgs field is generically regular semisimple.
Let $X'$ be a compact connected Riemann surface.
Let $(E',\delbar_{E'},\theta')$ be a stable Higgs bundle
of rank $2$ with $\deg(E')=0$
on $X'$ such that
(i) the spectral curve
$\Sigma(E',\theta')$ is irreducible,
(ii) $\tr\theta'=0$,
(iii) $\theta'$ is generically regular semisimple.
We take the normalization
$\kappa:X\lrarr \Sigma(E',\theta')$.
Let $p:X\lrarr X'$ be the morphism
obtained as the composite of
$\kappa$ and the projection
$\Sigma(E',\theta')\lrarr X'$.
We set
$(E,\delbar_E,\theta)=p^{\ast}(E',\delbar_{E'},\theta')$.
We have the involution on $\Sigma(E',\theta')$
induced by the multiplication of $-1$
on the cotangent bundle $T^{\ast}X'$.
It induces an involution $\rho$ on $X$.
We fix a Hermitian metric $h_{\det(E')}$ 
of $\det(E')$.

If $(E,\delbar_E,\theta)$ is polystable,
we can easily observe that 
the harmonic metrics $h_t'$
for the Higgs bundles $(E',\delbar_{E'},t\theta')$
are independent of $t$.
(See \S\ref{subsection;16.5.21.20}.)
We formally set $h^{\lim}_{E',\theta'}:=h_t'$.

If $(E,\delbar_E,\theta)$ is stable,
the above symmetry conditions for 
$(E,\delbar_E,\theta)$ and $\rho$ are satisfied.
The pull back $h_{\det(E)}=p^{\ast}h_{\det(E')}$
satisfies the condition
$\rho^{\ast}h_{\det(E)}=h_{\det(E)}$.
We have the metrics $h_{L_i}^{\lim}$ 
and 
$h^{\lim}_{E,\theta}
=h^{\lim}_{L_1}\oplus h^{\lim}_{L_2}$
as above.
Because we have
$\rho^{\ast}h^{\lim}_{E,\theta}=h^{\lim}_{E,\theta}$,
we have a unique Hermitian metric
$h^{\lim}_{E',\theta'}$ on $X'\setminus D(E',\theta')$
such that $p^{\ast}h^{\lim}_{E',\theta'}=h^{\lim}_{E,\theta}$.
We also have a characterization of
$h^{\lim}_{E',\theta'}$
as a harmonic metric for 
the stable filtered Higgs bundles on $(X',D(E',\theta'))$
induced by the limiting configuration of
$(E,\delbar_E,\theta)$.
(See \S\ref{subsection;16.5.21.20}.)
We have the following direct corollary.

\begin{cor}[The degree $0$ case of Corollary \ref{cor;16.5.19.10}]
For $t>0$, let $h'_t$ be the harmonic metrics
of the Higgs bundles $(E',\delbar_{E'},t\theta')$
such that $\det(h'_t)=h_{\det(E')}$.
Then, the sequence $h'_t$ is convergent
to $h^{\lim}_{E',\theta'}$
on any compact subset in $X'\setminus D(E',\theta')$.
\hfill\qed
\end{cor}

\subsubsection{Complement}
\label{subsection;15.10.12.1}

Let $(E,\delbar_E,\theta,h)$ be a Higgs bundle of rank $2$.
Suppose that $\theta$ is nilpotent,
i.e., $\Sigma(E,\theta)$ is the $0$-section in $T^{\ast}X$.
Then, we can easily observe that
$\lim_{t\to\infty}(E,\delbar_E,t\theta,h_t)$
is convergent to a polarized complex variation of Hodge structure.
We sketch it.
We have a polystable Higgs bundle
$(E_{\infty},\delbar_{\infty},\theta_{\infty})$
obtained as the limit
$\lim_{t\to\infty}(E,\delbar_E,t\theta)$
in the coarse moduli space of semistable Higgs bundles 
of degree $0$ on $X$.
Note that the norms
$\bigl|R(\nabla_{h_t})\bigr|_{h_t,g_X}=
 \bigl|t\theta\bigr|^2_{h_t,g_X}$ $(t>0)$
are uniformly bounded,
where $g_X$ is a K\"ahler metric of $X$.
(See Proposition \ref{prop;15.7.17.10}, for example.)
Hence, we obtain a harmonic metric $h_{\infty}$
of $(E_{\infty},\delbar_{\infty},\theta_{\infty})$
as the limit of a convergent subsequence $h_{t_i}$ $(t_i\to\infty)$.
If $(E_{\infty},\delbar_{\infty},\theta_{\infty})$
is not stable,
$\theta_{\infty}$ is trivial,
$(E_{\infty},\delbar_{\infty})$ 
is a direct sum of line bundles,
and $h_{\infty}$ is a flat metric.
If $(E_{\infty},\delbar_{\infty},\theta_{\infty})$
is stable,
then as in the case of the limit for $t\to 0$
studied by Simpson,
we can observe that
$(E_{\infty},\delbar_{\infty},\theta_{\infty})$
is a Hodge bundle,
and $h_{\infty}$ is equivariant
with respect to the natural grading $S^1$-action.
Hence, 
$(E_{\infty},\delbar_{\infty},\theta_{\infty},h_{\infty})$
comes from a polarized complex variation of Hodge structure.

In view of the moduli theoretic picture 
\cite{Hitchin-self-duality,s4},
we can see it as follows.
The energies of the family of
$(E,\delbar_E,t\theta,h_t)$ $(t>0)$ are bounded.
So, the family is relatively compact
in the moduli space of harmonic bundles.
Hence, when $t$ goes to $\infty$,
$(E,\delbar_E,\theta,h_t)$ goes
to a fixed point in the moduli space
induced by the natural $\cnum^{\ast}$-action
on the moduli space of Higgs bundles.
It means that $(E,\delbar_E,\theta,h_t)$ is 
convergent to a polarized complex variation of Hodge structure.

If $(E,\delbar_E,\theta)$ is not generically regular semisimple
but $\rank E=2$,
then the study of 
$\lim_{t\to\infty}(E,\delbar_E,t\theta,h_t)$
is easily reduced to the above case.

\subsection{Acknowledgement}

The author thanks the referee for his thorough reading
and valuable remarks.
This study has grown out of my effort to understand
(part of)
the intriguing works
\cite{Collier-Li}, 
\cite{GMN13},
\cite{KNPS}
and \cite{MSWW, MSWW2}.
I thank Carlos Simpson 
for inspiring discussions on many occasions.
Indeed, 
I am motivated to this study
by his question several years ago.
I owe a lot to his philosophy and his techniques.
I also thank him for many valuable comments
to improve this manuscript,
and for his patience in reading the preliminary versions of
this manuscript.
I hope that this study would be useful
for the project in \cite{KNPS}.
I learned the attractive study
\cite{MSWW, MSWW2}
from the stimulating talk given by Frederik Witt 
in the conference 
``The Geometry, Topology and 
Physics of Moduli Spaces of Higgs Bundles''
held at the Institute for Mathematical Sciences,
National University of Singapore in 2014,
for which the visit was supported by the Institute.
I thank Claude Sabbah for his kindness
and for many discussions.
I am grateful to
Yoshifumi Tsuchimoto and Akira Ishii
for their constant encouragement.
I thank Szilard Szabo for his interest to this study
and some discussions.

This work was partially supported by 
the Grant-in-Aid for Scientific Research (S)
(No. 24224001) and 
the Grant-in-Aid for Scientific Research (C) 
(No. 15K04843), 
Japan Society for the Promotion of Science.

\section{Asymptotic decoupling}
\label{section;15.7.16.10}

\subsection{Simpson's main estimate for harmonic bundles on discs}

\label{subsection;15.7.17.20}

Let $(E,\delbar_E,\theta,h)$ be a harmonic bundle
on a complex curve $X$.
We have the spectral curve 
$\Sigma(E,\theta)$ in the cotangent bundle of $X$.
We have an estimate of the norm of $\theta$,
depending only on $\Sigma(E,\theta)$
(Proposition \ref{prop;15.7.1.10}).
If $\Sigma(E,\theta)$ is decomposed 
into a disjoint union $\Sigma_1\sqcup\Sigma_2$,
then we have the corresponding
decomposition of the Higgs bundle
$(E,\delbar_E,\theta)
=(E_1,\delbar_{E_1},\theta_1)
 \oplus
 (E_2,\delbar_{E_2,}\theta_2)$.
The Hermitian product 
of sections of $E_1$ and $E_2$ given by $h$
should be small.
We have such estimates
depending only on $\Sigma(E,\theta)$
(Corollary \ref{cor;15.8.18.22}).
Because we use the arguments essentially given
in \cite{Simpson90},
such estimates are called Simpson's main estimate.
It is enough to consider the case where $X$ is a disc.

\subsubsection{Estimate of the sup norm 
 in terms of the eigenvalues}
\label{subsection;15.7.23.1}

For $R>0$,
we set $\Delta(R):=
 \bigl\{z\in\cnum\,\big|\,|z|<R\bigr\}$.
We consider a harmonic bundle $(E,\delbar_E,\theta,h)$
of rank $r$ on $\Delta(R)$.
We have the description
$\theta=f\,dz$,
where $f$ is a holomorphic endomorphism of $E$.
Fix $M>0$,
and suppose the following.
\begin{itemize}
\item
 For any $P\in \Delta(R)$,
 the eigenvalues $\gamma$ of
 $f_{|P}$ satisfy $|\gamma|<M$.
\end{itemize}

We recall the following proposition,
for which a proof was given 
in \cite{Mochizuki-KH1},
for example.
We include an outline of the proof 
for the convenience of the readers.

\begin{prop}
\label{prop;15.7.17.10}
Fix $0<R_1<R$.
Then, we have $C_1,C_2>0$ depending only on
$r, R_1,R$ such that
\[
 |f|_h\leq C_1M+C_2
\]
holds on $\Delta(R_1)$.
\end{prop}
\pf
Let $f^{\dagger}_h$ denote the adjoint of $f$
with respect to $h$.
As in \cite{Simpson90},
we have the following inequality
on $\Delta(R)$
(see also Lemma \ref{lem;15.8.2.20} below):
\[
 -\del_z\del_{\zbar}
 \log|f|^2_h
\leq
 -\frac{\bigl|[f,f_h^{\dagger}]\bigr|_h^2}{|f|_h^2}
\]
For any $P\in \Delta(R)$,
we set $g(P):=\sum_{i=1}^r |\alpha_i|^2$,
where $\alpha_1,\ldots,\alpha_r$
are the eigenvalues of $f_{|P}$.
There exists a constant $C_3>0$
depending only on $r$
such that
\[
 \bigl|
 [f,f^{\dagger}_h]_{|P}
 \bigr|_h
\geq
 C_3\bigl(|f_{|P}|_h^2-g(P)\bigr).
\]
We obtain
\[
 -\del_z\del_{\zbar}
 \log|f|_h^2(P)
\leq 
 -C_3^2\frac{\bigl||f_{|P}|_h^2-g(P)\bigr|^2}{|f_{|P}|_h^2}
\]
Note that if $|f_{|P}|_h^2\geq 2g(P)$ for some $P$,
then we have
\[
 -\frac{\bigl||f_{|P}|_h^2-g(P)\bigr|^2}{|f_{|P}|_h^2}
\leq
 -\frac{1}{4}|f_{|P}|_h^2.
\]
For any positive number $B$,
we have
\[
 -\del_z\del_{\zbar}
 \log\frac{B}{(R^2-|z|^2)^2}
=-\frac{2R^2}{B}\frac{B}{(R^2-|z|^2)^2}.
\]
Take $B$ satisfying
$B\geq 8R^2/C_3^2$
and $B\geq 2R^4rM^2$.
Then, we use an idea in the proof of 
a lemma of Ahlfors \cite{a}.
Set $\nbigz:=\bigl\{
 P\in\Delta(R)\,\big|\,
 |f_{|P}|_h^2>B(R^2-|z(P)|^2)^{-2}
 \bigr\}$.
Suppose that $\nbigz\neq\emptyset$,
and we shall derive a contradiction.
We have
$|f_{|P}|_h^2>BR^{-4}
\geq 2rM^2>2g(P)$
for any $P\in\nbigz$.
Hence, the following inequality holds on $\nbigz$:
\[
 -\del_z\del_{\zbar}\Bigl(
 \log|f|_h^2-\log\bigl(B(R^2-|z|^2)^{-2}\bigr)
 \Bigr)
\leq
 -\frac{C_3^2}{4}
 \Bigl(|f|_h^2
-B(R^2-|z|^2)^{-2}
\Bigr)
\leq 0
\]
Note that $\nbigz$ is relatively compact
in $\Delta(R)$.
So, we have
$|f|_h^2=B(R^2-|z|^2)^{-2}$
on the boundary of $\nbigz$,
and hence we obtain 
$|f|_h^2\leq B(R^2-|z|^2)^{-2}$
on $\nbigz$.
But, it contradicts with 
the choice of $\nbigz$.
So, we obtain $\nbigz=\emptyset$.
Namely, we have
$|f|_h^2\leq B(R^2-|z|^2)^{-2}$
on $\Delta(R)$.
Then, we obtain the claim of the proposition.
\hfill\qed

\subsubsection{Asymptotic orthogonality}

We continue to use the notation
in \S\ref{subsection;15.7.23.1}.
Suppose that we have
a finite subset $S\subset \cnum$
and a decomposition
$(E,\theta)=
 \bigoplus_{\alpha\in S}(E_{\alpha},\theta_{\alpha})$.
We have the expression $\theta_{\alpha}=f_{\alpha}dz$
for each $\alpha$,
where $f_{\alpha}$ is a holomorphic endomorphism of
$E_{\alpha}$.

\begin{assumption}
\label{assumption;15.7.17.100}
Fix $C_{10}>1$,
and we impose the following conditions.
\begin{itemize}
\item
We set 
$d:=\min\bigl\{|\alpha-\beta|\,\big|\,\alpha,\beta\in S,\alpha\neq\beta
 \bigr\}$.
Then, 
we have
$d\geq 1$ and $M\leq C_{10}d$.
Here, $M$ is the constant in the beginning of
{\rm \S\ref{subsection;15.7.23.1}}.
\item
For any $P\in\Delta(R)$,
the eigenvalues $\gamma$ of $f_{\alpha|P}$
satisfy $|\gamma-\alpha|\leq d/100$.
\end{itemize}
\end{assumption}

Let $\pi_{\alpha}$ be the projection of $E$ onto $E_{\alpha}$
with respect to the above decomposition.
Let $\pi_{\alpha}'$ denote the orthogonal projection
of $E$ onto $E_{\alpha}$.
We set $\rho_{\alpha}:=\pi_{\alpha}-\pi_{\alpha}'$.
The following proposition is 
a variant of the estimates given 
in \cite{Simpson90} and \cite{mochi2, Mochizuki-wild}
for harmonic bundles given on punctured discs
or the products of punctured discs,
but it is more useful for the purpose in this paper.

\begin{prop}
\label{prop;15.7.12.10}
Let $R_1$ be as in Proposition {\rm\ref{prop;15.7.17.10}}.
Fix $0<R_2<R_1$.
There exist positive constants 
$\epsilon_0$
and $C_{11}$
depending only on
$R$, $R_1$, $R_2$, $r$ and $C_{10}$,
such that 
$|\rho_{\alpha}|_h\leq 
 C_{11}\exp(-\epsilon_0d)$
on $\Delta(R_2)$.
\end{prop}
\pf
We recall a general inequality for holomorphic sections
of $\End(E)$
which commutes with $\theta$,
for the convenience of the readers.
\begin{lem}
\label{lem;15.8.2.20}
Let $s$ be any holomorphic section of $\End(E)$
such that $[\theta,s]=0$.
Then, we have the following inequality:
\begin{equation}
\label{eq;15.8.2.21}
 -\del_z\del_{\zbar}
 \log|s|_h^2
\leq
-\frac{\bigl|[f_h^{\dagger},s]\bigr|_h^2}{|s|_h^2}
\end{equation}
\end{lem}
\pf
We have the equality
$-\del_z\del_{\zbar}|s|_h^2
=-|\del_{z,h}s|_h^2
-h\bigl(s,\del_{\zbar}\del_{z,h}s\bigr)$.
Hence, we obtain the following:
\begin{multline}
\label{eq;15.8.2.10}
 -\del_z\del_{\zbar}\log|s|_h^2
=-\frac{\del_z\del_{\zbar}|s|_h^2}{|s|_h^2}
+\frac{\del_z|s|_h^2}{|s|_h^2}
\frac{\del_{\zbar}|s|_h^2}{|s|_h^2}
=-\frac{h(s,\del_{\zbar}\del_{z,h}s)}{|s|_h^2}
-\frac{\bigl|\del_{z,h}s\bigr|_h^2}{|s|_h^2}
+\frac{h(\del_{z,h}s,s)}{|s|_h^2}
\frac{h(z,\del_{z,h}s)}{|s|_h^2}
 \\
\leq
-\frac{h(s,\del_{\zbar}\del_{z,h}s)}{|s|_h^2}
=
-\frac{h(s,(\del_{\zbar}\del_{z,h}-\del_{z,h}\del_{\zbar})s)}{|s|_h^2}.
\end{multline}
By the Hitchin equation,
we have
$R(h)+[\theta,\theta^{\dagger}_h]
=R(h)+[f,f_h^{\dagger}]\,dz\,d\zbar=0$,
i.e.,
$R(h)=[f,f_h^{\dagger}]d\zbar\,dz$.
We obtain the following from 
(\ref{eq;15.8.2.10}):
\[
 -\del_{z}\del_{\zbar}
 \log|s|_h^2
\leq
 -\frac{h\bigl(s,
 \bigl[[f,f_h^{\dagger}],s\bigr]
 \bigr)}{|s|_h^2}
\]
Because $[f,s]=0$,
we obtain 
\[
 -\del_z\del_{\zbar}
 \log|s|_h^2
\leq
 -\frac{h\bigl(s,\bigl[f,[f_h^{\dagger},s]\bigr]\bigr)}{|s|_h^2}
=-\frac{h\bigl([f_h^{\dagger},s],[f_h^{\dagger},s]\bigr)}{|s|_h^2}
=-\frac{\bigl|[f_h^{\dagger},s]\bigr|_h^2}{|s|_h^2}
\]
Thus, we obtain Lemma \ref{lem;15.8.2.20}.
\hfill\qed

\vspace{.1in}

Let $r_{\alpha}:=\rank E_{\alpha}$.
Because $\del_{\zbar}\del_z\log r_{\alpha}=0$,
we obtain the following inequality on $\Delta(R)$
from (\ref{eq;15.8.2.21}):
\[
 -\del_z\del_{\zbar}
 \log\bigl(|\pi_{\alpha}|_h^2/r_{\alpha}\bigr)
\leq
 -\frac{\bigl|[f_h^{\dagger},\pi_{\alpha}]\bigr|_h^2}{|\pi_{\alpha}|_h^2}
\]

According to Proposition \ref{prop;15.7.17.10}
and Lemma \ref{lem;15.7.1.1} below,
there exists a positive constant $C_{12}$ 
depending only on
$R$, $R_1$, $r$ and $C_{10}$
such that
we have
$\bigl|\pi_{\alpha}\bigr|_h\leq C_{12}$
and 
$|\rho_{\alpha}|_h\leq C_{12}$
on $\Delta(R_1)$.
We set 
$k_{\alpha}:=\log(|\pi_{\alpha}|_h^2/r_{\alpha})$.
Note that 
$|\pi_{\alpha}|_h^2=|\pi_{\alpha}'|_h^2+|\rho_{\alpha}|_h^2$
and 
$|\pi_{\alpha}'|_h^2=r_{\alpha}$.
Hence, we have
$k_{\alpha}=\log(1+|\rho_{\alpha}|^2/r_{\alpha})$.
There exists a positive constant $C_{13}$
depending only on 
$R$, $R_1$, $r$ and $C_{10}$
such that 
$C_{13}^{-1}|\rho_{\alpha}|_h^2
\leq k_{\alpha}
\leq
 C_{13}|\rho_{\alpha}|_h^2$
on $\Delta(R_1)$.

According to Lemma \ref{lem;15.7.1.2} below,
there exists a constant $\epsilon_1>0$
depending only on
$R,R_1,r,C_{10}$
such that 
we have
$ \bigl|
 \bigl[f_h^{\dagger},\pi_{\alpha}\bigr] 
 \bigr|_h
\geq \epsilon_1d|\rho_{\alpha}|_h$
on $\Delta(R_1)$.
Hence, we have
a constant $\epsilon_2>0$
depending only on
$R,R_1,r,C_{10}$
such that the following holds
on $\Delta(R_1)$:
\[
 -\del_z\del_{\zbar}
 k_{\alpha}
\leq
 -\epsilon_2d^2k_{\alpha}
\]

For any positive number $\epsilon_3>0$,
we have the following on $\Delta(R_1)$:
\[
 -\del_z\del_{\zbar}
 \exp\bigl(
 \epsilon_3d|z|^2
 \bigr)
=-(d\epsilon_3+d^2|z|^2\epsilon_3^2)
 \exp(\epsilon_3d|z|^2)
\geq
 -(d\epsilon_3+d^2\epsilon_3^2R_1^2)
 \exp\bigl(\epsilon_3d|z|^2\bigr)
\]
Take $\epsilon_3>0$
such that
$\epsilon_3\leq (\epsilon_2R_1^{-1})/2$
and 
$\epsilon_3\leq d\epsilon_2/2$.
Then, we have
\[
 -\del_z\del_{\zbar}
 \exp\bigl(
 \epsilon_3d|z|^2
 \bigr)
\geq
 -\epsilon_2d^2
  \exp\bigl(
 \epsilon_3d|z|^2
 \bigr).
\]
We take $C_{14}>0$
depending only on 
$R,R_1,r,C_{10}$
such that 
$k_{\alpha}(P)<C_{14}$
for $|z(P)|=R_1$.
We set 
\[
 \nbigz:=\Bigl\{
 P\in\Delta(R_1)\,\Big|\,
 k_{\alpha}(P)>
 C_{14}
 \exp\bigl(
  \epsilon_3d|z|^2
-\epsilon_3dR_1^2
 \bigr)(P)
 \Bigr\}.
\]
Suppose that $\nbigz$ is non-empty.
By the choice of $C_{14}$,
$\nbigz$ is relatively compact
in $\Delta(R_1)$.
So, we have
$k_{\alpha}=
 C_{14}
 \exp\bigl(
  \epsilon_3d|z|^2
-\epsilon_3dR_1^2
 \bigr)$
on the boundary of $\nbigz$.
We also have the following inequality
on $\nbigz$:
\[
 -\del_z\del_{\zbar}\Bigl(
 k_{\alpha}-
 C_{14}
 \exp\bigl(
  \epsilon_3d|z|^2
-\epsilon_3dR_1^2
 \bigr)
 \Bigr)
\leq
 -\epsilon_2d^2\Bigl(
  k_{\alpha}-
 C_{14}
 \exp\bigl(
  \epsilon_3d|z|^2
-\epsilon_3dR_1^2
 \bigr)
 \Bigr)
\leq 0
\]
We obtain
$k_{\alpha}
\leq
 C_{14}
 \exp\bigl(
  \epsilon_3d|z|^2
-\epsilon_3dR_1^2
 \bigr)$
on $\nbigz$,
which contradicts with the construction of $\nbigz$.
Hence, we obtain that $\nbigz=\emptyset$,
i.e.,
$k_{\alpha}
\leq
 C_{14}
 \exp\bigl(
 \epsilon_3d|z|^2
-\epsilon_3dR_1^2)
=C_{14}\exp\bigl(
 -\epsilon_3(R_1^2-|z|^2)d
 \bigr)$
holds on $\Delta(R_1)$.
We obtain the following on
$\Delta(R_2)$:
\[
 k_{\alpha}
\leq C_{14}
 \exp\bigl(-\epsilon_3(R_1^2-R_2^2)d\bigr)
\]
Then, 
by setting $C_{11}:=C_{14}>0$
and $\epsilon_0:=\epsilon_3(R_1^2-R_2^2)>0$,
we obtain 
the claim of Proposition \ref{prop;15.7.12.10}.
\hfill\qed

\vspace{.1in}
For any endomorphism $g$ of $E$,
let $g_{h}^{\dagger}$
denote the adjoint of $g$
with respect to $h$.
We also denote it by $g^{\dagger}$
if there is no risk of confusion.

\begin{cor}
\label{cor;15.7.22.10}
There exist positive constants
$C_{20}$ and $\epsilon_{20}$,
depending only on 
$R$, $R_1$, $R_2$, $r$ and $C_{10}$ 
such that 
$\bigl|
 [f,(\pi_{\alpha})^{\dagger}_h]
 \bigr|_h
=\bigl|
 [f_h^{\dagger},\pi_{\alpha}]
 \bigr|_h
\leq
 C_{20}\exp(-\epsilon_{20} d)$
on $\Delta(R_2)$.
\end{cor}
\pf
We have
$|\rho_{\alpha}|_h=|\pi_{\alpha}-\pi_{\alpha}'|_h
\leq C_{11}\exp(-\epsilon_{0}d)$.
We also have
$\bigl|
 (\rho_{\alpha})_h^{\dagger}
 \bigr|_h
=\bigl|
 (\pi_{\alpha})_h^{\dagger}
 -\pi_{\alpha}'
 \bigr|_h
\leq C_{11}\exp(-\epsilon_0d)$.
We obtain
\[
 \bigl|
[f,\pi_{\alpha}^{\dagger}]
 \bigr|_h
=\bigl|
 [f,\pi_{\alpha}^{\dagger}-\pi_{\alpha}]
 \bigr|_h
\leq
 \bigl|
 [f,\rho_{\alpha}^{\dagger}]
 \bigr|_h
+\bigl|
[f,\rho_{\alpha}]
 \bigr|_h
\leq
 2(C_1M+C_2)
 C_{11}\exp(-\epsilon_0d).
\]
Then, we obtain the desired estimate.
\hfill\qed

\begin{cor}
\label{cor;15.8.18.22}
There exist positive constants $C_{21}$
and $\epsilon_{21}$ 
depending only on $R$, $R_1$, $R_2$, $r$ and $C_{10}$,
such that the following holds:
\begin{itemize}
\item
Take any $\alpha,\beta\in S$ with $\alpha\neq\beta$.
Let $s_{\alpha}$ and $s_{\beta}$ be local sections of
$E_{\alpha}$ and $E_{\beta}$.
Then, we have
\begin{equation}
\label{eq;15.8.18.20}
 \bigl|
 h(s_{\alpha},s_{\beta})
 \bigr|
\leq
 C_{21}\exp(-\epsilon_{21}d)
 \cdot
 |s_{\alpha}|_h\cdot
 |s_{\beta}|_h
\end{equation}
\end{itemize}
\end{cor}
\pf
We have
$h(s_{\alpha},s_{\beta})
=h(\pi_{\alpha}s_{\alpha},s_{\beta})
=h(s_{\alpha},\pi_{\alpha}^{\dagger}s_{\beta})
=h(s_{\alpha},(\pi_{\alpha}^{\dagger}-\pi_{\alpha})s_{\beta})$.
We have 
$\pi_{\alpha}^{\dagger}-\pi_{\alpha}
=\rho_{\alpha}^{\dagger}-\rho_{\alpha}$,
where $\rho_{\alpha}^{\dagger}$ denote the adjoint of
$\rho_{\alpha}$ with respect to $h$.
Then, we obtain (\ref{eq;15.8.18.20})
from Proposition \ref{prop;15.7.12.10}.
\hfill\qed

\subsubsection{Appendix: Preliminary from linear algebra}
\label{subsection;15.6.30.1}

Fix positive constants $G_i$  $(i=1,2)$.
Fix a positive integer $r>0$.
Let $U$ be any $r$-dimensional $\cnum$-vector space
with a Hermitian metric $h$.
Let $f$ be an endomorphism on $U$.
Suppose that we have
a finite subset $S\subset \cnum$
and a decomposition
$(U,f)=\bigoplus_{\alpha\in S} (U_{\alpha},f_{\alpha})$,
which is not necessarily orthogonal.
We impose the following conditions.
\begin{itemize}
\item
We set 
$d:=\min\bigl\{|\alpha-\beta|\,\big|\,\alpha,\beta\in S,\alpha\neq\beta
 \bigr\}$.
Then, 
we have $d\geq 1$ and $|f|_{h}\leq G_1d+G_2$.
\item
The eigenvalues $\gamma$ of $f_{\alpha}$
satisfy $|\gamma-\alpha|\leq d/100$.
\end{itemize}

Let $\pi_{\alpha}$ be the projection
with respect to the decomposition
$U=\bigoplus U_{\alpha}$.
Let $\pi_{\alpha}'$ denote the orthogonal projection
of $U$ onto $U_{\alpha}$.
We set $\rho_{\alpha}:=\pi_{\alpha}-\pi_{\alpha}'$.
To clarify the argument,
we recall the following lemma from \cite{Simpson90}
which was used in \cite{mochi2,Mochizuki-wild}.

\begin{lem}
\label{lem;15.7.1.1}
There exists a positive constant $B_1$
depending only on $r$ and $G_i$ $(i=1,2)$
such that
\[
\bigl|\pi_{\alpha}\bigr|_h\leq B_1,
\quad
 |\rho_{\alpha}|_h\leq B_1.
\]
\end{lem}
\pf
Because 
$\rho_{\alpha}$ and $\pi_{\alpha}'$ are orthogonal,
we have $|\rho_{\alpha}|_h\leq |\pi_{\alpha}|_h$.
Hence, it is enough to obtain the estimate for $\pi_{\alpha}$.
Let $\id_U$ denote the identity on $U$.
We have
\[
 \pi_{\alpha}
 =\frac{1}{2\pi\sqrt{-1}}
 \int_{\gamma_{\alpha}}\bigl(\zeta\id_U-f\bigr)^{-1}d\zeta.
\]
Here,
$\gamma_{\alpha}$ denotes the loop
$\gamma_{\alpha}(\theta)=
 \alpha+de^{\sqrt{-1}\theta}/10$
$(0\leq\theta\leq 2\pi)$.
There exist positive constants $B_i$ $(i=2,3)$
depending only on 
$r$ and $G_i$ $(i=1,2)$ such that 
\[
 \bigl|
 (\gamma_{\alpha}(\theta)\id_U-f)^{-1}
 \bigr|_h
\leq d^{-1}B_2
 \Bigl(
 \bigl(G_1+d^{-1}G_2\bigr)^r+1
 \Bigr)
\leq
 d^{-1}B_3.
\]
Thus, we obtain the claim of the lemma.
\hfill\qed

\vspace{.1in}

Let $f^{\dagger}_h$ denote the adjoint of $f$
with respect to $h$.
To clarify our argument,
we recall the following lemma
from \cite{Simpson90}
which was used in \cite{mochi2,Mochizuki-wild}.

\begin{lem}
\label{lem;15.7.1.2}
We have $\delta>0$ 
depending only on $r$ and $G_i$ $(i=1,2)$
such that 
$\bigl|
 \bigl[
 f_h^{\dagger},\pi_{\alpha}
 \bigr]\bigr|_h
 \geq
 \delta\cdot d\cdot|\rho_{\alpha}|_h$.
\end{lem}
\pf
We take a numbering of the elements of $S$,
i.e.,
$S=\{\alpha_1,\ldots,\alpha_m\}$.
We impose $\alpha_1=\alpha$.
We set $F_{j}(U):=\bigoplus_{i\leq j} U_{\alpha_i}$
and $F_{<j}(U):=\bigoplus_{i<j}U_{\alpha_i}$.
Let $U_{j}'$ be the orthogonal complement
of $F_{<j}(U)$ in $F_j(U)$.
We have the orthogonal decomposition
$F_{j}(U)=\bigoplus_{i\leq j} U_{i}'$.
Because $f(F_j)\subset F_j$,
we have the decomposition
$f=\sum_{i\leq j}f_{ij}$,
where $f_{ij}:U_{j}'\lrarr U_{i}'$.
As the adjoint,
we have
$f^{\dagger}_h=\sum_{i\leq j} (f_{ij})^{\dagger}_h$,
where
$(f_{ij})^{\dagger}_h:
 U_i'\lrarr U_j'$.
We set $f^{\dagger}_{ij}:=(f_{ji})^{\dagger}_h$.
Then, we have
$f^{\dagger}_h=\sum_{i\geq j}f^{\dagger}_{ij}$,
where $f_{ij}^{\dagger}:U_{j}'\lrarr U_{i}'$.

For each $i$,
we take an orthonormal base
$e^{(i)}_1,\ldots,e^{(i)}_{r_i}$ of $U_i'$
for which 
$f^{\dagger}_{ii}$ is represented by a lower triangular matrix $A_i$.
Let $\Gamma_i\in\End(U_i')$ be determined by
$\Gamma_i(e^{(i)}_k)=\alpha^{(i)}_ke^{(i)}_k$,
where $\alpha^{(i)}_k$ denote
the $(k,k)$-entry of $A_i$.
Then, $f^{\dagger}_{ii}-\Gamma_i$ is nilpotent.
We put $\Gamma=\sum\Gamma_i\in\End(U)$.
Then, $f^{\dagger}_h-\Gamma$ is nilpotent,
and we have 
$|f^{\dagger}_h-\Gamma|_h\leq G_1d+G_2$.

Let $\Pi$ denote the orthogonal projection of
$\End(E)=\bigoplus \Hom(U_j',U_i')$
onto $\bigoplus_{j>i} \Hom(U_j',U_i')$.
Let $F_1$ be the endomorphism on
$\bigoplus_{j>i}\Hom(U_j',U_i')$ 
given by
$F_1(B):=\Pi([f^{\dagger}_h,B])$.
Let $F_2$ be the endomorphism on
$\bigoplus _{j>i}\Hom(U_j',U_i')$
given by
$F_2(B):=[\Gamma,B]=\Pi([\Gamma,B])$.

Let us observe that
$F_1-F_2$ is nilpotent.
We have the nilpotent endomorphism $\Ftilde$
on $\End(E)$ given by
$\Ftilde(B)=[f^{\dagger}_h-\Gamma,B]$,
which preserves 
$\bigoplus_{j\leq i}\Hom(U_j',U_i')$.
Then, $F_1-F_2$ is equal to
the endomorphism on
$\bigoplus_{j>i}\Hom(U_j',U_i')\simeq
 \End(E)\big/\bigoplus_{j\leq i}\Hom(U_j',U_i')$
induced by $\Ftilde$,
and hence $F_1-F_2$ is nilpotent.

There exists a constant $B_{11}>0$
depending only on $r$ 
such that 
$|F_1-F_2|_h\leq B_{11}(G_1d+G_2)$.
Moreover,
$F_2$ is invertible,
and 
there exists a constant $B_{12}>0$
depending only on $r$ 
such that 
$|F_2^{-1}|_h\leq B_{12}d^{-1}$.
Hence, we obtain that 
$|F_1^{-1}|_h\leq B_{13}d^{-1}(1+(G_1+G_2/d)^{r})$
for a constant $B_{13}>0$
depending only on $r$.

Note that we have 
$\Pi([f_h^{\dagger},\pi_{\alpha}'])=0$.
Hence,
we have the following:
\[
 \bigl|
 [f_h^{\dagger},\pi_{\alpha}]
 \bigr|_h
\geq
 \Bigl|
 \Pi([f_h^{\dagger},\pi_{\alpha}])
 \Bigr|_h
=\Bigl|
 \Pi([f_h^{\dagger},\rho_{\alpha}])
 \Bigr|_h
=\bigl|
 F_1(\rho_{\alpha})
 \bigr|_h
\geq
 B^{-1}_{13}
 d
\bigl(1+(G_1+G_2d^{-1})^r\bigr)^{-1}
|\rho_{\alpha}|_h
\]
Thus, we obtain the claim of the lemma.
\hfill\qed

\subsection{Asymptotic decoupling of harmonic bundles on discs}
\label{subsection;15.7.17.30}

We continue to use the setting in 
\S\ref{subsection;15.7.17.20}.
We further impose the following condition:
\begin{itemize}
\item
$\rank E_{\alpha}=1$
for each $\alpha\in S$.
\end{itemize}
In other words,
we assume to have holomorphic functions
$g_{\alpha}$ $(\alpha\in S)$
on $\Delta(R)$
such that
$\theta=\bigoplus g_{\alpha}\id_{E_{\alpha}}dz$.
(Recall the generically regular semisimple condition
in \S\ref{subsection;15.7.22.1}.)
By the condition,
we have
$|g_{\alpha}(P)-\alpha|\leq d/100$.
In this setting,
we explain that Simpson's main estimate implies
the asymptotic decoupling of the Hitchin equation.

Let $g_{\cnum}$ denote the Euclidean metric $dz\,d\zbar$
of $\cnum$.
For any section $s$ of $\End(E)\otimes\Omega^{p,q}$
on $\Delta(R')$ $(R'>0)$,
let $|s|_{h,g_{\cnum}}$ denote the function on $\Delta(R')$
by taking the norm of $s$ 
at each $P\in \Delta(R')$
with respect to $h$ and $g_{\cnum}$.

\subsubsection{Decay of the curvatures}

Let $R(h)$ denote the curvature
of the Chern connection of $(E,\delbar_E,h)$.
We obtain the following ``asymptotic decoupling''
of the Hitchin equation
$R(h)+[\theta,\theta_h^{\dagger}]=0$.

\begin{thm}
\label{thm;15.7.12.11}
There exist positive constants
$C_{30}$ and $\epsilon_{30}$,
depending only on 
$R$, $R_1$, $R_2$, $r$ and $C_{10}$ 
such that
$\bigl|
 R(h)
\bigr|_{h,g_{\cnum}}
=\bigl|
[\theta,\theta_h^{\dagger}]
 \bigr|_{h,g_{\cnum}}
\leq C_{30}
 \exp(-\epsilon_{30}d)$
on $\Delta(R_2)$.
\end{thm}
\pf
It is enough to obtain the estimate for
$[\theta,\theta^{\dagger}_h]$.
We have the decomposition
$f^{\dagger}_h
=\sum_{\alpha,\beta} \pi_{\alpha}\circ f^{\dagger}_h\circ\pi_{\beta}$.
By Corollary \ref{cor;15.7.22.10},
if $\alpha\neq \beta$,
we have $C_{31}>0$,
depending only on 
$R$, $R_1$, $R_2$, $r$ and $C_{10}$,
such that
\[
 \bigl|
 \pi_{\alpha}\circ f^{\dagger}_h\circ\pi_{\beta}
 \bigr|_h
=\bigl|
 [\pi_{\alpha},f^{\dagger}_h]\circ\pi_{\beta}
 \bigr|_h
\leq
 C_{31}
 \exp(-\epsilon_{20}d).
\]
For $\alpha\neq \beta$,
we have
$[f_{\alpha},\pi_{\beta}\circ f^{\dagger}_h\circ\pi_{\beta}]=0$.
Because $\rank E_{\alpha}=1$,
we also have
$\bigl[
 f_{\alpha},\pi_{\alpha}\circ f^{\dagger}_h\circ \pi_{\alpha}
 \bigr]=0$.
Then, 
we obtain the estimate for
$[\theta,\theta^{\dagger}_h]$.
\hfill\qed

\subsubsection{The connections and the projections}
\label{subsection;16.5.23.2}

Let $\del_{E}$ denote the $(1,0)$-part
of the Chern connection associated to $h$ and $\delbar_E$.
According to Proposition \ref{prop;15.7.12.10},
the decomposition $E=\bigoplus E_{\alpha}$
is almost orthogonal.
Let us see that
such an almost orthogonality holds
at the level of the first derivative
in the sense that $\del_E\pi_{\alpha}$ is very small.

\begin{prop}
\label{prop;15.7.1.10}
Take $0<R_3<R_2$.
There exist positive constants $\epsilon_{40}$ and $C_{40}$
depending only on 
$R$, $R_1$, $R_2$, $R_3$, $r$ and $C_{10}$
such that 
$\bigl|\del_E\pi_{\alpha}\bigr|_{h,g_{\cnum}}=
\bigl|\delbar_E\pi^{\dagger}_{\alpha}\bigr|_{h,g_{\cnum}}
\leq
 C_{40} \exp(-\epsilon_{40}d)$
on $\Delta(R_3)$.
\end{prop}
\pf
It is enough to obtain the estimate for
$\del_E\pi_{\alpha}$.
In the following,
 the constants may depend only on 
$R$, $R_1$, $R_2$, $r$ and $C_{10}$.
We have $C_{41}>0$ such that
the following holds on $\Delta(R_2)$:
\[
 \bigl|
 \delbar_E\del_{E}\pi_{\alpha}
 \bigr|_{h,g_{\cnum}}
=\bigl|
 [R(h),\pi_{\alpha}]
 \bigr|_{h,g_{\cnum}}
\leq
 C_{41}
 \exp(-\epsilon_{30} d).
\]
Because $\del_E \pi_{\alpha}^{\dagger}=0$,
we have the following on $\Delta(R_2)$:
\[
 \bigl|
 \delbar_E\del_E(\pi_{\alpha}-\pi_{\alpha}^{\dagger})
 \bigr|_{h,g_{\cnum}}
\leq
 C_{41}
 \exp(-\epsilon_{30}d)
\]
We already have
$\bigl|
 \pi_{\alpha}-\pi_{\alpha}^{\dagger}
 \bigr|_h
\leq C_{42}\exp(-\epsilon_0 d)$ for 
a constant $C_{42}>0$.
We may assume $\epsilon_{30}<\epsilon_0$.
Hence, we have a constant $C_{43}>0$
such that 
$\bigl\|
 \pi_{\alpha}-\pi_{\alpha}^{\dagger}
 \bigr\|_{h,g_{\cnum},L_2^p}
\leq
 C_{43}
 \exp(-\epsilon_{30} d)$
for a large $p>1$.
Here,
$\|\cdot\|_{h,g_{\cnum},L_2^p}$
denote the $L_2^p$-norm with respect to
$h$ and $g_{\cnum}$ on $\Delta(R_2)$.
Then, we obtain the claim of the lemma.
\hfill\qed

\vspace{.1in}
Recall that $\pi_{\alpha}'$
denote the orthogonal projection of $E$
onto $E_{\alpha}$.
Let us see that it is almost holomorphic.

\begin{prop}
\label{prop;15.7.1.11}
We have a positive constant $C_{50}$,
depending only on 
$R$, $R_1$, $R_2$, $R_3$, $r$ and $C_{10}$
such that the following holds
on $\Delta(R_3)$:
\[
 \bigl|
 \delbar_E\pi_{\alpha}'
 \bigr|_{h,g_{\cnum}}
=\bigl|
 \del_E\pi_{\alpha}'
 \bigr|_{h,g_{\cnum}}
\leq
 C_{50}
 \exp(-\epsilon_{40}d).
\]
\end{prop}
\pf
It is enough to prove the estimate for
$\delbar_E\pi_{\alpha}'$.
Let $E=E_{\alpha}\oplus E_{\alpha}^{\bot}$
denote the orthogonal decomposition.
We may naturally regard
$\rho_{\alpha}:=\pi_{\alpha}-\pi_{\alpha}'$
as a morphism
$E_{\alpha}^{\bot}\lrarr E_{\alpha}$.
We may also regard
$\rho_{\alpha}^{\dagger}:=
\pi_{\alpha}^{\dagger}-\pi_{\alpha}'$
as a morphism
$E_{\alpha}\lrarr E_{\alpha}^{\bot}$.

We have the induced holomorphic structure
on $E_{\alpha}^{\bot}\simeq E/E_{\alpha}$.
Let $\delbar_E^{(0)}$ be the holomorphic structure
on $E_{\alpha}\oplus E_{\alpha}^{\bot}$
obtained as the direct sum.
We have 
$\delbar_E=\delbar_E^{(0)}+\kappa$,
where
$\kappa$ is a section of
$\Hom(E_{\alpha}^{\bot},E_{\alpha})\otimes\Omega^{0,1}$.
We have
\begin{equation}
\label{eq;15.8.8.2}
 \delbar_E(\pi_{\alpha}-\pi_{\alpha}^{\dagger})
=\delbar_E\rho_{\alpha}
-\delbar_E\rho_{\alpha}^{\dagger}
=\delbar_E^{(0)}\rho_{\alpha}
-\delbar_E^{(0)}\rho_{\alpha}^{\dagger}
-\kappa\circ\rho_{\alpha}^{\dagger}
+\rho_{\alpha}^{\dagger}\circ\kappa
\end{equation}
Note that 
$\delbar_E^{(0)}(\rho_{\alpha})$,
$\delbar_E^{(0)}(\rho_{\alpha}^{\dagger})$,
$\kappa\circ\rho_{\alpha}^{\dagger}$
and 
$\rho_{\alpha}^{\dagger}\circ\kappa$
are sections of
$\Hom(E_{\alpha}^{\bot},E_{\alpha})\otimes\Omega^{0,1}$,
$\Hom(E_{\alpha},E_{\alpha}^{\bot})\otimes\Omega^{0,1}$,
$\Hom(E_{\alpha},E_{\alpha})\otimes\Omega^{0,1}$
and
$\Hom(E_{\alpha}^{\bot},E_{\alpha}^{\bot})\otimes\Omega^{0,1}$,
respectively.
Note that the bundles are orthogonal
with respect to $h$ and $g_{\cnum}$.
In general,
for any orthogonal decomposition of bundles
$\nbigv=\bigoplus\nbigv_i$
and for any section $s=\sum s_i$ of $\nbigv$,
the norms of $s_i$ are smaller than 
the norm of $s$.
Hence, we obtain the following
from (\ref{eq;15.8.8.2}):
\begin{equation}
\label{eq;15.8.8.1}
 \bigl|
 \delbar^{(0)}_E\rho_{\alpha}
 \bigr|_{h,g_{\cnum}}
\leq
 \bigl|
 \delbar_E(\pi_{\alpha}-\pi_{\alpha}^{\dagger})
 \bigr|_{h,g_{\cnum}}
=\bigl|
 \delbar_E(\pi_{\alpha}^{\dagger})
 \bigr|_{h,g_{\cnum}}
\end{equation}
We also have
$\delbar_E\pi_{\alpha}'
=-\delbar_E\rho_{\alpha}
=-\delbar_E^{(0)}\rho_{\alpha}$.
Hence, 
we obtain
$\bigl|
 \delbar_E\pi_{\alpha}'
 \bigr|_{h,g_{\cnum}}
\leq
 \bigl|
 \delbar_E(\pi_{\alpha}^{\dagger})
 \bigr|_{h,g_{\cnum}}$
from (\ref{eq;15.8.8.1}).
Then, the claim of Proposition \ref{prop;15.7.1.11}
follows from the estimate
in Proposition \ref{prop;15.7.1.10}.
\hfill\qed

\subsubsection{The decay of the curvatures on the line bundles}

Let $\delbar_{\alpha}$ denote the holomorphic structure
of $E_{\alpha}$.
Let $h_{\alpha}$ be the restriction of $h$ to $E_{\alpha}$.
Let $\del_{\alpha}$ denote the $(1,0)$-part of
the Chern connection of
$(E_{\alpha},\delbar_{\alpha},h_{\alpha})$.
Let $R(h_{\alpha})$ denote the curvature 
of the connection $\delbar_{\alpha}+\del_{\alpha}$.
We have
$\del_{\alpha}s=
 \pi_{\alpha}'\circ\del_E s$
and $\delbar_{\alpha}s=\delbar_Es$
for any section $s$ of $E_{\alpha}$.

\begin{prop}
\label{prop;15.7.12.12}
We have a positive constant $C_{60}$
depending only on 
$R$, $R_1$, $R_2$, $R_3$, $r$ and $C_{10}$
such that
$\bigl|
 R(h_{\alpha})
 \bigr|_{h_{\alpha},g_{\cnum}}
\leq
 C_{60}\exp(-\epsilon_{40}d)$
on $\Delta(R_3)$.
\end{prop}
\pf
In the following,
 the constants may depend only on 
$R$, $R_1$, $R_2$, $R_3$, $r$ and $C_{10}$.
Let $s$ be any section of $E_{\alpha}$.
We have
\[
 \delbar_E\circ\pi_{\alpha}'\bigl(
 \del_E s\bigr)
=\delbar_E\circ
 \pi_{\alpha}'\circ\del_E\bigl(
 \pi_{\alpha}s\bigr)
=\delbar_E\circ\pi_{\alpha}(\del_E s)
+\delbar_E\circ\pi_{\alpha}'\circ(\del_E\pi_{\alpha})(s)
=\pi_{\alpha}(\delbar_E\del_E s)
+\delbar_E\bigl(
 \pi_{\alpha}'(\del_E\pi_{\alpha})s
 \bigr).
\]
We also have
\[
 \pi_{\alpha}'\circ\del_E\circ\delbar_E s
=\pi_{\alpha}'\circ\del_E(\pi_{\alpha}\delbar_E s)
=\pi_{\alpha}\del_E\delbar_E s
+\pi_{\alpha}'\del_E(\pi_{\alpha})\delbar_E s.
\]
Hence, it is enough to obtain an estimate of
$\delbar_E(\pi_{\alpha}'\circ\del_E\pi_{\alpha})
=\delbar_E(\pi_{\alpha}')\circ\del_E(\pi_{\alpha})
+\pi_{\alpha}'\circ\delbar_E\del_E\pi_{\alpha}$.
By Proposition \ref{prop;15.7.1.10}
and Proposition \ref{prop;15.7.1.11},
we have $C_{61}>0$
such that 
$\bigl|
 \delbar_E(\pi_{\alpha}')\circ\del_E(\pi_{\alpha})
 \bigr|_{h,g_{\cnum}}
\leq
 C_{61}\exp(-\epsilon_{40}d)$.
We also have
$\bigl|
 \pi_{\alpha}'\circ
\delbar_E\del_E\pi_{\alpha}
 \bigr|_{h,g_{\cnum}}
=\bigl|
 \pi_{\alpha}'\circ
 [R(h),\pi_{\alpha}]
 \bigr|_{h,g_{\cnum}}
\leq 
C_{62}
 \exp(-\epsilon_{40}d)$
for a constant $C_{62}>0$.
Then, the claim of the lemma follows.
\hfill\qed

\subsubsection{Approximation of the flat connections}
\label{subsection;16.5.23.10}

We consider the flat connection
$\DD^1:=\delbar_E+\del_E+\theta+\theta^{\dagger}$
on $E$.
We also have the connection
$\DD^1_0$ 
on $E=\bigoplus_{\alpha\in S} E_{\alpha}$
given by
\[
 \DD^1_0=\bigoplus_{\alpha\in S}
 \Bigl(
 \bigl(\delbar_{\alpha}+\del_{\alpha}\bigr)
+(g_{\alpha}dz+\overline{g_{\alpha}}d\zbar)
 \id_{E_{\alpha}}
 \Bigr).
\]
Here, $g_{\alpha}$ are holomorphic functions
such that $\theta=\bigoplus g_{\alpha}\id_{E_{\alpha}}dz$,
as in the beginning of \S\ref{subsection;15.7.17.30}.

\begin{cor}
\label{cor;15.7.17.40}
There exists a positive constant $C_{70}$,
depending only on 
$R$, $R_1$, $R_2$, $R_3$, $r$ and $C_{10}$
such that 
$\bigl|
 \DD^1-\DD^1_0
 \bigr|_{h,g_{\cnum}}
\leq
 C_{70}
 \exp(-\epsilon_{40}d)$
on $\Delta(R_3)$.
\end{cor}
\pf
It is enough to obtain an estimate of
$\del_E-\bigoplus_{\alpha\in S} \del_{\alpha}$.
We have a constant $C_{71}>0$,
depending only on 
$R$, $R_1$, $R_2$, $R_3$, $r$ and $C_{10}$
such that 
$\bigl|
 \pi_{\alpha}'\circ
 (\del_E\pi_{\alpha})
 \bigr|_{h,g_{\cnum}}
\leq
 C_{71}\exp(-\epsilon_{40}d)$.
Hence, for any section $s$ of $E_{\alpha}$,
we have
\[
\bigl|
 \pi_{\alpha}\circ\del_{E}s-\del_{\alpha}s
 \bigr|_{h,g_{\cnum}}
=\bigl|
 \pi_{\alpha}'\circ\pi_{\alpha}\circ\del_Es
-\pi_{\alpha}'\del_E(\pi_{\alpha}s)
 \bigr|_{h,g_{\cnum}}
=\bigl|
 \pi_{\alpha}'\circ
 \del_E\pi_{\alpha}(s)
 \bigr|_{h,g_{\cnum}}
\leq
 C_{71}\exp(-\epsilon_{40}d)
 |s|_h.
\]
Then, the claim of the lemma follows.
\hfill\qed

\subsubsection{Higher derivatives}
\label{subsection;16.5.23.20}

We take a numbering
$\{\alpha_1,\ldots,\alpha_r\}$ on $S$.
For each $i$,
we can take a holomorphic frame $u_{i}$ of $E_{\alpha_i}$
such that
$|u_{i|0}|_h=1$
and
$|\del_{\alpha_i}u_{i}|_{h,g_{\cnum}}\leq 
 C_{80}\exp(-\epsilon_{80}d)$
on $\Delta(R_3)$
for some positive constants
$C_{80}$ and $\epsilon_{80}$,
depending only on $R$, $R_1$, $R_2$, $R_3$, $r$ and $C_{10}$.
Let us sketch how to obtain such sections
in an elementary way.
For a real coordinate $z=x+\sqrt{-1}y$,
we can take a section
$s_i$ of $E_{\alpha_i}$ on $\Delta(R_{2})$
such that
$\bigl(\delbar_{\alpha_i}+\del_{\alpha_i}\bigr)s_i
=s_i\cdot \nu_{i}dx=s_i\cdot \nu_i(dz+d\zbar)/2$
such that $\nu_i(x,0)=0$.
The curvature form is given by
$\del_y\nu_i\,dy\,dx$.
By Proposition \ref{prop;15.7.12.12},
we have $|\nu_i|\leq C_{80}'\exp(-\epsilon'_{80}d)$.
Take $R_2'$ such that $R_3<R_2'<R_2$.
We can take
a function $\rho_i$ on $\Delta(R_2')$
satisfying
$\del_{\zbar}\rho_i=\nu_i/2$,
$|\rho_i|\leq C_{80}''\exp(-\epsilon''_{80}d)$
and
$\|\del_z\rho_i\|_{L^p(\Delta(R_2'))}
 \leq C_{80,p}''\exp(-\epsilon''_{80})$
for any $p\geq 1$.
We set $v_i:=s_ie^{-\rho_i}$ on $\Delta(R_2')$.
Then, we have
$\delbar_{\alpha_i}v_i=0$
and $\del_{\alpha_i}v_i=v_i\cdot \kappa_idz$,
where
$\|\kappa_i\|_{L^p(\Delta(R_2'))}\leq 
 C_{80,p}^{(3)}\exp(-\epsilon^{(3)}_{80}d)$.
Because $\del_{\zbar}\kappa_id\,\zbar\,dz$
is the curvature of $E_{\alpha}$,
we have
 $\bigl|\del_{\zbar}\kappa_i\bigr|
\leq
 C^{(4)}_{80}\exp(-\epsilon^{(4)}_{80}d)$.
Hence, we obtain
$|\kappa_i|\leq
 C^{(5)}_{80}\exp(-\epsilon^{(5)}_{80}d)$
on $\Delta(R_3)$.
By adjusting the norm of $v_i$ at the origin,
we obtain the desired section $u_i$.
Note that we have
$\bigl|
 \log|u_{i}|_h
 \bigr|
\leq
 C_{81}\exp(-\epsilon_{81}d)$
for some positive constants
$C_{81}$ and $\epsilon_{81}$
depending only on $R$, $R_1$, $R_2$, $R_3$, $r$ and $C_{10}$.

We obtain a frame $\vecu=(u_1,\ldots,u_r)$ of $E$
on $\Delta(R_3)$.
Let $H$ be the $r$-square Hermitian matrix valued
function given by
$H_{ij}=h(u_i,u_j)$.
Let $\Theta$ be the holomorphic $r$-square matrix valued function 
such that $\Theta_{ii}=g_{\alpha_i}$
and $\Theta_{ij}=0$ $(i\neq j)$.
We have
$\theta\vecu=\vecu\cdot\Theta dz$.
Let $\Theta^{\dagger}$ be the $r$-square matrix valued function
given by
$\Theta^{\dagger}=\Hbar^{-1}\lefttop{t}\overline{\Theta}\Hbar$.
We have
$\theta^{\dagger}_h\vecu=\vecu\Theta^{\dagger}\,d\zbar$.

Because $H_{ii}=|u_i|_{h}^2$,
we have
$\bigl|\log H_{ii}\bigr|
\leq 2C_{81}\exp(-\epsilon_{81}d)$,
as remarked above.
We also have positive constants $C_{82}$
and $\epsilon_{82}$ 
depending only on 
$R$, $R_1$, $R_2$, $R_3$, $r$ and $C_{10}$
such that the following holds:
\begin{equation}
 \label{eq;15.7.24.32}
 |H_{ij}|\leq C_{82}\exp(-\epsilon_{82}d),
\quad\quad(i\neq j)
\end{equation}
\begin{equation}
 \label{eq;15.7.24.30}
 \bigl|
 \del_z H_{ij}
 \bigr|
=\bigl|
 \del_{\zbar} H_{ij}
 \bigr|
\leq 
 C_{82}\exp(-\epsilon_{82}d)
\end{equation}
\begin{equation}
\label{eq;15.7.24.31}
 \bigl|
 \del_{\zbar}\del_z H_{ij}
 \bigr|
\leq 
 C_{82}\exp(-\epsilon_{82}d)
\end{equation}
Indeed,  
(\ref{eq;15.7.24.32}) follows from
Corollary \ref{cor;15.8.18.22}.
We have
$\del h(u_i,u_j)=
 h\bigl(
 (\del_E\pi_{\alpha_i})u_i,u_j
 \bigr)
+h\bigl(
 \pi_{\alpha_i}(\del_Eu_i),u_j
 \bigr)$.
As in the proof of Corollary \ref{cor;15.7.17.40},
we have 
$\bigl|
 \pi_{\alpha_i}\del_Eu_i-\del_{\alpha_i}u_i
 \bigr|_{h,g_{\cnum}}
\leq C_{82}'\exp(-\epsilon'_{82}d)|u_i|_h$.
Then, we obtain (\ref{eq;15.7.24.30})
from Proposition \ref{prop;15.7.1.10}
and our choice of $u_i$.
We obtain (\ref{eq;15.7.24.31})
from the estimate for the curvature
$\delbar\bigl(\Hbar^{-1}\del \Hbar\bigr)$
and (\ref{eq;15.7.24.30}).

\begin{lem}
\label{lem;15.7.24.60}
Take any $\vecell=(\ell_1,\ell_2)\in
 \bigl(
 \seisuu_{\geq 0}\times\seisuu_{\geq 0}
 \bigr)\setminus\{(0,0)\}$
and  $R_4<R_3$.
We also take any $p>1$.
Then,
we have positive constants $C_{83,\vecell,p}$
and $\epsilon_{83,\vecell,p}$
depending only on 
$R$, $R_1$, $R_2$, $R_3$, $R_4$, $r$, $C_{10}$,
$\vecell$ and $p$,
such that the following holds:
\[
 \bigl\|
 \del_z^{\ell_1}\del_{\zbar}^{\ell_2}H
_{|\Delta(R_4)}
 \bigr\|_{g_{\cnum},L^p}
\leq
 C_{83,\vecell,p}
 \exp\bigl(
 -\epsilon_{83,\vecell,p}
 d
 \bigr).
\]
\end{lem}
\pf
The proof is given in a standard inductive argument
using the Hitchin equation and the elliptic regularity.
Because $\Theta_{ii}$ in the Hitchin equation can be large,
we give a rather detailed argument.

Let us consider the case $\vecell=(2,0)$.
By (\ref{eq;15.7.24.30}) and (\ref{eq;15.7.24.31}),
we have positive constants
$C_{84}$ and $\epsilon_{84}$
such that 
$\|\del_z H_{|\Delta(R_4)}\|_{L_1^p}
\leq
 C_{84}\exp(-\epsilon_{84}d)$.
So, we obtain
$\|\del_z\del_z H_{|\Delta(R_4)}\|_{L^p}
\leq
 C_{85}\exp(-\epsilon_{85}d)$
for some $C_{85}>0$ and $\epsilon_{85}>0$.
Similarly,
we can obtain the estimate for
$\|\del_{\zbar}\del_{\zbar}H_{\Delta(R_4)}\|_{L^p}$.

Suppose the claim has already been proved
for any $R_4$ and $p$ if $\ell_1+\ell_2<k$.
We consider $\vecell=(\ell_1,\ell_2)$
satisfying $\ell_1+\ell_2=k$ and $\ell_i>0$.
The Hitchin equation is described as follows:
\begin{equation}
 \label{eq;15.7.24.50}
 \Hbar^{-1}
 \del_{\zbar}\del_{z}\Hbar
-\Hbar^{-1}\del_{\zbar}\Hbar
 \cdot
 \Hbar^{-1}\del_z\Hbar
-\bigl[
 \Theta,\Hbar^{-1}(\lefttop{t}\overline{\Theta})\Hbar
 \bigr]
=0
\end{equation}
For each $\ell$,
we have $C_{86,\ell}$
depending on $C_{10}$ and $\ell$,
such that the following holds:
\[
 \Bigl|
 \del_z^{\ell}\Theta_{ii}
 \Bigr|
=
 \Bigl|
 \del_{\zbar}^{\ell}\overline{\Theta}_{ii}
 \Bigr|
\leq
 C_{86,\ell}d
\]
Applying $\del_z^{\ell_1-1}\del_{\zbar}^{\ell_2-1}$
to (\ref{eq;15.7.24.50}),
we obtain the following:
\[
 \Hbar^{-1}\del^{\ell_1}_z\del^{\ell_2}_{\zbar}\Hbar
-\bigl[
 \del_z^{\ell_1-1}\Theta,
 \Hbar^{-1}
 \lefttop{t}\overline{(\del_z^{\ell_2-1}\Theta)}
 \Hbar
 \bigr]
+G=0
\]
Here,
$G$ is expressed as a linear combination of
terms which contains derivatives of $\Hbar$.
Hence, we have
$\|G_{|\Delta(R_4)}\|_{L^p}
\leq C_{87}\exp(-\epsilon_{87}d)$.
We also have 
$\Bigl\|
 \bigl[
 \del_z^{\ell_1-1}\Theta,
 \Hbar^{-1}
 \lefttop{t}\overline{(\del_z^{\ell_2-1}\Theta)}
 \Hbar
 \bigr]_{|\Delta(R_4)}
 \Bigr\|_{L^p}
\leq
 C_{88}\exp(-\epsilon_{88}d)$,
by 
(\ref{eq;15.7.24.32})
and  (\ref{eq;15.7.24.50}).
Hence, we obtain the desired estimate 
if $\ell_1+\ell_2=k$ and $\ell_i\neq 0$.
We can deal with the cases
$\vecell=(k,0),(0,k)$
by using the elliptic regularity.
Then, by an inductive argument,
we can obtain the claim of Lemma \ref{lem;15.7.24.60}.
\hfill\qed

\begin{cor}
\label{cor;15.7.24.100}
Take any $\vecell=(\ell_1,\ell_2)\in
 \bigl(
 \seisuu_{\geq 0}\times\seisuu_{\geq 0}
 \bigr)\setminus\{(0,0)\}$
and  $R_4<R_3$.
Then,
we have positive constants $C_{90,\vecell}$
and $\epsilon_{90,\vecell}$
depending only on 
$R$, $R_1$, $R_2$, $R_3$, $R_4$, $r$, $C_{10}$
and $\vecell$
such that the following holds:
\[
\sup
_{\Delta(R_4)}
\bigl|
 \del_z^{\ell_1}\del_{\zbar}^{\ell_2}H
\bigr|
\leq
 C_{90,\vecell}
 \exp\bigl(
 -\epsilon_{90,\vecell}
 d
 \bigr).
\]
\hfill\qed
\end{cor}

\subsection{Hitchin WKB-problem}

\subsubsection{Preliminary}

Let $V$ be an $r$-dimensional complex vector space.
For Hermitian metrics $h_j$ $(j=1,2)$,
we can take a base $e_1,\ldots,e_r$ of $V$
which is orthogonal with respect to both 
$h_1$ and $h_2$.
We have the real numbers $\alpha_i$
determined by
$\kappa_i:=\log|e_i|_{h_2}-\log|e_i|_{h_1}$.
We impose 
$\kappa_1\geq \kappa_2\geq\cdots\geq \kappa_r$.
Then, we set
\[
 \vec{d}(h_1,h_2):=
 \bigl(
 \kappa_1,\ldots,\kappa_r
 \bigr)\in\real^r.
\]

Let $V_j$ $(j=1,2)$ be $r$-dimensional complex vector spaces
with Hermitian metrics $h_j$.
Let $f:V_1\lrarr V_2$ be a linear isomorphism.
We define the Hermitian metric
$f^{\ast}h_2$ on $V_1$
by $f^{\ast}h_2(u,v)=h_2(f(u),f(v))$.
We recall the following lemma from \cite{KNPS},
which can be easily proved.
\begin{lem}
\label{lem;15.8.18.1}
We have orthonormal frames
$\vecq(s)=(q_1(s),\ldots,q_r(s))$ $(s=0,1)$ on $V_s$
such that 
$f(q_j(0))=e^{\beta_j}q_j(1)$ $(j=1,\ldots,r)$,
where $\beta_j$ are real numbers
satisfying
$\beta_1\geq\beta_2\geq\cdots\geq \beta_r$.
In this case,
we have
\begin{equation}
\label{eq;15.7.22.40}
 \vec{d}(h_1,f^{\ast}h_2)
=(\beta_1,\ldots,\beta_r).
\end{equation}
\hfill\qed
\end{lem}

We set 
$\|f\|_{\op}:=
\sup\Bigl\{
 \bigl|f(u)\bigr|_{h_2}\,\Big|\,
 u\in V_1,
 |u|_{h_1}=1
 \Bigr\}$.
By Lemma \ref{lem;15.8.18.1},
we have 
$\beta_1=\log\|f\|_{\op}$.
We also have
\[
 \sum_{j=1}^k\beta_j
=\log\bigl\|
 \bigwedge^kf
 \bigr\|_{\op},
\]
where $\bigwedge^kf:\bigwedge^kV_1\lrarr \bigwedge^kV_2$
are the induced morphisms.
Hence, we have the following formula,
as noted in \cite{KNPS}:
\begin{equation}
 \label{eq;15.7.22.50}
 \beta_k=
 \log\big\|\bigwedge^kf\bigr\|_{\op}
-\log\bigl\|\bigwedge^{k-1}f\bigr\|_{\op}. 
\end{equation}

\subsubsection{Hitchin WKB-problem}

Let $X$ be any complex curve.
Let $\phi_i$ $(i=1,\ldots,r)$
be holomorphic $1$-forms on $X$.
Let $[0,1]:=\{s\in\real\,|\,0\leq s\leq 1\}$ be
the closed interval.
Let $\gamma:[0,1]\lrarr X$ be a $C^{1}$-map.
Suppose that it is a non-critical path
in the sense of \cite{KNPS},
i.e., the following holds:
\begin{itemize}
\item
 At any point $s\in[0,1]$,
we have
 $\gamma^{\ast}\Re(\phi_i)_s
 \neq
 \gamma^{\ast}\Re(\phi_j)_s$ 
$(i\neq j)$.
\end{itemize}
We have the expression
$\gamma^{\ast}(\phi_i)_s=a_i(s)\,ds$
for some $C^{\infty}$-functions $a_i:[0,1]\lrarr \cnum$.
We may assume
$\Re a_i(s)<\Re a_j(s)$ $(i<j)$ for any $s$.
We set
\[
 \alpha_i:=-\int_0^1\Re(a_i)\,ds.
\]
We have
$\alpha_1>\alpha_2>\cdots>\alpha_r$.

\vspace{.1in}

Let $(E,\delbar_E,\theta,h)$ be a harmonic bundle 
of rank $r$ on $X$.
We suppose the following.
\begin{itemize}
\item
We have the decomposition
$(E,\delbar_E,\theta)=
 \bigoplus_{i=1}^r(E_i,\delbar_{E_i},t\phi_i\id_{E_i})$
for some $t>0$,
where $\rank E_i=1$.
\end{itemize}

We have the associated flat connection
$\DD^1=\delbar_E+\del_E+\theta+\theta^{\dagger}$
on $E$.
Let $\Pi_{\gamma}:E_{|\gamma(0)}\lrarr E_{|\gamma(1)}$ 
be the parallel transport of $\DD^1$
along $\gamma$.
We have the metrics $h_{\gamma(0)}$ 
on the fibers $E_{|\gamma(\kappa)}$ $(\kappa=0,1)$
induced by $h$.
We obtain the metric 
$\Pi_{\gamma}^{\ast}h_{\gamma(1)}$ on $E_{|\gamma(0)}$
induced by $h_{\gamma(1)}$ and $\Pi_{\gamma}$.

The following theorem was conjectured in 
\cite{KNPS},
and some cases were verified in \cite{Collier-Li}.

\begin{thm}
\label{thm;15.7.20.100}
There exist positive constants
$t_0$, $\epsilon_0$ and $C_0$,
which may depend only on
$X$, $\phi_1,\ldots,\phi_r$ and $\gamma$,
such that the following holds
if $t\geq t_0$:
\[
\Bigl|
 \frac{1}{t}
 \vec{d}(h_{\gamma(0)},\Pi_{\gamma}^{\ast}h_{\gamma(1)})
-(2\alpha_1,\ldots,2\alpha_r)
\Bigr|
\leq
 C_0\exp(-\epsilon_0t)
\]
\end{thm}
\pf
In the following,
the constants $C_i$ and $\epsilon_i$
may depend only on 
$X$, $\phi_1,\ldots,\phi_r$ and $\gamma$,
unless otherwise specified.

We can take finite points 
$s_0=0<s_1<\cdots <s_{N-1}<s_N=1$
and coordinate neighbourhoods
$(U_k,z_k)$ around $\gamma(s_k)$
such that the following holds:
\begin{itemize}
\item
We can take $R_k>0$ $(k=0,\ldots,N)$ such that
$U_k$ contains the disc
$\Delta_k(R_k):=\{|z_k|<R_k\}$ 
and that
$\Delta_k(R_k/2):=\{|z_k|<R_k/2\}$
$(k=0,\ldots,N)$
give a covering of $\gamma([0,1])$.
\item
We have the expressions
$\phi_{i|U_k}=f_{k,i}dz_k$.
Set 
$d_k:=\min\bigl\{|f_{k,i}(0)-f_{k,j}(0)|\,\big|\,i\neq j\bigr\}$.
Then, we have
$|f_{k,i}(z_k)-f_{k,i}(0)|\leq d_k/100$
on $\Delta_k(R_k)$.
\end{itemize}
If $t_0$ is sufficiently large,
we have $t_0d_k\geq 1$ for any $k$.
Then,
Assumption \ref{assumption;15.7.17.100}
is satisfied for the harmonic bundles
$(E,\delbar_E,\theta,h)_{|\Delta_k(R_k)}$
if $t\geq t_0$.
So, we can apply the results 
in \S\ref{subsection;15.7.17.20}--\ref{subsection;15.7.17.30}
to each
$(E,\delbar_E,\theta,h)_{|\Delta_k(R_k)}$.

Let $h_{E_i}$ be the restriction of $h$ to $E_i$.
Let $\nabla_i$ be the Chern connection of
$(E_i,\delbar_{E_i},h_{E_i})$,
and we have the following connection on $E$:
\[
 \DD^1_0=\bigoplus_{i=1}^r
 \bigl(\nabla_i
+t(\phi_i+\overline{\phi_i})\id_{E_i}
 \bigr)
\]
We fix a K\"{a}hler metric $g_X$ of $X$.
By Corollary \ref{cor;15.7.17.40},
we have the following estimate
on the union of $\Delta_k(R_k/2)$
with respect to $h$ and $g_X$
for some positive constants
$C_1$ and $\epsilon_1$:
\[
 \Bigl|
 \DD^1-\DD^1_0
 \Bigr|_{h,g_X}
\leq C_1\exp(-\epsilon_1t)
\]

We have the vector bundle
$\gamma^{\ast}E
=\bigoplus_{i=1}^r\gamma^{\ast}E_i$
with the metric $\gamma^{\ast}h$
and the connections 
$\gamma^{\ast}\DD^1$
and $\gamma^{\ast}\DD^1_0$.
We take any orthonormal frames
$u_i$ of $\gamma^{\ast}E_i$
such that $\gamma^{\ast}\nabla_iu_i=0$.
They give a frame 
$\vecu=(u_1,\ldots,u_r)$ of $\gamma^{\ast}E$.
We have the following estimates
for some positive constants $C_2$ and $\epsilon_2$:
\[
\Bigl|
 \gamma^{\ast}h(u_i,u_j)
\Bigr|
\leq
C_2\exp(-\epsilon_2t),
\quad
(i\neq j).
\]
The connection
$\gamma^{\ast}\DD^1_0$ is represented 
by the diagonal matrix $A\,ds$
with respect to the frame $\vecu$,
where the $(i,i)$-entry of $A$ is 
$2t\Re a_i(s)$.
Hence, we have
\[
 \gamma^{\ast}\DD^1\vecu
=\vecu\bigl(A(s)+B_0(s)+B_1(s)\bigr)\,ds
\]
where $B_m(s)$ $(m=0,1)$ 
are $r$-square matrix valued $C^{\infty}$-functions
such that 
(i) $B_0(s)_{ij}=0$ $(i\neq j)$
and $B_1(s)_{ij}=0$ $(i=j)$,
(ii) there exist $\epsilon_3>0$ and $C_3>0$
such that 
\[
 \bigl|B_m(s)\bigr|\leq C_3\exp(-\epsilon_3t)
\quad  (m=0,1).
\]
We may assume that
$C_3\exp(-\epsilon_3t)$
is sufficiently small for any $t\geq t_0$
if $t_0$ is sufficiently large.
Then, applying Corollary \ref{cor;15.7.27.10} 
in \S\ref{subsection;15.7.17.101} below,
we have a $C^1$-function $G:[0,1]\lrarr M_r(\cnum)$
and a $C^0$-function $H:[0,1]\lrarr M_r(\cnum)$
such that the following holds
for some positive constants $C_4$ and $\epsilon_4$:
\begin{itemize}
\item
 The $C^1$-norm of $G$
 is dominated by 
 $C_4\exp(-\epsilon_4t)$.
\item
 $H(s)$  are diagonal,
 and 
 the $C^0$-norm of $H$ is dominated by 
 $C_4\exp(-\epsilon_4t)$.
\item
 The connection $\gamma^{\ast}\DD^1$ is represented by
$(A+B_0+H)\,ds$
with respect to the frame $\vecv:=\vecu(I+G)$,
where $I$ denotes the $r$-square identity matrix.
\end{itemize}
The frame $\vecv \exp\Bigl(-\int_0^s(A+B_0+H)d\tau\Bigr)$
of $\gamma^{\ast}E$ is flat with respect to
$\gamma^{\ast}\DD^1$.
The parallel transport $\Pi_{\gamma}$ of $\gamma^{\ast}\DD^1$
from $E_{\gamma(0)}$ to $E_{\gamma(1)}$
is represented by the following matrix
with respect to $\vecu(0)$  and $\vecu(1)$:
\[
 \bigl(I+G(1)\bigr)
 \exp\Bigl(
 -\int_0^1(A+B_0+H)d\tau
 \Bigr)
 \bigl(I+G(0)\bigr)^{-1}
\]
Here, $I$ denotes the $r$-square identity matrix.

Let $\vecp(s)=(p_1(s),\ldots,p_r(s))$ $(s=0,1)$
be the orthonormal frames of $E_{|\gamma(s)}$,
induced by
the frames $\vecu(s)$
and the Gram-Schmidt process.
We have 
$\vecp(s)=\vecu(s)\cdot(I+K(s))$,
where we have positive constants
$C_5$ and $\epsilon_5$
such that
$|K(s)|\leq C_5\exp(-\epsilon_5t)$.
Let $L(s)$ $(s=0,1)$ be determined by
$L(s)=(I+K(s))(I+G(s))-I$.
Then,
$\Pi_{\gamma}$ is represented 
by the following matrix
with respect to
the orthonormal bases
$\vecp(0)$ and $\vecp(1)$:
\[
 Z_{\gamma}:=
 \bigl(I+L(1)\bigr)
  \exp\Bigl(
 -\int_0^1(A+B_0+H)d\tau
 \Bigr)
 \bigl(I+L(0)\bigr)^{-1}
\]

For any $r$-square matrix $Y$,
we set
\[
 \|Y\|_{\op}:=
 \sup
\Bigl\{
 |Y\vecv|
\,\Big|\,
 \vecv\in\cnum^r,\,|\vecv|=1
 \Bigr\}.
\]
We clearly have 
$\|Y_1Y_2\|_{\op}
\leq
 \|Y_1\|_{\op}
 \|Y_2\|_{\op}$.
So, 
we have
positive constants $C_6$ and $\epsilon_6$,
such that the following holds:
\[
 \log\|Z_{\gamma}\|_{\op}
\leq
 \log
 \Bigl\|
  \exp\Bigl(
 -\int_0^1(A+B_0+H)d\tau
 \Bigr)
 \Bigr\|_{\op}
+\rho_1,
\quad\quad
 |\rho_1|\leq C_6\exp(-\epsilon_6t)
\]
We also have positive constants $C_7$ and $\epsilon_7$,
such that the following holds:
\[
 \log
 \Bigl\|
  \exp\Bigl(
 -\int_0^1(A+B_0+H)d\tau
 \Bigr)
 \Bigr\|_{\op}
\leq
  2t\alpha_1+\rho_2,
\quad\quad
 |\rho_2|\leq
 C_7\exp(-\epsilon_7t)
\]
We also have 
positive constants $C_8$ and $\epsilon_8$,
such that the following holds:
\[
 \bigl|
 \Pi_{\gamma}p_1(0)
 \bigr|_{h_{\gamma(1)}}
=\exp(2t\alpha_1)\cdot(1+\rho_3),
\quad\quad
 |\rho_3|\leq
 C_8\exp(-\epsilon_8t)
\]
Hence, we have
positive constants $C_9$ and $\epsilon_9$,
such that the following holds:
\[
  \log\|Z_{\gamma}\|_{\op}
\geq
 2t\alpha_1+\rho_4,
\quad\quad
 |\rho_4|\leq C_{9}\exp(-\epsilon_{9}t)
\]
Therefore,
we have 
positive constants $C_{10}$ and $\epsilon_{10}$,
such that the following holds:
\[
\Bigl|
 \log\|Z_{\gamma}\|_{\op}
-2t\alpha_1
\Bigr|
\leq C_{10}\exp(-\epsilon_{10}t)
\]
By applying the argument
to $\bigwedge^kZ_{\gamma}$,
we obtain positive constants $C_{11}$ and $\epsilon_{11}$,
such that the following holds
for any $k$:
\begin{equation}
 \label{eq;15.7.22.30}
\Bigl|
 \log\|
 \bigwedge^kZ_{\gamma}\|_{\op}
-2t\sum_{j=1}^k\alpha_j
\Bigr|
\leq C_{11}\exp(-\epsilon_{11}t)
\end{equation}
By using (\ref{eq;15.7.22.50}),
we can deduce the claim of the theorem
from (\ref{eq;15.7.22.30}).
\hfill\qed

\subsection{Appendix: A singular perturbation theory}
\label{subsection;15.7.17.101}

We explain a singular perturbation theory 
which is available in our situation,
and which seems slightly different from those
in \cite{Collier-Li} and \cite{KNPS}.
Note that we applied Corollary \ref{cor;15.7.27.10}
in the proof of Theorem \ref{thm;15.7.20.100}
to find a family of small gauge transforms
with which the family of connections
$\gamma^{\ast}\DDD$ on $\closedclosed{0}{1}$
are transformed to the connections
whose connection matrices are diagonal.

\subsubsection{Preliminary}

We set $[0,1]:=\bigl\{s\in\real\,\big|\,0\leq s\leq 1\bigr\}$.
For any non-negative integer $\ell$,
let $C^{\ell}([0,1])$ denote the space of 
$\cnum$-valued $C^{\ell}$-functions on $[0,1]$.
We set $\|f\|_0:=\sup_{s\in [0,1]}|f(s)|$
for any $f\in C^0([0,1])$.

Let $M_r(\cnum)$ denote the space of
$r$-square matrices.
Let $M_r(\cnum)_0$ denote the space of
the $r$-square diagonal matrices,
i.e.,
$M_r(\cnum)_0=\bigl\{
(a_{ij})\in M_r(\cnum)\,\big|\,
 a_{ij}=0\,\,(i\neq j)
\bigr\}$.
Let $M_r(\cnum)_1$ denote the set of
the off-diagonal matrices,
i.e.,
$M_r(\cnum)_1:=\bigl\{
(a_{ij})\in M_r(\cnum)\,\big|\,
 a_{ij}=0\,\,(i=j)
\bigr\}$.
We have 
$M_r(\cnum)
=M_r(\cnum)_0\oplus M_r(\cnum)_1$.

For any non-negative integer $\ell$,
let $C^{\ell}([0,1],M_r(\cnum))$
denote the space of 
$C^{\ell}$-maps $X:[0,1]\lrarr M_r(\cnum)$.
Similarly,
let $C^{\ell}([0,1],M_r(\cnum)_{\kappa})$ $(\kappa=0,1)$
be the space of 
$C^{\ell}$-maps $X:[0,1]\lrarr M_r(\cnum)_{\kappa}$.
We set $\|X\|_0:=\sup_{i,j}\|X_{ij}\|_0$ 
for $X\in C^0([0,1],M_r(\cnum))$.

\subsubsection{Statement}
\label{subsection;15.7.1.22}

Fix $C_0>0$.
We consider 
$a_j,b_j\in C^0([0,1])$ $(j=1,\ldots,r)$
satisfying the following conditions:
\begin{itemize}
\item
 $\Re a_1(s)<\Re a_2(s)<\cdots <\Re a_r(s)$
 for any $s$.
\item
 $\bigl|b_j(s)\bigr|\leq C_0$ for any $s$.
\end{itemize}
For any $t\geq 0$,
we put
$\alpha^t_j(s):=ta_j(s)+b_j(s)$.
Let $A^t$ denote the $r$-square matrix
whose $(j,j)$-entries are $\alpha^t_j$.

\begin{prop}
\label{prop;15.7.19.1}
There exist constants $C_1>0$ and $\epsilon_1>0$,
depending only on $C_0$,
such that the following holds:
\begin{itemize}
\item
For any $t\geq 0$
and  any $B\in C^0([0,1],M_r(\cnum)_1)$
satisfying $\|B\|_0\leq \epsilon_1$,
 we can take
 $G^t\in C^1\bigl([0,1],M_r(\cnum)_1\bigr)$
 and 
 $H^t\in C^0\bigl([0,1],M_r(\cnum)_0\bigr)$
 such that 
 (i) $\|G^t\|_0+\|\del_sG^t+[A^t,G^t]\|_0+\|H^t\|_0
 \leq C_1\|B\|_0$,
 (ii) we have
\begin{equation}
 \label{eq;15.7.19.2}
A^t+B=(I+G^t)^{-1}(A^t+H^t)(I+G^t)+(I+G^t)^{-1}\del_sG^t.
\end{equation}
Here, $I$ denotes the $r$-square identity matrix.
\end{itemize}
\end{prop}
We shall prove the proposition in 
\S\ref{subsection;15.7.19.20}--\ref{subsection;15.7.19.21}.
Indeed, we shall give a more refined claim 
(see Corollary \ref{cor;15.7.1.21}).

We give a reformulation of 
Proposition \ref{prop;15.7.19.1}.
Recall that when we have a vector bundle $V$ with
a connection $\nabla$ and a frame 
$\vecw=(w_1,\ldots,w_m)$,
we have the matrix-valued $1$-form 
$A=(A_{ij})$
determined by
$\nabla w_j=\sum A_{ij}w_i$,
and we describe the relation by
$\nabla\vecw=\vecw A$.

\begin{cor}
\label{cor;15.7.27.10}
Let $E$ be a $C^1$-vector bundle on $[0,1]$
with a frame $\vecv=(v_1,\ldots,v_r)$.
Let $B\in C^0([0,1],M_r(\cnum)_1)$
satisfying $\|B\|_0\leq \epsilon_1$.
Take $t\geq 0$.
Let $\nabla^t$ be the connection on $E$
given as follows:
\[
 \nabla^t\vecv=\vecv\bigl(A^t(s)+B(s)\bigr)ds.
\]
Then, we can take
$G^t\in C^1([0,1],M_r(\cnum)_1)$
and 
$H^t\in C^0([0,1],M_r(\cnum)_0)$
such that 
(i) $\|G^t\|_0+\|\del_sG^t+[A^t,G^t]\|_0+\|H^t\|_0\leq C_1\|B\|_0$,
(ii) for the frame $\vecu^t=\vecv\cdot (I+G^t)^{-1}$,
we have
$\nabla^t\vecu^t
=\vecu^t\cdot (A^t+H^t)ds$.
Here, $\epsilon_1$ and $C_1$
are constants as in Proposition {\rm\ref{prop;15.7.19.1}}.
\hfill\qed
\end{cor}

\begin{rem}
In Proposition {\rm\ref{prop;15.7.19.1}},
the off-diagonal part $B$ is required to be sufficiently small.
It is satisfied in the Hitchin WKB-problem if $t$ is sufficiently large,
as observed in Corollary {\rm\ref{cor;15.7.17.40}}.
However, 
it is not satisfied 
in the Riemann-Hilbert WKB-problem {\rm\cite{KNPS}}
even if $t$ is large, in general.
But, we may still apply Proposition {\rm\ref{prop;15.7.19.1}}
after dividing the path to shorter paths
\footnote{This was remarked by C. Simpson.}
as did in {\rm\cite{KNPS}}.
It also seems possible to apply Proposition {\rm\ref{prop;15.7.19.1}}
to the WKB-problem for family of $\lambda$-connections
without going to small paths.
\hfill\qed
\end{rem}

\subsubsection{Some linear maps}
\label{subsection;15.7.19.20}

Let $C^1([0,1],M_r(\cnum)_1)_{\del}$
denote the subspace of
$C^1([0,1],M_r(\cnum)_1)$
which consists of the functions
$X=(X_{ij}):[0,1]\lrarr M_r(\cnum)$
such that
$X_{ij}(1)=0$ $(i<j)$
and $X_{ij}(0)=0$ $(i>j)$.
We have the linear map
$D_0^t:C^1([0,1],M_r(\cnum)_1)_{\del}
\lrarr
 C^0([0,1],M_r(\cnum)_1)$
given by
\[
 D_0^t(X):=
 \del_sX+[A^t,X],
\quad
\mbox{\rm i.e., }
 D_0^t(X)_{ij}
=\del_sX_{ij}
+(\alpha_i^t-\alpha_j^t)X_{ij}.
\]

For $i\neq j$,
we put
\[
 F_{i,j}^t(s):=
\int_0^s\bigl(
 \alpha^t_i(\tau)-\alpha^t_j(\tau)
 \bigr)d\tau.
\]
We have the map
$I_0^t:C^0\bigl([0,1],M_r(\cnum)_1\bigr)
\lrarr C^1(\bigl[0,1],M_r(\cnum)_1\bigr)_{\del}$
given as follows:
\[
 I_0^t(X)_{i,j}:=
 \left\{
 \begin{array}{ll}
 \int_1^s\exp(-F_{i,j}^t(s)+F_{i,j}^t(\tau))
 X_{i,j}(\tau)d\tau
 & (i<j)\\
 \mbox{{}}\\
 \int_0^s\exp(-F_{i,j}^t(s)+F_{i,j}^t(\tau))
 X_{i,j}(\tau)d\tau
 & (i>j)
 \end{array}
 \right.
\]

Together with the identity map
on $C^0([0,1],M_r(\cnum)_0)$,
we also obtain the following maps:
\[
 I_0^t:
 C^0([0,1],M_r(\cnum)_0)
\oplus
 C^0([0,1],M_r(\cnum)_1)
\lrarr
 C^0([0,1],M_r(\cnum)_0)
\oplus
 C^1([0,1],M_r(\cnum)_1)_{\del}
\]
\[
 D_0^t:
 C^0([0,1],M_r(\cnum)_0)
\oplus
 C^1([0,1],M_r(\cnum)_1)_{\del}
\lrarr
  C^0([0,1],M_r(\cnum)_0)
\oplus
 C^0([0,1],M_r(\cnum)_1)
\]
Then, 
$I_0^t$ and $D_0^t$
are mutually inverse.

\begin{lem}
There exists a constant $K_1>0$,
depending only on $C_0$,
such that the following holds
for any $t\geq 0$ and for any 
 $(Z,W)\in 
 C^0([0,1],M_r(\cnum)_0)
 \oplus
 C^0([0,1],M_r(\cnum)_1)$:
\[
 \bigl\|
 I_0^t(Z,W)
 \bigr\|_{0}
\leq 
 K_1
 \bigl\|
 (Z,W)
 \bigr\|_0.
\]
\end{lem}
\pf
The estimate for $(i,i)$-entries are obvious
by the construction.
Let us consider the estimate for $(i,j)$-entries $(i\neq j)$.
We set
$Q^t_{i,j}(s):=
 t\int_{0}^s\bigl(
 a_i(\tau)-a_j(\tau)
 \bigr)\,d\tau$
and 
$R_{i,j}(s):=\int_0^s\bigl(
 b_i(\tau)-b_j(\tau)
 \bigr)d\tau$.
We have
$F_{i,j}^t=Q^t_{i,j}+R_{i,j}$.
Because $\Re(a_k)<\Re(a_{\ell})$ $(k<\ell)$,
we have the following for $s_1\leq s_2$:
\[
 \Re Q^t_{i,j}(s_1)-\Re Q^t_{i,j}(s_2)\geq 0\,\,(i<j),
\quad
 \Re Q^t_{i,j}(s_1)-\Re Q^t_{i,j}(s_2)\leq 0\,\,(i>j).
\]
Hence, we have a constant $K_1'$,
depending only on $C_0$,
such that 
$\Bigl|
 \exp\bigl(
 -F_{i,j}^t(s)+F_{i,j}^t(\tau)
 \bigr)
 \Bigr|
\leq K_1'$
holds in the cases
(i) $s\leq \tau$ and $i<j$,
(ii) $s\geq \tau$ and $i>j$.
Then, we obtain the following in the case $i<j$:
\[
 \bigl|
 I_0^t(X)_{i,j}
 \bigr|
\leq
 \int_s^1
 \Bigl|
 \exp\bigl(-F_{i,j}^t(s)+F_{i,j}^t(\tau)\bigr)X_{i,j}(\tau)
 \Bigr|\,d\tau
\leq
 K_1'\int_s^1|X_{i,j}(\tau)|d\tau
\leq
 K_1'\|X_{i,j}\|_0
\]
Hence, we obtain
$\|I_0^t(X)_{i,j}\|_0
\leq K_1'\|X_{i,j}\|_0$
in the case $i<j$.
Similarly, we obtain
$\|I_0^t(X)_{i,j}\|_0
\leq K_1'\|X_{i,j}\|_0$
in the case $i>j$.
\hfill\qed

\subsubsection{Proof of Proposition \ref{prop;15.7.19.1}}
\label{subsection;15.7.19.21}

We take a small number $\epsilon>0$.
We set
\[
 \nbigh_{\epsilon}:=\bigl\{
 H\in C^0([0,1],M_r(\cnum)_0)\,\big|\,
 \|H\|_0\leq\epsilon
 \bigr\},
\quad\quad
 \nbigg^t_{\epsilon}:=
 \bigl\{
 G\in  C^1([0,1],M_r(\cnum)_1)_{\del}\,\big|\,
 \|D_0^tG\|_0\leq\epsilon
 \bigr\}.
\]
Note that $\|G\|_0\leq K_1\|D_0^t G\|_0$
for any $G\in\nbigg^t_{\epsilon}$.

Let $I\in M_r(\cnum)$ denote the identity matrix.
If $\epsilon K_1<1/2$,
then $(I+G)(s)$ are invertible
for any $G\in \nbigg^t_{\epsilon}$.
So, we have the maps
$J^t:\nbigh_{\epsilon}\times\nbigg^t_{\epsilon}
\lrarr
 C^0([0,1],M_r(\cnum)_0)
\oplus
 C^0([0,1],M_r(\cnum)_1)$
given by
\[
 J^t(H,G):=(I+G)^{-1}(A^t+H)(I+G)
+(I+G)^{-1}\del_sG
-A^t.
\]
Let 
$T_{(H,G)}J^t$
denote the derivative at $(H,G)$.
We have
$T_{(0,0)}J^t(X,Y)=X+D_0^tY$.
More generally,
we have
\begin{multline}
 T_{(H,G)}J^t(X,Y)
=(I+G)^{-1}X(I+G)
+(I+G)^{-1}(A^t+H)Y
-(I+G)^{-1}Y(I+G)^{-1}(A^t+H)(I+G)
 \\
-(I+G)^{-1}Y(I+G)^{-1}\del_sG
+(I+G)^{-1}\del_sY
\end{multline}
We regard
$T_{(H,G)}J^t$
as maps
$C^0([0,1],M_r(\cnum)_0)
\oplus
 C^1([0,1],M_r(\cnum)_1)_{\del}
\lrarr
C^0([0,1],M_r(\cnum)_0)
\oplus
 C^0([0,1],M_r(\cnum)_1)$.
We obtain the following family of endomorphisms
on 
$C^0([0,1],M_r(\cnum)_0)
\oplus
 C^0([0,1],M_r(\cnum)_1)$:
\[
 T_{(H,G)}J^t\circ
 \bigl(
 T_{(0,0)}J^t
 \bigr)^{-1},
\quad\quad
((H,G)\in \nbigh_{\epsilon}\times\nbigg^t_{\epsilon}).
\]

\begin{lem}
\label{lem;15.8.10.1}
There exists a constant $C_2>0$,
depending only on $C_0$,
such that the following holds
for any $t\geq 0$
and any 
$(Z,W)\in C^0([0,1],M_r(\cnum)_0)
\oplus
 C^0([0,1],M_r(\cnum)_1)$:
\[
 \Bigl\|
 T_{(H,G)}J^t
\circ
\bigl(
 T_{(0,0)}J^t\bigr)^{-1}(Z,W)
-(Z,W)
\Bigr\|_0
\leq
C_2\bigl(\|H\|_0+\|D_0^tG\|_0\bigr)
 \bigl(
 \|Z\|_0+\|W\|_0
 \bigr)
\]
\end{lem}
\pf
We consider 
\begin{equation}
\label{eq;15.6.30.10}
 (I+G)^{-1}(A^t)Y
-(I+G)^{-1}Y(I+G)^{-1}(A^t)(I+G)
-(I+G)^{-1}Y(I+G)^{-1}\del_sG+(I+G)^{-1}\del_sY
-D_0^tY.
\end{equation}
We have
$(I+G)^{-1}\del_sY
=(I+G)^{-1}D_0^tY
-(I+G)^{-1}[A^t,Y]$.
Hence, (\ref{eq;15.6.30.10})
is rewritten as
\begin{multline}
 (I+G)^{-1}(A^t)Y
-(I+G)^{-1}Y(I+G)^{-1}(A^t)(I+G)
-(I+G)^{-1}Y(I+G)^{-1}\del_sG
 \\
-(I+G)^{-1}[A^t,Y]
+\bigl(
 (I+G)^{-1}-I
 \bigr)D_0^tY
\end{multline}
It is equal to the following:
\begin{multline}
\label{eq;15.7.1.20}
 (I+G)^{-1}Y(A^t)
-(I+G)^{-1}Y(I+G)^{-1}(A^t)(I+G)
-(I+G)^{-1}Y(I+G)^{-1}\del_sG
 \\
+\bigl(
 (I+G)^{-1}-I
 \bigr)D_0^tY
\end{multline}
We have the following:
\begin{multline}
-(I+G)^{-1}Y(I+G)^{-1}(A^t)(I+G)
-(I+G)^{-1}Y(I+G)^{-1}\del_sG
 \\
=-(I+G)^{-1}Y(I+G)^{-1}(A^t)
-(I+G)^{-1}Y(I+G)^{-1}
 \bigl(\del_sG+(A^t)G\bigr)
 \\
=-(I+G)^{-1}Y(I+G)^{-1}(A^t)
-(I+G)^{-1}Y(I+G)^{-1}G(A^t)
-(I+G)^{-1}Y(I+G)^{-1}D_0^tG
\end{multline}
Hence, (\ref{eq;15.7.1.20})
is equal to the following:
\begin{multline}
 (I+G)^{-1}Y(A^t)
-(I+G)^{-1}Y(I+G)^{-1}(A^t)
 -(I+G)^{-1}Y(I+G)^{-1}G(A^t)
 \\
-(I+G)^{-1}Y(I+G)^{-1}D_0^tG
+\bigl(
 (I+G)^{-1}-I
 \bigr)D_0^tY
\\
=
-(I+G)^{-1}Y(I+G)^{-1}D_0^tG
+\bigl(
 (I+G)^{-1}-I
 \bigr)D_0^tY
\end{multline}
In all, we obtain the following:
\begin{multline}
 T_{(H,G)}J^t\circ
 \bigl(T_{(0,0)}J^t\bigr)^{-1}(Z,W)
-(Z,W)
=
-(I+G)^{-1}I_0^t(W)(I+G)^{-1}D_0^tG
-\bigl((I+G)^{-1}-I\bigr)W \\
+(I+G)^{-1}Z(I+G)-Z
+\bigl[
 (I+G)^{-1}H(I+G),(I+G)^{-1}I_0^t(W)
 \bigr]
\end{multline}
Because
$\|G\|_0\leq K_1\|D_0^tG\|_0$,
we have
$\|(I+G)^{-1}-I\|_0\leq 
 C_3\|D_0^tG\|_0$
and 
$\|(I+G)^{-1}X(I+G)-X\|_0
\leq C_4(\|D_0^tG\|_0\|X\|_0)$
for positive constants $C_i$ $(i=3,4)$
depending only on $C_0$.
We also have the following
for a positive constant $C_5$
depending only on $C_0$:
\[
\Bigl\|
 \bigl[
 (I+G)^{-1}H(I+G),
 (I+G)^{-1}Y
 \bigr]
\Bigr\|_0
\leq
C_5\bigl(\|D_0^tG\|_0+\|H\|_0\bigr)\|D_0^tY\|_0
\]
Hence, we obtain the claim of the lemma.
\hfill\qed

\vspace{.1in}

\begin{cor}
\label{cor;15.7.1.21}
There exist positive constants
$\epsilon_{10}>0$ and $C_{10}>0$,
depending only on $C_0$,
with the following property:
\begin{itemize}
\item
For any  $t\geq 0$
and any $B\in C^0([0,1],M_r(\cnum)_1)$
such that 
$\|B\|_0\leq \epsilon_{10}$,
we have
a unique 
$(H^t,G^t)\in \nbigh_{\epsilon}\times\nbigg^t_{\epsilon}$
such that 
$J^t(H^t,G^t)=B$
and that
$\|H^t\|_0+\|D_0^tG^t\|_0\leq C_{10}\|B\|_0$.
\end{itemize}
\end{cor}
\pf
We set 
$\nbiggbar_{\epsilon}:=
 \bigl\{
 \Gbar\in C^0([0,1],M_r(\cnum)_1)\,\big|\,
 \|\Gbar\|_0\leq\epsilon
 \bigr\}
\subset
 C^0([0,1],M_r(\cnum)_1)$.
We have the bijections
$I_0^t:\nbiggbar_{\epsilon}\simeq
 \nbigg_{\epsilon}^t$ for any $t\geq 0$.
We consider the maps
\[
 \nbigf^t:=
 J^t\circ
 (T_{(0,0)}J^t)^{-1}:
 \nbigh_{\epsilon}\times
 \nbiggbar_{\epsilon}
\lrarr
 C^0([0,1],M_r(\cnum)_0)
\oplus
 C^0([0,1],M_r(\cnum)_1).
\]
Let $T\nbigf^t:
 C^0([0,1],M_r(\cnum))
\lrarr 
 C^0([0,1],M_r(\cnum))$ 
denote the derivative of $\nbigf^t$.
According to Lemma \ref{lem;15.8.10.1},
we have a positive constant $C_3$ such that 
the operator norms of 
$T_{(H,\Gbar)}\nbigf^t-\id$
are dominated by 
$C_3\bigl(
 \|H\|_0+\|\Gbar\|_0\bigr)$
for any $(H,\Gbar)\in\nbigh_{\epsilon}\times\nbiggbar_{\epsilon}$.
The constant $C_3$ may depend only on $C_0$,
and the estimate is uniform for $t$.
By the inverse function theorem
(see \cite{Lang}, for instance),
there exist positive constants
$\epsilon_{10}$ and $C_{10}$,
depending only on $C_0$,
with the following property:
\begin{itemize}
\item
For any  $t\geq 0$
and any $B\in C^0([0,1],M_r(\cnum)_1)$
such that 
$\|B\|_0\leq \epsilon_{10}$,
we have a unique
$(H^t,\Gbar^t)\in \nbigh_{\epsilon}\times\nbiggbar_{\epsilon}$
such that 
$\nbigf^t(H^t,\Gbar^t)=B$
and that
$\|H^t\|_0+\|\Gbar^t\|_0\leq C_{10}\|B\|_0$.
\end{itemize}
By setting $G^t:=(T_{(0,0)}J^t)^{-1}\Gbar^t$,
we obtain the claim of Corollary \ref{cor;15.7.1.21}.
We also finish the proof of Proposition \ref{prop;15.7.19.1}.
\hfill\qed

\subsection{Appendix: The case of Hermitian-Einstein metrics}
\label{subsection;16.5.23.30}

For $R>0$,
we set $\Delta(R):=
 \bigl\{z\in\cnum\,\big|\,|z|<R\bigr\}$.
Let $(E,\delbar_E,\theta)$ be a Higgs bundle of rank $r$
on $\Delta(R)$.
We fix a Hermitian metric $h_{\det(E)}$ of $\det(E)$.
Let $h$ be a Hermitian-Einstein metric of 
$(E,\delbar_E,\theta)$
such that $\det(h)=h_{\det(E)}$,
i.e.,
$R(h)^{\bot}+[\theta,\theta^{\dagger}_h]=0$,
where $R(h)^{\bot}$ denotes the trace-free part of 
the curvature $R(h)$.
We have an obvious generalization of the results in 
\S\ref{subsection;15.7.17.20}
and \S\ref{subsection;15.7.17.30},
which we state explicitly in this subsection.

We have a $C^{\infty}$-function $\nu$ such that
$\delbar\del\nu=R(h_{\det(E)})/r$
such that the $L_2^p$-norm of $\nu$
is dominated by $C_{0,p}\|R(h_{\det(E)})\|_{\infty}$ $(p>1)$,
where $\|R(h_{\det(E)})\|_{\infty}$ denotes the sup norm
of $R(h_{\det(E)})$ with respect to the Euclidean metric,
and $C_{0,p}$ denotes the constant depending on $p$
and the radius $R$.
Then, the metric $\htilde=he^{-\nu}$ is 
a harmonic metric for 
the Higgs bundle $(E,\delbar_E,\theta)$.
Note that
$h$ and $\htilde$ induce the same metrics
on the vector bundles $\End(E)\otimes\Omega^{p,q}$.

We have the description $\theta=f\,dz$.
Let $M$ be the constant as in \S\ref{subsection;15.7.23.1}.
We obtain the following from 
Proposition \ref{prop;15.7.17.10}.
\begin{cor}
Fix $0<R_1<R$.
Then, we have $C_1,C_2>0$ depending only on
$r, R_1,R$ such that
$|f|_h\leq C_1M+C_2$
holds on $\Delta(R_1)$.
\hfill\qed
\end{cor}
We impose the conditions
as in Assumption \ref{assumption;15.7.17.100}.
We obtain the following from Proposition \ref{prop;15.7.12.10}.
\begin{cor}
Let $R_1$ be as in Proposition {\rm\ref{prop;15.7.17.10}}.
Fix $0<R_2<R_1$.
There exist positive constants 
$\epsilon_0$
and $C_{11}$
depending only on
$R$, $R_1$, $R_2$, $r$ and $C_{10}$,
such that 
$|\rho_{\alpha}|_h\leq 
 C_{11}\exp(-\epsilon_0d)$
on $\Delta(R_2)$.
\hfill\qed
\end{cor}
We impose the condition that $\rank E_{\alpha}=1$
as in \S\ref{subsection;15.7.17.30}.
We have $R(\htilde)=R(h)^{\bot}$
and $\theta^{\dagger}_h=\theta^{\dagger}_{\htilde}$.
Hence,
we obtain the following from Theorem \ref{thm;15.7.12.11}.
\begin{cor}
There exist positive constants
$C_{30}$ and $\epsilon_{30}$,
depending only on 
$R$, $R_1$, $R_2$, $r$ and $C_{10}$ 
such that
$\bigl|
 R(h)^{\bot}
\bigr|_{h,g_{\cnum}}
=\bigl|
[\theta,\theta_h^{\dagger}]
 \bigr|_{h,g_{\cnum}}
\leq C_{30}
 \exp(-\epsilon_{30}d)$
on $\Delta(R_2)$.
\hfill\qed
\end{cor}

Because $h$ and $\htilde$ give the same connection
on $\End(E)$ as the Chern connections,
we obtain the following from Proposition \ref{prop;15.7.1.10}
and Proposition \ref{prop;15.7.1.11}.

\begin{cor}
Take $0<R_3<R_2$.
There exist positive constants $\epsilon_{40}$ and $C_{40}$
depending only on 
$R$, $R_1$, $R_2$, $R_3$, $r$ and $C_{10}$
such that 
$\bigl|\del_{E,h}\pi_{\alpha}\bigr|_{h,g_{\cnum}}=
\bigl|\delbar_E\pi^{\dagger}_{\alpha}\bigr|_{h,g_{\cnum}}
\leq
 C_{40} \exp(-\epsilon_{40}d)$
and 
$\bigl|
 \delbar_E\pi_{\alpha}'
 \bigr|_{h,g_{\cnum}}
=\bigl|
 \del_{E,h}\pi_{\alpha}'
 \bigr|_{h,g_{\cnum}}
\leq
 C_{40}
 \exp(-\epsilon_{40}d)$
on $\Delta(R_3)$.
\hfill\qed
\end{cor}

Let $h_{\alpha}$ (resp. $h_{0,\alpha}$)
denote the restriction of $h$ (resp. $\htilde$)
to $L_{\alpha}$.
Because $h_{\alpha}=\htilde_{\alpha}e^{\nu/r}$,
we obtain the following from Proposition
\ref{prop;15.7.12.12}.
\begin{cor}
We have a positive constant $C_{60}$
depending only on 
$R$, $R_1$, $R_2$, $R_3$, $r$ and $C_{10}$
such that
$\bigl|
 R(h_{\alpha})-R(h_{\det(E)})/r
 \bigr|_{h_{\alpha},g_{\cnum}}
\leq
 C_{60}\exp(-\epsilon_{40}d)$
on $\Delta(R_3)$.
\hfill\qed
\end{cor}

We have the projectively flat connection
$\DD^1:=\delbar_E+\del_{E,h}+\theta+\theta_h^{\dagger}$
induced by $h$.
We also have the flat connection
$\DDtilde^1:=
 \delbar_E+\del_{E,\htilde}
+\theta+\theta_{\htilde}^{\dagger}$.
They are related as
$\DD^1=\DDtilde^1+\del\nu\cdot\id_E$.
We have the Chern connections
$\delbar_{\alpha}+\del_{\alpha,h}$
on $E_{\alpha}$
induced by $h_{\alpha}$.
Similarly, we have the Chern connections
$\delbar_{\alpha}+\del_{\alpha,\htilde}$
on $E_{\alpha}$
induced by $\htilde_{\alpha}$.
They are related as
$\del_{\alpha,h}=
 \del_{\alpha,\htilde}+\del\nu\cdot\id_{E_{\alpha}}$.
Hence, we obtain the following from 
Corollary \ref{cor;15.7.17.40}.
\begin{cor}
We set 
$\DD^1_0:=
 \bigoplus\bigl(
 \delbar_{\alpha}+\del_{\alpha,h}
+(g_{\alpha}dz+\overline{g_{\alpha}}d\zbar)
 \cdot \id_{E_{\alpha}}
 \bigr)$
as in {\rm\S\ref{subsection;16.5.23.10}}.
Then, there exists a positive constant $C_{70}$,
depending only on 
$R$, $R_1$, $R_2$, $R_3$, $r$ and $C_{10}$
such that 
$\bigl|
 \DD^1-\DD^1_0
 \bigr|_{h,g_{\cnum}}
\leq
 C_{70}
 \exp(-\epsilon_{40}d)$
holds
on $\Delta(R_3)$.
\hfill\qed
\end{cor}

We also have an obvious generalization
of the estimates of the higher derivatives 
in \S\ref{subsection;16.5.23.20},
which we omit to describe.

\section{Local models}
\label{section;15.7.15.1}

\subsection{Review on unramifiedly good filtered Higgs bundles}

\subsubsection{Filtered bundles on curves}
\label{subsection;15.7.20.100}

Let $X$ be a complex curve 
with a discrete subset $D$.
We recall the concept of filtered bundles
on $(X,D)$ \cite{Simpson90}.
Let $\real^D$ denote the set of maps
$D\lrarr \real$.
Elements of $\real^D$
are denoted by
$\veca=(a_P|\,P\in D)$.
Let $\nbigo_X(\ast D)$ be the sheaf of meromorphic functions
on $X$ whose poles are contained in $D$.
Let $\nbige$ be a locally free $\nbigo_X(\ast D)$-module
of rank $r$.
A filtered bundle over $\nbige$
is a family of coherent $\nbigo_X$-submodules
$\nbigp_{\ast}\nbige=
 \bigl(
 \nbigp_{\veca}\nbige\subset\nbige\,\big|\,
 \veca\in \real^D \bigr)$
with the following property.
\begin{itemize}
\item
 $\nbigp_{\veca}\nbige$ 
 are lattices of $\nbige$,
 i.e.,
 $\nbigp_{\veca}\nbige\otimes_{\nbigo_X}\nbigo_X(\ast D)
=\nbige$.
\item
 The stalk of $\nbigp_{\veca}\nbige$ at $P\in D$
 depends only on $a_P\in\real$.
 We denote it by
 $\nbigp^P_{a_P}(\nbige_P)$.
\item
 We have
 $\nbigp^P_a(\nbige_P)\subset\nbigp^P_b(\nbige_P)$
 for $a\leq b$.
 For any $a\in\real$,
 there exists $\epsilon>0$
 such that 
 $\nbigp^P_a(\nbige_P)=
 \nbigp^P_{a+\epsilon}(\nbige_P)$.
\item
 For $n\in\seisuu$,
 we have
 $\nbigp^P_a(\nbige_P)
 \otimes_{\nbigo_{X,P}}\nbigo_X(nP)_P
=\nbigp^P_{a+n}(\nbige_P)$. 
\end{itemize}
Such $\nbigp_{\ast}\nbige$ is called
a filtered bundle over $\nbige$ on $(X,D)$.
For any $\nbigo_X(\ast D)$-submodule 
$\nbigg\subset\nbige$,
we have the induced filtered bundle
$\nbigp_{\ast}\nbigg$ over $\nbigg$
given by
$\nbigp_{\veca}\nbigg
=\nbigp_{\veca}\nbige\cap \nbigg$
in $\nbige$.

\paragraph{Parabolic structure}

Let $\nbigp_{\ast}\nbige$ be a filtered bundle on $(X,D)$.
Let $\veczero\in\real^D$ denote the element such that 
the $P$-th component of $\veczero$ are $0$
for any $P\in D$.
We have the locally free $\nbigo_X$-module
$\nbigp_{\veczero}\nbige$ on $X$.
For each $Q\in X$, let $\nbigp_{\veczero}\nbige_{|Q}$
denote the fiber of 
the vector bundle $\nbigp_{\veczero}\nbige$ over $Q$.
For $P\in D$ and for $-1<a\leq 0$,
we set 
\[
F^P_a(\nbigp_{\veczero}\nbige_{|P})
:=\Image\Bigl(
 \nbigp^P_{a}(\nbige_P)
\lrarr
 \nbigp_{\veczero}\nbige_{|P}
 \Bigr).
\]
In this way,
we obtain a filtration $F^P$ of 
$\nbigp_{\veczero}\nbige_{|P}$
indexed by $\{-1<a\leq 0\}$.
We have
$F^P_{0}(\nbigp_{\veczero}\nbige_{|P})=\nbigp_{\veczero}\nbige_{|P}$.
Note that 
we have the natural identification
$\nbigp_{\veczero}\nbige_{|P}
\simeq
 \nbigp^{P}_{0}(\nbige_P)
 \big/
 \nbigp^{P}_{-1}(\nbige_P)$
because 
$\nbigp^P_{a+n}(\nbige)
=\nbigp^P_{a}(\nbige)\otimes_{\nbigo_{X,P}}
 \nbigo_{X,P}(nP)$,
and hence we have
$F^P_{-1}(\nbigp_{\veczero}\nbige_{|P})=0$.
For any $-1<a\leq 0$,
there exists $\epsilon>0$
such that
$F^P_{a}(\nbigp_{\veczero}\nbige_{|P})=
 F^P_{a+\epsilon}\nbigp_{\veczero}\nbige_{|P}$.
Such a family of filtrations
$(F^P\,|\,P\in D)$
is called a parabolic structure on the vector bundle
$\nbigp_{\veczero}\nbige$ along $D$.
We can easily observe that
filtered bundles on $(X,D)$
are equivalent to
vector bundles on $X$ with a parabolic structure along $D$.

\subsubsection{Unramifiedly good filtered Higgs bundles}

Let $X$ and $D$ be as in \S\ref{subsection;15.7.20.100}.
Let $\nbigp_{\ast}\nbige$ be a filtered bundle on $(X,D)$.
Let $\theta$ be a Higgs field of $\nbige$,
i.e., $\theta$ is a holomorphic section of
$\End(\nbige)\otimes\Omega^1_X$.
The filtered bundle with a Higgs field 
$(\nbigp_{\ast}\nbige,\theta)$
is called unramifiedly good at $P$,
if the following holds.
\begin{itemize}
\item
 We have a finite subset 
 $\nbigi(P)\subset \nbigo_X(\ast D)_P$
 and a decomposition of the stalk
$\nbigp^P_a\nbige_P=
 \bigoplus_{\gminia\in \nbigi(P)}
 \nbigp^P_a\nbige_{P,\gminia}$
such that 
\[
 (\theta-d\gminia\id)
 \nbigp^P_a\nbige_{P,\gminia}
\subset
 \nbigp^P_a \nbige_{P,\gminia}
 \otimes\Omega^1_X(\log D)_P.
\]
Here, $d\gminia$ denotes the exterior derivative
of $\gminia\in\nbigi(P)\subset\nbigo_X(\ast D)_P$.
We have
$\nbigp^P_a\nbige_{P,\gminia}
\subset \nbigp^P_b\nbige_{P,\gminia}$
for $a\leq b$. 
The induced map
$\nbigi(P)\lrarr
 \nbigo_X(\ast D)_P/\nbigo_{X,P}$
is assumed to be injective.
\end{itemize}
If $(\nbigp_{\ast}\nbige,\theta)$ is unramifiedly good at any $P\in D$,
it is called an unramifiedly good filtered Higgs bundle.

\begin{rem}
More generally,
a filtered bundle with a Higgs field
$(\nbigp_{\ast}\nbige,\theta)$ is called good at $P$
if we have a neighbourhood $X_P$ of $P$,
a Galois covering map $X'_P\lrarr X_P$ ramified along $P$,
and an unramifiedly good filtered Higgs bundle 
$(\nbigp^P_{\ast}\nbige',\theta')$ on $(X'_P,\pi^{-1}(P))$
such that
$(\nbigp_{\ast}\nbige,\theta)_{|X_P}$
is the descent of 
$(\nbigp^P_{\ast}\nbige',\theta')$.
In this paper, we consider only unramified ones.
\hfill\qed
\end{rem}

\subsubsection{Unramified wild harmonic bundles and 
the associated unramifiedly good filtered Higgs bundles}

Let $X$ and $D$ be as in \S\ref{subsection;15.7.20.100}.
Let $(E,\delbar_E,\theta,h)$ be a harmonic bundle on $X\setminus D$.
It is called wild and unramified over $(X,D)$
if the following holds
for any $P\in D$.
\begin{itemize}
\item
 Let $(U_P,z)$ be a holomorphic coordinate
 neighbourhood of $P$ in $X$ with $z(P)=0$.
 Then, we have a finite subset 
 $\nbigi(P)\subset z^{-1}\cnum[z^{-1}]$
 and a decomposition
\[
 (E,\delbar_E,\theta)_{|U_P\setminus\{P\}}=
 \bigoplus_{\gminia\in\nbigi(P)}
 (E_{\gminia},\delbar_{E_{\gminia}},\theta_{\gminia}),
\]
and the coefficients 
$a_{\gminia,j}(z)$
of the characteristic polynomials
 $\det(t\id-g_{\gminia})
=\sum a_{\gminia,j}(z)t^j$
 are holomorphic at $z=0$,
where $g_{\gminia}\in\End(E_{\gminia})$ 
are determined by
$\theta_{\gminia}-d\gminia\id_{E_{\gminia}}
=g_{\gminia}\,dz/z$.
\end{itemize}
A unramified wild harmonic bundle
is called tame
if $\nbigi(P)=\{0\}$ in $\nbigo_X(\ast D)_P/\nbigo_{X,P}$
at each $P\in D$.

The tame case of the following proposition
was established in \cite{Simpson90},
and generalized to the wild case 
in \cite{Mochizuki-wild}.

\begin{prop}
Let $(E,\delbar_E,\theta,h)$
be an unramified wild harmonic bundle
on $(X,D)$.
Then, we have an unramifiedly good filtered Higgs bundle
$(\nbigp^h_{\ast}E,\theta)$ on $(X,D)$
with the following property.
\begin{itemize}
\item
 We have
 $(\nbigp^h_{\veca}E,\theta)_{|X\setminus D}
 =(E,\theta)$
 for any $\veca\in\real^D$.
\item
 Let $U$ be any open subset of $X$.
 Let $f$ be a holomorphic section of $E$
 on $U\setminus D$.
 Then, $f$ is a section of
 $\nbigp^h_{\veca}E$ on $U$
 if and only if 
$|f|_h=O(|z_P|^{-a_P-\epsilon})$ for any $\epsilon>0$
 around each $P\in U\cap D$,
 where $z_P$ denotes a holomorphic coordinate 
 around $P$ with $z_P(P)=0$.
\hfill\qed
\end{itemize}
\end{prop}

In this situation,
$h$ is called adapted to
the filtered bundle $\nbigp^h_{\ast}E$.

\subsubsection{Filtered Higgs bundles on compact Riemann surfaces}

Suppose that
$X$ is a compact connected Riemann surface
with a finite subset $D$.
Let $(\nbigp_{\ast}\nbige,\theta)$ be 
an unramifiedly good filtered Higgs bundle on $(X,D)$.
Let $\nbigg\subset\nbige$ be any locally free
$\nbigo_X(\ast D)$-submodule.
We set
$\nbigp_{\veca}\nbigg
:=\nbigp_{\veca}\nbige\cap\nbigg$
in $\nbige$ for any $\veca\in\real^D$,
and then we obtain a filtered bundle 
$\nbigp_{\ast}\nbigg$ over $\nbigg$.
The degree of $\nbigp_{\ast}\nbigg$
is defined as follows
\cite{Maruyama-Yokogawa,Mehta-Seshadri}
for any $\vecb\in\real^D$:
\begin{equation}
 \label{eq;15.8.18.2}
 \deg(\nbigp_{\ast}\nbigg)
=\deg(\nbigp_{\vecb}\nbigg)
-\sum_{P\in D}
 \sum_{b_P-1<a\leq b_P}
 a\cdot
 \dim_{\cnum}\bigl(
 \nbigp^P_{a}(\nbigg_P)
 \big/
 \nbigp^P_{<a}(\nbigg_P)
 \bigr)
\end{equation}
Here,
we set
 $\nbigp^P_{<a}(\nbigg_{P})
=\bigcup_{c<a}
 \nbigp^P_{c}(\nbigg_{P})$,
and we regard
$\nbigp^P_{a}(\nbigg_P)
 \big/
 \nbigp^P_{<a}(\nbigg_P)$
as finite dimensional $\cnum$-vector spaces
in a natural way.
It is easy to see that the right hand side of
(\ref{eq;15.8.18.2}) is independent of the choice of $\vecb$.
The unramifiedly good filtered Higgs bundle
$(\nbigp_{\ast}\nbige,\theta)$ is called stable
(resp. semistable)
if the following inequality holds
for any locally free $\nbigo_X(\ast D)$-submodule
$\nbigg\subset\nbige$
with 
(i) $\theta\nbigg\subset
     \nbigg\otimes\Omega^1_X$,
(ii) $0<\rank\nbigg<\rank \nbige$:
\[
 \frac{\deg(\nbigp_{\ast}\nbigg)}{\rank\nbigg}
<
\frac{\deg(\nbigp_{\ast}\nbige)}{\rank\nbige}
\quad\quad
\left(
\mbox{\rm resp. }
 \frac{\deg(\nbigp_{\ast}\nbigg)}{\rank\nbigg}
\leq
\frac{\deg(\nbigp_{\ast}\nbige)}{\rank\nbige}
\right)
\]
The unramifiedly good filtered Higgs bundle
$(\nbigp_{\ast}\nbige,\theta)$ 
on $(X,D)$ is called polystable
if it is the direct sum
 of stable ones
$(\nbigp_{\ast}\nbige_i,\theta_i)$ $(i=1,\ldots,m)$
with $\deg(\nbigp_{\ast}\nbige_i)=\deg(\nbigp_{\ast}\nbige)$.
In the following proposition, the tame case was proved 
by Simpson \cite{Simpson90},
and see \cite{Mochizuki-wild} for the wild case,
for example.

\begin{prop}
Let $(E,\delbar_E,\theta,h)$ be an unramifiedly good
wild harmonic bundle on $(X,D)$.
Then, the associated filtered bundle
$(\nbigp^h_{\ast}E,\theta)$ on $(X,D)$
is poly-stable with $\deg(\nbigp^h_{\ast}E)=0$.
\hfill\qed
\end{prop}

In the following proposition, the tame case was proved 
by Simpson \cite{Simpson90},
and the wild case was proved 
by Biquard and Boalch \cite{biquard-boalch}.
\begin{prop}
\label{prop;15.7.23.10}
Let $(\nbigp_{\ast}\nbige,\theta)$
be an unramifiedly good filtered Higgs bundle
on $(X,D)$ with $\deg(\nbigp_{\ast}\nbige)=0$.
If $(\nbigp_{\ast}\nbige,\theta)$ is stable,
we have a harmonic metric $h$
of $(E,\theta)=(\nbige,\theta)_{|X\setminus D}$
with an isomorphism
$(\nbigp^h_{\ast}E,\theta)
\simeq
 (\nbigp_{\ast}\nbige,\theta)$.
Such a metric $h$ is unique up to 
the multiplication of positive constants.
\hfill\qed
\end{prop}

\subsubsection{Filtered bundles of rank one on curves}
\label{subsection;15.7.21.10}

Let $X$ be any complex curve
with a discrete subset $D$.
Let $L$ be a line bundle on $X$.
Suppose that 
a tuple of real numbers
$\vecb=(b_P\,|\,P\in D)$
is attached.
Then, we have the filtered bundle
$\nbigp^{\vecb}_{\ast}L$ over 
the meromorphic bundle $L(\ast D)$
given as follows.
For $\veca=(a_P\,|\,P\in D)\in\real^{D}$,
we have 
the integers $n(a_P)$ $(P\in D)$
such that
$a_P-1<n(a_P)+b_P\leq a_P$,
and we set 
\[
 \nbigp^{\vecb}_{\veca}L
:=L\Bigl(
 \sum_{P\in D} n(a_P)P
 \Bigr).
\]
If $X$ is compact,
we have
$\deg(\nbigp_{\ast}^{\vecb}L)
=\deg(L)-\sum_{P\in D} b_{P}$.

\subsubsection{Comparison of filtered bundles}

Let $X$ be a compact connected Riemann surface
with a finite subset $D$.
Let $\nbigp_{\ast}\nbige$ be a filtered bundle 
over a locally free $\nbigo_X(\ast D)$-module $\nbige$
on $(X,D)$.
Fix a point $P\in D$.
Set $D_1:=D\setminus P$.
For any $b\in\real$,
let $\nbigp^P_b\nbige$ denote 
the $\nbigo_{X}(\ast D_1)$-module
such that
$\nbigp^P_b\nbige(\ast P)=\nbige$
and
the stalk of $\nbigp^P_b\nbige$ at $P$
is $\nbigp^P_b(\nbige_P)$.
For any $b\in \real$
and $\vecc\in\real^{D_1}$,
we set
$\nbigp_{\vecc}(\nbigp^P_b\nbige):=
 \nbigp_{(b,\vecc)}\nbige$.
For any $b\in\real$,
we have the filtered bundle
$\nbigp_{\ast}(\nbigp^P_b\nbige)
=\bigl(
 \nbigp_{\vecc}\nbigp^P_b\nbige\,\big|\,
 \vecc\in \real^{D_1}
 \bigr)$
over 
$\nbigp^P_b\nbige$
on $(X,D_1)$.

We have the set
$\bigl\{a\,\big|\,0<a\leq 1,\,\,
 \nbigp^P_a(\nbige_P)/\nbigp^{P}_{<a}(\nbige_P)
\neq 0
\bigr\}
=\{a_1,\ldots,a_{\ell}\}$.
We assume $a_i<a_{i+1}$.
The following observation is given
in \cite{Maruyama-Yokogawa}.

\begin{lem}
\label{lem;15.8.7.10}
We have the following equality:
\begin{equation}
\label{eq;15.8.7.1}
 \int_0^1
 \deg\bigl(\nbigp_{\ast}(\nbigp^P_b\nbige)\bigr)\,db
=\deg\bigl(\nbigp_{\ast}(\nbigp^P_0\nbige)\bigr)
-\sum_{i=1}^{\ell}
 (a_i-1)\dim\bigl(
\nbigp^P_{a_i}(\nbige_P)/\nbigp^P_{<a_i}(\nbige_P)
 \bigr)
\end{equation}
Here, 
we regard
$\deg\bigl(\nbigp_{\ast}(\nbigp^P_b\nbige)\bigr)$
as a measurable function in variable $b$.
In particular,
we have
$\deg(\nbigp_{\ast}\nbige)
=\int_{0}^1
 \deg\bigl(\nbigp_{\ast}(\nbigp^P_b\nbige)\bigr)\,db$.
\end{lem}
\pf
We obtain (\ref{eq;15.8.7.1}) by a direct computation.
Because 
$\nbigp^P_{a_i}(\nbige_P)/\nbigp^P_{<a_i}(\nbige_P)
\simeq
\nbigp^P_{a_i-1}(\nbige_P)/\nbigp^P_{<a_i-1}(\nbige_P)$,
we can observe the right hand side of
(\ref{eq;15.8.7.1})
is equal to $\deg(\nbigp_{\ast}\nbige)$.
\hfill\qed

\vspace{.1in}

We recall a lemma for the comparison
of filtered bundles
for the convenience of the readers.
\begin{lem}
\label{lem;15.8.7.11}
Let $\nbige$ be a locally free $\nbigo_X(\ast D)$-module.
Let $\nbigp_{\ast}^i\nbige$ $(i=1,2)$ be
filtered bundles over $\nbige$
such that
$\nbigp^1_{\veca}\nbige\subset\nbigp^2_{\veca}\nbige$
for any $\veca\in\real^D$.
If $\deg(\nbigp^1_{\ast}\nbige)=\deg(\nbigp^2_{\ast}\nbige)$,
then we have
$\nbigp^1_{\veca}\nbige=\nbigp^2_{\veca}\nbige$ 
for any $\veca\in\real^D$.
\end{lem}
\pf
We use an induction on $|D|$.
If $D=\emptyset$,
the claim is clear.
Take $P\in D$,
and set $D_1:=D\setminus \{P\}$.
We obtain the filtered bundles
$\nbigp^i_{\ast}\nbigp^P_{b}\nbige$ $(i=1,2)$
for any $b\in\real$ as above.
By using Lemma \ref{lem;15.8.7.10},
we can easily deduce 
$\deg\bigl(\nbigp^1_{\ast}(\nbigp^P_b\nbige)\bigr)
=\deg\bigl(\nbigp^2_{\ast}(\nbigp^P_b\nbige)\bigr)$
for any $b\in\real$
from the equality
$\deg(\nbigp^1_{\ast}\nbige)
=\deg(\nbigp^2_{\ast}\nbige)$.
By the hypothesis of the induction,
we obtain 
$\nbigp^1_{\vecc}(\nbigp^P_b\nbige)
=\nbigp^2_{\vecc}(\nbigp^P_b\nbige)$
for any $\vecc$.
Then, the claim of the lemma follows.
\hfill\qed

\subsection{Filtered Higgs bundles of rank $2$
on $(\proj^1,\infty)$}
\label{subsection;15.7.28.1}

Let $\zeta$ be the standard coordinate 
on $\cnum\subset\proj^1$.
We set
$\Vtilde=
 \nbigo_{\proj^1}(\ast\infty)\vtilde_1
\oplus
 \nbigo_{\proj^1}(\ast\infty)\vtilde_2$.
Take a non-zero complex number $\alpha$
and a positive integer $m$.
We have the Higgs field
$\thetatilde$ on $\Vtilde$ 
given by
\[
 \thetatilde(\vtilde_1,\vtilde_2)
=(\vtilde_1,\vtilde_2)
 \left(
 \begin{array}{cc}
 \alpha\zeta^md\zeta & 0 \\
 0 & -\alpha\zeta^m d\zeta
 \end{array}
 \right).
\]

Fix $0\leq \ell\leq m$.
Let $\Etilde_{\ell}\subset \Vtilde$ be 
the $\nbigo_{\proj^1}(\ast\infty)$-submodule
generated by $\etilde_1=\vtilde_1+\vtilde_2$
and $\etilde_2=\zeta^{\ell}\vtilde_2$.
Because
$\thetatilde(\Etilde_{\ell})
\subset
 \Etilde_{\ell}\otimes\Omega^1_{\proj^1}$,
we obtain the meromorphic Higgs bundle
$(\Etilde_{\ell},\thetatilde)$.
We have
\[
 \thetatilde(\etilde_1,\etilde_2)
=(\etilde_1,\etilde_2)
 \left(
 \begin{array}{cc}
 \alpha \zeta^md\zeta
 & 0 \\
 -2\alpha\zeta^{m-\ell}d\zeta
 & -\alpha\zeta^md\zeta
 \end{array}
 \right).
\]

Take $c\in\real$,
and set $c_1:=c$ and $c_2:=-c-\ell$.
We have the parabolic Higgs bundle
$(\nbigp^{c}_{\ast}\Etilde_{\ell},\thetatilde)$
on $(\proj^1,\infty)$
with $\deg(\nbigp^{c}_{\ast}\Etilde_{\ell})=0$
given as follows
\[
 \nbigp^c_{b}(\Etilde_{\ell})_{|\proj^1\setminus\{\infty\}}
=\Bigl(
 \nbigo_{\proj^1}([b-c_1]\infty)\vtilde_1
\oplus
 \nbigo_{\proj^1}([b-c_2]\infty)\vtilde_2
 \Bigr)_{|\proj^1\setminus\{\infty\}}
\]
Here, $[x]:=\max\{n\in\seisuu\,|\,n\leq x\}$.
We shall impose that $c_1\geq c_2$,
i.e., $c\geq -\ell/2$.

\subsubsection{Stability condition}

We set $L_i:=\nbigo_{\proj^1}(\ast\infty)\vtilde_i\subset \Vtilde$
$(i=1,2)$.
We have
$\Etilde_{\ell}\cap L_i=
 \nbigo_{\proj^1}(\ast\infty)\cdot \zeta^{\ell}\vtilde_i$.
We have the induced parabolic bundle
$\nbigp^c_{\ast}(\Etilde_{\ell}\cap L_i)$.
We have
$\deg\bigl(\nbigp^c_{\ast}(\Etilde_{\ell}\cap L_i)\bigr)
=-\ell-c_i$.
Hence,
$(\nbigp^{c}_{\ast}\Etilde_{\ell},\thetatilde)$ is stable
if and only if
$-\ell-c_i<0$ $(i=1,2)$,
i.e.,
$c_i<0$ $(i=1,2)$.
We also have that
$(\nbigp^{c}_{\ast}\Etilde_{\ell},\thetatilde)$
is semistable if and only if
$-\ell-c_1=0$ or $-\ell-c_2=0$,
i.e.,
$c_1=0,-\ell$.
Under the assumption
$c_1\geq c_2$,
the semistability is equivalent to
$c_1=0$.

\subsubsection{The determinant bundles}

We have
$\det(\Etilde_{\ell})=
 \nbigo_{\proj^1}(\ast\infty)
 \cdot \etilde_1\wedge \etilde_2
=\nbigo_{\proj^1}(\ast\infty)
 \cdot \zeta^{\ell}\vtilde_1\wedge\vtilde_2$.
The induced filtered bundle
$\det(\nbigp^{c}_{\ast}\Etilde_{\ell},\thetatilde)$
over $\det(\Etilde_{\ell})$ 
is equal to
$\nbigp^0_{\ast}
 \bigl(
 \nbigo_{\proj^1}\cdot
 \etilde_1\wedge\etilde_2
 \bigr)$,
where
$\nbigp^0_{b}
 \bigl(
 \nbigo_{\proj^1}\cdot
 \etilde_1\wedge\etilde_2
 \bigr)
=\nbigo_{\proj^1}([b])\cdot
 \etilde_1\wedge\etilde_2$
for $[b]:=\max\{n\in\seisuu\,|\,n\leq b\}$.
In particular, 
they are independent of $c$.
We fix the Hermitian metric 
$h_{\det(\Etilde_{\ell})}$ of $\det(\Etilde_{\ell})$
given by 
$h_{\det(\Etilde_{\ell})}\bigl(
 \etilde_1\wedge \etilde_2,\etilde_1\wedge \etilde_2
\bigr)=1$.

\subsubsection{Harmonic metrics in the case $\ell=0$}
\label{subsection;15.7.28.10}

If $\ell=0$,
$(\nbigp_{\ast}^c\Etilde_0,\thetatilde)$ cannot be stable.
And, 
$(\nbigp_{\ast}^c\Etilde_0,\thetatilde)$ is semistable
if and only if $c=0$.
Indeed, it is polystable in that case,
i.e.,
we have the decomposition
$(\nbigp_{\ast}^0\Etilde_0,\thetatilde)
=(\nbigp_{\ast}^0L_1,\thetatilde_1)
\oplus
 (\nbigp_{\ast}^0L_2,\thetatilde_2)$,
where $\thetatilde_1$
and $\thetatilde_2$
are the multiplications of $\alpha\zeta^md\zeta$
and $-\alpha\zeta^md\zeta$,
respectively.
We have the harmonic metrics 
$h_{L_i}$ for $(\nbigp_{\ast}^0L_i,\thetatilde_i)$.
We impose 
$h_{L_1}\otimes h_{L_2}=h_{\det(\Etilde_0)}$.
Then, we have a harmonic metric
$h_{\Etilde_0}=h_{L_1}\oplus h_{L_2}$
for $(\nbigp_{\ast}^0\Etilde_0,\thetatilde)$.
Note that we have an ambiguity given by automorphisms
of $(\nbigp_{\ast}^0\Etilde_0,\thetatilde)$,
i.e.,
for any $\alpha>0$,
$\alpha h_{L_1}\oplus \alpha^{-1}h_{L_2}$
is also a harmonic metric 
for $(\nbigp_{\ast}^0\Etilde_0,\thetatilde)$.
We also note that 
$h_{L_i}(\vtilde_i,\vtilde_i)$
are constants.

\subsubsection{Harmonic metrics in the case $\ell>0$
 and their homogeneous property}
\label{subsection;15.8.6.12}

Suppose that $\ell>0$.
According to Proposition \ref{prop;15.7.23.10},
if the unramifiedly good filtered Higgs bundle
$(\nbigp^{c}_{\ast}\Etilde_{\ell},\thetatilde)$ is stable,
i.e.,
$-\ell<c<0$,
we have the harmonic metric $h_{c,\ell}$
of $(\Etilde_{\ell},\thetatilde)_{|\cnum}$
adapted to 
the filtered Higgs bundle 
$(\nbigp^c_{\ast}\Etilde_{\ell},\thetatilde)$
such that 
$\det(h_{c,\ell})=h_{\det \Etilde_{\ell}}$.

\vspace{.1in}

For any $\tau\in\cnum^{\ast}$,
let $\varphi_{\tau}:\proj^1\lrarr\proj^1$
be given by
$\varphi_{\tau}(\zeta)=\tau^2\zeta$.
We have the isomorphism
$\Phi_{\tau}:
 \varphi_{\tau}^{\ast}\Etilde_{\ell}
\simeq
 \Etilde_{\ell}$
given by
$\tau^{\ell}\varphi_{\tau}^{\ast}\etilde_1
\longleftrightarrow
 \etilde_1$
and 
$\tau^{-\ell}\varphi_{\tau}^{\ast}\etilde_2
\longleftrightarrow
 \etilde_2$.
Under the isomorphism,
we have
$\tau^{\ell}\varphi_{\tau}^{\ast}\vtilde_i
\longleftrightarrow
 \vtilde_i$ $(i=1,2)$.
Hence, $\Phi_{\tau}$ induces an isomorphism of
the filtered bundles
$\varphi_{\tau}^{\ast}\nbigp_{\ast}\Etilde_{\ell}
\simeq
 \nbigp_{\ast}\Etilde_{\ell}$.
Under the isomorphism,
we have
$\varphi_{\tau}^{\ast}\thetatilde
=\tau^{2(m+1)}\thetatilde$.
Hence,
the Hermitian metric
$\varphi_{\tau}^{\ast}h_{c,\ell}$
gives a harmonic metric of
$(\nbigp^c_{\ast}\Etilde_{\ell},\tau^{2(m+1)}\thetatilde)$.
Note that
$\det\varphi^{\ast}h_{c,\ell}
=h_{\det \Etilde_{\ell}}$.

\begin{prop}
If $|\tau|=1$,
we have
$\varphi_{\tau}^{\ast}h_{c,\ell}
=h_{c,\ell}$.
\end{prop}
\pf
If $|\tau|=1$,
both the Hermitian metrics
$\varphi_{\tau}^{\ast}h_{c,\ell}$
and 
$h_{c,\ell}$
are harmonic metrics
of $(\Etilde_{\ell},\thetatilde)$.
They are adapted to the filtered bundle
$\nbigp^c_{\ast}\Etilde_{\ell}$,
and they satisfy 
$\det(h_{c,\ell})=\det(\varphi_{\tau}^{\ast}h_{c,\ell})
=h_{\det \Etilde_{\ell}}$.
Hence,
we have
$h_{c,\ell}=\varphi_{\tau}^{\ast}h_{c,\ell}$
by the uniqueness.
\hfill\qed

\begin{cor}
\label{cor;15.8.6.1}
The functions
$h_{c,\ell}(\etilde_i,\etilde_i)$
$(i=1,2)$
depend only on $|\zeta|$.
The function
$\zeta^{-\ell}h_{c,\ell}(\etilde_1,\etilde_2)$
depend only on $|\zeta|$.
The functions
$h_{c,\ell}(\vtilde_i,\vtilde_j)$ $(i,j\in\{1,2\})$
depend only on $|\zeta|$.
\hfill\qed
\end{cor}

\begin{rem}
We clearly have the homogeneous property 
of the harmonic metrics
such as Corollary {\rm\ref{cor;15.8.6.1}}
even in the case $\ell=0$.
\hfill\qed
\end{rem}

\subsubsection{The norms of $\vtilde_j$}
\label{subsection;15.7.13.50}

By the norm estimate of wild harmonic bundles
\cite{Mochizuki-wild},
we have constants 
$C_{i,c}>0$ $(i=1,2)$
depending on $c$
such that 
$C_{1,c}|\zeta|^{c_j}
\leq
 \bigl|
 \vtilde_j
 \bigr|_{h_{c,\ell}}
\leq
 C_{2,c}|\zeta|^{c_j}$,
where $c_1=c$ and $c_2=-c-\ell$.
Let us refine it in our situation.

\begin{prop}
\label{prop;15.8.6.20}
There exist positive constants
$b_c$, $C_{3}$ and $\epsilon_{3}$
such that the following holds
on $\{|\zeta|>1\}$:
\begin{equation}
\label{eq;15.8.6.10}
 \Bigl|
 \log
 |\vtilde_1|_{h_{c,\ell}}
-\log\bigl(
 b_c|\zeta|^{c_1}
 \bigr)
 \Bigr|
\leq
 C_{3}
 \exp(-\epsilon_{3}|\zeta|^{m+1})
\end{equation}
\begin{equation}
\label{eq;15.8.6.11}
 \Bigl|
 \log
 |\vtilde_2|_{h_{c,\ell}}
-\log\bigl(
 b^{-1}_c|\zeta|^{c_2}
 \bigr)
 \Bigr|
\leq
 C_{3}
 \exp(-\epsilon_{3}|\zeta|^{m+1})
\end{equation}
\begin{equation}
\label{eq;15.8.6.30}
\Bigl|
 \zeta\del_{\zeta}\log|\vtilde_i|_{h_{c,\ell}}
-c_i/2
\Bigr|
 \leq
 C_{3}\exp\bigl(-\epsilon_3|\zeta|^{m+1}\bigr)
\quad
(i=1,2)
\end{equation}
Here, $b_c$ may depend on $c$ ,
but $C_{3}$ and $\epsilon_{3}$
are independent of $c$.
\end{prop}
\pf
Let $\eta:=\zeta^{-1}$ be the coordinate around $\infty$.
Set $U_{\infty}=\{|\eta|\leq 1\}\subset\proj^1$.

\begin{lem}
\label{lem;15.7.16.1}
There exist $C_4>0$ and $\epsilon_4>0$
which are independent of $c$,
such that the following holds
on $U_{\infty}\setminus\{\infty\}$:
\[
 \Bigl|
 \del_{\eta}
 \del_{\etabar}\bigl(
 \log |\vtilde_j|_{h_{c,\ell}}
 \bigr)
\Bigr|
\leq
 C_4\exp\bigl(-\epsilon_4|\eta|^{-m-1}
 \bigr)
\]
\end{lem}
\pf
Let $h_{L_j,c,\ell}$ denote the restriction of
$h_{c,\ell}$ to $L_j$.
Let $R(h_{L_j,c,\ell})$ denote the curvature
of the Chern connection of
$(L_j,h_{L_j,c,\ell})$.
Let $g_{\cnum}$ be the Euclidean metric
$d\zeta\,d\zetabar$.
By Proposition \ref{prop;15.7.12.12},
we have positive constants $C_5$ and $\epsilon_5$
which are independent of $c$,
such that the following holds:
\[
 \bigl|
 R(h_{L_j,c,\ell})
 \bigr|_{h_{L_j,c,\ell},g_{\cnum}}
\leq
 C_5
 \exp\bigl(-\epsilon_5|\zeta|^{m+1}\bigr)
\]
Because 
$\delbar\del\log\bigl|
 \vtilde_j
 \bigr|_{h_{c,\ell}}
=R(h_{L_j,c,\ell})$,
we obtain the claim of the lemma.
\hfill\qed

\vspace{.1in}

We set $Y_{\etatilde}:=
 \bigl\{\etatilde\in\cnum\,\big|\,
 |\etatilde|<1
 \bigr\}$.
For any $\kappa<1$,
we have the isomorphism
$\Psi_{\kappa}:
 Y_{\etatilde}
\lrarr 
 \{|\eta|<\kappa\}$
given by
$\Psi_{\kappa}(\etatilde)=\kappa\etatilde$.
We have the following on $Y_{\etatilde}\setminus\{0\}$:
\[
 \bigl|
 \del_{\etatilde}
 \del_{\overline{\etatilde}}
 \Psi_{\kappa}^{\ast}(\log|\vtilde_1|_{h_{c,\ell}})
 \bigr|
\leq
 C_4
 \kappa^{2}
\exp(-\epsilon_4\kappa^{-m-1}|\etatilde|^{-m-1})
\]

Take a large $p>1$.
For each $c$ and $\kappa$,
we can take an $\real$-valued $L_2^p$-function
$G_{\kappa,c}$
on a neighbourhood of the closure of $Y_{\etatilde}$
such that 
(i) $G_{\kappa,c}$ is a function of $|\eta|$,
(ii) there exist positive constants
 $C_{6}$ and $\epsilon_{6}$
 such that 
the $L_2^p$-norm of $G_{\kappa,c}$ on $Y_{\etatilde}$
is dominated by
$C_{6}\exp(-\epsilon_{6}\kappa^{-m-1})$,
(iii) $G_{\kappa,c}(0)=0$,
(iv) the following holds on
$Y_{\etatilde}$:
\[
 \del_{\etatilde}
 \del_{\overline{\etatilde}}
 \bigl(
  \Psi_{\kappa}^{\ast}(\log|\vtilde_1|_{h_{c,\ell}})
 -G_{\kappa,c}
 \bigr)
=0
\]
Because
$\Psi_{\kappa}^{\ast}(\log|\vtilde_1|_{h_{c,\ell}})
 -G_{\kappa,c}
-\log|\etatilde|^{-c_1}$
is bounded,
we have a holomorphic function 
$g_{\kappa,c}$ on 
$Y_{\etatilde}$
such that
$\Psi_{\kappa}^{\ast}(\log|\vtilde_1|_{h_{c,\ell}})
 -G_{\kappa,c}
-\log|\etatilde|^{-c_1}
=\Re(g_{\kappa,c})$.
Because 
$\Psi_{\kappa}^{\ast}(\log|\vtilde_1|_{h_{c,\ell}})
 -G_{\kappa,c}
-\log|\etatilde|^{-c_1}$
depends only on $|\eta|$,
we obtain that
$g_{\kappa,c}$ is constant.
We also obtain
\[
 \etatilde\del_{\etatilde}
 \Psi_{\kappa}^{\ast}(\log|\vtilde_1|_{h_{c,\ell}})
+c_1/2
-\etatilde_{\del_{\etatilde}}
 G_{\kappa,c}=0.
\]

Hence,
on $\{|\eta|<\kappa\}$,
the function
$\log|\vtilde_1|_{h_{c,\ell}}
-(\Psi_{\kappa}^{-1})^{\ast}G_{\kappa,c}
-\log|\eta|^{-c_1}$
is a constant.
We have positive constants
$C_7$ and $\epsilon_7$ such that 
the following holds 
on $\{\kappa/2<|\eta|<\kappa\}$:
\[
 \bigl|
 (\Psi_{\kappa}^{-1})^{\ast}G_{\kappa,c}
 \bigr|
\leq
 C_7\exp\bigl(-\epsilon_7|\eta|^{-m-1}\bigr)
\]
\[
 \bigl|
 \eta\del_{\eta}
  (\Psi_{\kappa}^{-1})^{\ast}G_{\kappa,c}
 \bigr|
\leq
  C_7\exp\bigl(-\epsilon_7|\eta|^{-m-1}\bigr)
\]

We also obtain that the function
$F:=\log|\vtilde_1|_{h_{c,\ell}}
-\log|\eta|^{-c_1}$ gives a continuous function
on $\{|\eta|<1\}$.
We set $b_c:=\exp(F(0))$.
Then, we obtain (\ref{eq;15.8.6.10}),
and (\ref{eq;15.8.6.30}) for $i=1$.
By using 
$\bigl|
 \vtilde_1\wedge\vtilde_2
 \bigr|_{h_{c,\ell}}=|\zeta|^{-\ell}$
and the estimate (\ref{eq;15.7.13.1}) below,
we obtain the estimate
(\ref{eq;15.8.6.11}).
We also obtain (\ref{eq;15.8.6.30}) for $i=2$
with a similar argument.
\hfill\qed

\subsubsection{Asymptotic orthogonality and some complements}

As studied in \cite{Mochizuki-wild},
we have the following estimate
on $\{|\zeta|\geq 1\}$,
which also follows from 
Corollary \ref{cor;15.8.18.22}:
\begin{equation}
 \label{eq;15.7.13.1}
\bigl|
 h_{c,\ell}(\vtilde_1,\vtilde_2)
\bigr|
\leq
 K\exp(-\delta|\zeta|^{m+1})
 \cdot
 \bigl|\vtilde_1\bigr|_{h_{c,\ell}}
\cdot
 \bigl|\vtilde_2\bigr|_{h_{c,\ell}}
\leq 
 K'\exp(-\delta|\zeta|^{m+1})
\end{equation}
Here, $K$, $K'$ and $\delta$
are positive constants
which are independent of $c$.

Let $\delbar_{c,\ell}$ denote the holomorphic structure
of $\Etilde_{\ell}$.
Let $\del_{c,\ell}$ denote the $(1,0)$-part of
the Chern connection of $(\Etilde_{\ell},h_{c,\ell})$.
Let $g_{\cnum}$ denote the Euclidean metric of 
$\cnum_{\zeta}$.

\begin{lem}
\label{lem;15.8.18.31}
We have positive constants $K_1$
which are independent of $c$,
such that 
$\bigl|
 \zeta\del_{c,\ell}\vtilde_i
 \bigr|_{h_{c,\ell},g_{\cnum}}
\leq
 K_1\cdot |\vtilde_i|_{h_{c,\ell}}$
on $\{|\zeta|\geq 1\}$.
\end{lem}
\pf
In this proof,
the constants $K_i$ and $\delta_i$ are independent of $c$.
To simplify the description,
we denote $h_{c,\ell}$ by $h$.
By the asymptotic orthogonality (\ref{eq;15.7.13.1}),
we have a positive constant $K_2$
such that the following holds
on $\{|\zeta|\geq 1\}$:
\[
 \bigl|
 \del_{c,\ell}\vtilde_1
 \bigr|_{h,g_{\cnum}}
\leq
 K_2\Bigl(
 |\vtilde_1|_{h}^{-1}\cdot
 \bigl|
 h(\del_{c,\ell}\vtilde_1,\vtilde_1)
 \bigr|_{g_{\cnum}}
+|\vtilde_2|_{h}^{-1}\cdot
 \bigl|
 h(\del_{c,\ell}\vtilde_1,\vtilde_2)
 \bigr|_{g_{\cnum}}
 \Bigr)
\]
By Proposition \ref{prop;15.8.6.20},
we have a constant $K_3>0$
such that the following holds
on $\{|\zeta|\geq 1\}$:
\[
 |\vtilde_1|_{h}^{-1}
 \bigl|
 h(\del_{c,\ell}\vtilde_1,\vtilde_1)
 \bigr|_{g_{\cnum}}
=
 |\vtilde_1|_{h}\cdot
 \Bigl|
 \del\log|\vtilde_1|^2_{h}
 \Bigr|_{g_{\cnum}}
\leq
 K_3
 |\vtilde_1|_h\cdot
 |\zeta|^{-1}
\]
Let $\pi_1$ denote the projection of 
$\Vtilde=L_1\oplus L_2$ onto $L_1$.
We also have the following on 
$\{|\zeta|\geq 1\}$:
\begin{equation}
\label{eq;15.8.18.30}
 |\vtilde_2|_h^{-1}\cdot
 \bigl|
 h(\del_{c,\ell}\vtilde_1,\vtilde_2)
 \bigr|_{g_{\cnum}}
\leq
 |\vtilde_2|_h^{-1}\cdot
 \bigl|
 h\bigl((\del_{c,\ell}\pi_1)\cdot \vtilde_1,\vtilde_2\bigr)
 \bigr|
+
 |\vtilde_2|_h^{-1}\cdot
 \bigl|
 h\bigl(\pi_1(\del_{c,\ell}\vtilde_1),\vtilde_2\bigr)
 \bigr|
\end{equation}
By Proposition \ref{prop;15.7.1.10},
we have positive constants $K_4$ and $\delta_4$
such that
\begin{equation}
 \label{eq;15.8.18.32}
  \bigl|\del_{c,\ell}\pi_1\bigr|_{h,g_{\cnum}}
\leq
 K_4\exp(-\delta_4|\zeta|^{m+1}).
\end{equation}
By (\ref{eq;15.7.13.1}) and Proposition \ref{prop;15.7.12.10},
we have positive constants $K_5$ and $\delta_5$
such that
\begin{equation}
 \label{eq;15.8.18.33}
 \bigl|
 h\bigl(\pi_1(\del_{c,\ell}\vtilde_1),\vtilde_2\bigr)
 \bigr|
\leq
 K_5\exp(-\delta_5|\zeta|^{m+1})
 \bigl|\del_{c,\ell}\vtilde_1\bigr|_{h,g_{\cnum}}
\cdot
 \bigl|\vtilde_2\bigr|_{h}.
\end{equation}
Hence, 
we obtain
the following on $\{|\zeta|\geq 1\}$:
\[
 \bigl|\del_{c,\ell}\vtilde_1\bigr|_{h,g_{\cnum}}
\leq
 K_2\bigl(
 K_3|\zeta|^{-1}
+K_4\exp(-\delta_4|\zeta|^{m+1})
 \bigr)|\vtilde_1|_{h}
+K_2K_5\exp(-\delta_5|\zeta|^{m+1})
  \bigl|\del_{c,\ell}\vtilde_1\bigr|_{h,g_{\cnum}}
\]
Hence, we obtain the desired estimate for 
$\del_{c,\ell}\vtilde_1$.
Similarly, we obtain the estimate for
$\del_{c,\ell}\vtilde_2$.
\hfill\qed

\begin{lem}
\label{lem;15.8.18.40}
We have positive constants $K_6$ and $\delta_6$
which are independent of $c$,
such that the following holds:
\[
 \Bigl|
 \del_{\zeta}h_{c,\ell}(\vtilde_1,\vtilde_2)
 \Bigr|
\leq
 K_6\exp(-\delta_6|\zeta|^{m+1}),
\quad\quad
 \Bigl|
 \del_{\zetabar}h_{c,\ell}(\vtilde_1,\vtilde_2)
 \Bigr|
\leq
 K_6\exp(-\delta_6|\zeta|^{m+1}),
\]
\[
 \Bigl|
 \del_{\zeta}\del_{\zetabar}
 h_{c,\ell}(\vtilde_1,\vtilde_2)
 \Bigr|
\leq
  K_6\exp(-\delta_6|\zeta|^{m+1}).
\]
\end{lem}
\pf
In this proof,
the constants $K_i$ and $\delta_i$ are independent of $c$.
We obtain the estimate for 
$\del_{\zeta}h_{c,\ell}(\vtilde_1,\vtilde_2)$
from (\ref{eq;15.8.18.30}--\ref{eq;15.8.18.33})
and Lemma \ref{lem;15.8.18.31}.
We obtain the estimate for 
$\del_{\zetabar}h_{c,\ell}(\vtilde_1,\vtilde_2)$
in a similar way.
Let $R(h_{c,\ell})$ denote the curvature of
$(\Etilde_{\ell},h_{c,\ell})$.
We have the following:
\[
 \del\delbar
 h_{c,\ell}(\vtilde_1,\vtilde_2)
=h_{c,\ell}(\del_{c,\ell}\vtilde_1,\del_{c,\ell}\vtilde_2)
+h_{c,\ell}(\vtilde_1,R(h_{c,\ell})\vtilde_2)
\]
We have positive constants
$K_7$ and $\delta_7$ such that
$|R(h_{c,\ell})|_{h_{c,\ell},g_{\cnum}}
\leq
 K_7\exp(-\delta_n|\zeta|^{m+1})$.
Hence, 
we have
\[
 \bigl|
h_{c,\ell}(\vtilde_1,R(h_{c,\ell})\vtilde_2)
 \bigr|_{g_{\cnum}}
\leq
 K_8\exp(-\delta_8|\zeta|^{m+1})
\]
for positive constants $K_8$ and $\delta_8$.
Let $\pi_i$ denote the projection of
$\Vtilde=L_1\oplus L_2$ to $L_i$.
We have the following:
\begin{multline}
 h_{c,\ell}(\del_{c,\ell}\vtilde_1,\del_{c,\ell}\vtilde_2)
=h_{c,\ell}\bigl(
 \pi_1(\del_{c,\ell}\vtilde_1),
 \pi_2(\del_{c,\ell}\vtilde_2)
 \bigr)
+h_{c,\ell}\bigl(
 (\del_{c,\ell}\pi_1)\vtilde_1,
 \pi_2(\del_{c,\ell}\vtilde_2)
 \bigr)
 \\
+h_{c,\ell}\bigl(
 \pi_1(\del_{c,\ell}\vtilde_1),
 (\del_{c,\ell}\pi_2)\vtilde_2
 \bigr)
+h_{c,\ell}\bigl(
  (\del_{c,\ell}\pi_1)\vtilde_1,
 (\del_{c,\ell}\pi_2)\vtilde_2
 \bigr)
\end{multline}
Hence, as in the case of the estimate for
$h_{c,\ell}(\del_{c,\ell}\vtilde_1,\vtilde_2)$,
we have positive constants $K_9$ and $\delta_9$
such that 
$
\bigl|
 h_{c,\ell}(\del_{c,\ell}\vtilde_1,\del_{c,\ell}\vtilde_2)
\bigr|_{g_{\cnum}}
\leq
 K_9\exp(-\delta_9|\zeta|^{m+1})$.
Thus, we obtain the desired estimate for
$\del_{\zeta}\del_{\zetabar}h_{c,\ell}(\vtilde_1,\vtilde_2)$.
\hfill\qed

\subsubsection{Convergence of some sequences}

Set $t:=\tau^{2(m+1)}$ for $\tau\in\cnum^{\ast}$.
We use the notation in \S\ref{subsection;15.8.6.12}.
We have the family 
of the harmonic metrics
$h_{t,c,\ell}:=
 \varphi_{\tau}^{\ast}h_{c,\ell}$
for $(\nbigp^c_{\ast}\Etilde_{\ell},t\thetatilde)$
satisfying 
$\det(h_{t,c,\ell})=h_{\det(E)}$.

Under the isomorphism
$\Phi_{\tau}:
 \varphi_{\tau}^{\ast}\Etilde_{\ell}
\simeq
 \Etilde_{\ell}$,
we have
$\tau^{\ell}\varphi_{\tau}^{\ast}\vtilde_i
\longleftrightarrow
 \vtilde_i$ $(i=1,2)$.
Take any $T>0$.
We have positive constants
$C$ and $\epsilon$,
which are independent of $c$ and $T$,
such that the following holds
on $\{|\zeta|\geq T\}$
for any $t$ satisfying $|t|>T^{-1}$:
\[
\Bigl|
\log\bigl|\vtilde_1\bigr|_{h_{t,c,\ell}}
-\log\bigl(
b_c|t|^{(\ell+2c)/2(m+1)}
 |\zeta|^{c} 
\bigr)
\Bigr|
\leq 
C
\exp\bigl(-\epsilon|\zeta|^{m+1}|t|\bigr)
\]
\[
 \Bigl|
\log
 \bigl|\vtilde_2 \bigr|_{h_{t,c,\ell}}
-\log\bigl(
 b_c^{-1}|t|^{-(\ell+2c)/2(m+1)}
 |\zeta|^{-c-\ell}
 \bigr)
\Bigr|
\leq
 C
\exp\bigl(
 -\epsilon
 |\zeta|^{m+1}|t|
 \bigr)
\]
\[
 \Bigl|
 h_{t,c,\ell}
 \bigl(\vtilde_1,\vtilde_2\bigr)
\Bigr|
\leq
 C\exp\bigl(
 -\epsilon|\zeta|^{m+1}|t|
 \bigr)
\cdot
 \bigl|
 \vtilde_{1}
 \bigr|_{h_{t,c,\ell}}
\cdot
 \bigl|
 \vtilde_{2}
 \bigr|_{h_{t,c,\ell}}
\]
\[
 \bigl|
 R(h_{t,c,\ell})
 \bigr|_{h_{t,c,\ell},g_{\cnum}}
\leq
 C
 \exp\bigl(-\epsilon|\zeta|^{m+1}|t|\bigr)
\]
Here,
$g_{\cnum}$ denote the standard Euclidean metric
on $\cnum$.

Let $h^{\lim}_{c,\ell}$ be the Hermitian metric
of $\Vtilde_{|\cnum^{\ast}}$
given by
\[
 h^{\lim}_{c,\ell}(\vtilde_1,\vtilde_1)=|\zeta|^{2c},
\quad
 h^{\lim}_{c,\ell}(\vtilde_2,\vtilde_2)=|\zeta|^{-2c-2\ell},
\quad
 h^{\lim}_{c,\ell}(\vtilde_1,\vtilde_2)=0.
\]
For any $\gamma>0$,
the automorphism $\Psi_{\gamma}$ on $\Vtilde_{|\cnum^{\ast}}$
is given by
$\Psi_{\gamma}=\gamma\id_{L_1}\oplus \gamma^{-1}\id_{L_2}$.
We define
$\Psi_{\gamma}^{\ast}h_{t,c,\ell}(u_1,u_2)
:=h_{t,c,\ell}(\Psi_{\gamma}u_1,\Psi_{\gamma}u_2)$.
The following is clear.

\begin{prop}
\label{prop;15.8.6.100}
Set $\gamma(t):=b_c^{-1}t^{-(\ell+2c)/2(m+1)}$.
Then, 
we have the convergence
$\lim_{|t|\to\infty}
 \Psi_{\gamma(t)}^{\ast}h_{t,c,\ell}=h^{\lim}_{c,\ell}$
on $\cnum^{\ast}$.
For any fixed $T>0$,
we have a constant $\delta_{T}>0$
depending on $T$,
such that the order of the convergence 
is $\exp(-\delta_{T}|t|)$
on $\{|\zeta|\geq T\}$.
\hfill\qed
\end{prop}

\subsection{Family of harmonic metrics in the case $\ell>0$}
\label{subsection;15.7.28.2}

We continue to use the notation in \S\ref{subsection;15.7.28.1}.
We study the dependence of the harmonic metrics
$h_{c,\ell}$ on $c$ in the case $\ell>0$.

\subsubsection{Continuity of $b_{\vecc}$
with respect to the parabolic weights}

\begin{prop}
\label{prop;15.7.13.41}
$b_c$ is continuous with respect to $c$.
\end{prop}
\pf
Fix $-\ell<c^0<0$.
It is enough to study the continuity  at $c^0$.
We give a preliminary.
Take a neighbourhood $U$ of $c^0$.
We have a family of 
Hermitian metrics $h^0_{c,\ell}$ $(c\in U)$
of $\Etilde_{\ell}$
satisfying the following conditions:
\begin{itemize}
\item
We have $h^0_{c^0,\ell}=h_{c^0,\ell}$
on $\cnum$.
\item
We have 
$h^0_{c,\ell}(\vtilde_1,\vtilde_2)
=h_{c^0,\ell}(\vtilde_1,\vtilde_2)$,
$h^0_{c,\ell}(\vtilde_1,\vtilde_1)=
 h_{c^0,\ell}(\vtilde_1,\vtilde_1)
 |\zeta|^{c-c^0}$
and 
$h^0_{c,\ell}(\vtilde_2,\vtilde_2)=
 h_{c^0,\ell}(\vtilde_2,\vtilde_2)
 |\zeta|^{c^0-c}$
on $\{|\zeta|\geq 1\}$.
\item
 $\bigl|\vtilde_1\wedge\vtilde_2\bigr|_{h^0_{c,\ell}}
=|\zeta|^{\ell}$.
\item
 We have
 $\lim_{c\to c^0}h^0_{c,\ell}=h^0_{c^0,\ell}$
 in the $C^{\infty}$-sense 
 on any compact subset in $\cnum$.
\end{itemize}
Note that the conditions
are compatible on 
$\bigl\{|\zeta|\geq 1
 \bigr\}$.
We fix a $C^{\infty}$-K\"ahler metric 
$g_{\proj^1}$ of $\proj^1$.

\begin{lem}
We have a constant $C>0$
such that the following holds
on $\proj^1\setminus\{\infty\}$
for any $c\in U$:
\[
 \bigl|
 R(h^0_{c,\ell})
 \bigr|_{h^0_{c,\ell},g_{\proj^1}}
\leq C,
\quad
 \bigl|
 [\thetatilde,\thetatilde^{\dagger}_{h^0_{c,\ell}}]
 \bigr|_{h^0_{c,\ell},g_{\proj^1}}
\leq C
\]
\end{lem}
\pf
In the proof,
$C_i$ and $\epsilon_i$
are positive constants,
which are independent of $c$.
It is enough to consider the issue
on $\bigl\{|\zeta|\geq 1\bigr\}$.
The estimate for 
$[\thetatilde,\thetatilde^{\dagger}_{h^0_{c,\ell}}]$
follows from (\ref{eq;15.7.13.1})
and the construction of $h^0_{c,\ell}$.
Let us study the estimate for $R(h^0_{c,\ell})$.

Let $H$ be the $M_2(\cnum)$-valued function
on $\{|\zeta|\geq 1\}$
given by
$H_{ij}=h_{c^0,\ell}(\vtilde_i,\vtilde_j)$.
By Lemma \ref{lem;15.8.18.40},
we have positive constants
$C_{1}$ and $\epsilon_{1}$
such that the following holds:
\[
 \bigl|
 \del H_{12}
 \bigr|_{g_{\proj^1}}
\leq 
 C_{1}\exp(-\epsilon_{1}|\eta|^{-m-1}),
\quad\quad
 \bigl|
 \del H_{21}
 \bigr|_{g_{\proj^1}}
\leq 
 C_{1}\exp(-\epsilon_{1}|\eta|^{-m-1}),
\]
Set $c_1^0=c^0$ and $c_2^0=-c^0-\ell$.
By Lemma \ref{lem;15.8.18.31},
we have a positive constant $C_{2}$
such that the following holds:
\[
 \bigl|
 \del H_{ii}
 \bigr|_{g_{\proj^1}}
 \leq
 C_{2}|\eta|^{-2c^0_i-1}
\,\,\,(i=1,2)
\]

Let $\Gamma_c$ be the $M_2(\cnum)$-valued function
given as 
$\Gamma_{c,11}=|\zeta|^{c-c_0}$,
$\Gamma_{c,22}=|\zeta|^{c_0-c}$,
and $\Gamma_{c,ij}=0$ $(i\neq j)$.
Then, $R(h^0_{c,\ell})$ on $\{|\zeta|\geq 1\}$
is represented by
the following matrix with respect to the frame
$(\vtilde_1,\vtilde_2)$:
\begin{multline}
 \delbar\bigl(
 (\Gamma_c\Hbar\Gamma_c)^{-1}
 \del(\Gamma_c\Hbar\Gamma_c)
 \bigr)
=\delbar(\Gamma_c^{-1}\Hbar^{-1})
 \cdot \Gamma_c^{-1}\del\Gamma\cdot\Hbar\Gamma_c
-\Gamma_c^{-1}\Hbar^{-1}
 (\Gamma_c^{-1}\del\Gamma_c)
 \delbar(\Hbar\Gamma)
 \\
+\delbar\Gamma_c^{-1}
 (\Hbar^{-1}\del\Hbar)\Gamma_c
-\Gamma_c^{-1}\Hbar^{-1}\del\Hbar\cdot
 \delbar\Gamma_c
+\Gamma_c^{-1}
 \delbar(\Hbar^{-1}\del\Hbar)\cdot\Gamma_c
\end{multline}
Because
$\bigl|
 R(h_{c^0,\ell})
\bigr|_{h_{c^0,\ell},g_{\proj^1}}
\leq
 C_{3}
\exp(-\epsilon_{3}|\eta|^{-m-1})$,
we have
\[
 \bigl|
 \delbar(\Hbar^{-1}\del\Hbar)
 \bigr|_{g_{\proj^1}}
\leq C_{4}
 \exp(-\epsilon_{4}|\eta|^{-m-1}),
\]
and hence
$\bigl|
 \Gamma_c^{-1}\delbar(\Hbar^{-1}\del\Hbar)\cdot\Gamma_c
 \bigr|_{g_{\proj^1}}
\leq
 C_{5}
 \exp(-\epsilon_{5}|\eta|^{-m-1})$.
We also have
\begin{equation}
\label{eq;15.7.13.11}
 \delbar(\Gamma_c^{-1}\Hbar^{-1})
 \cdot \Gamma_c^{-1}\del\Gamma_c\cdot\Hbar\Gamma_c
-\Gamma_c^{-1}\Hbar^{-1}
 (\Gamma_c^{-1}\del\Gamma_c)
 \delbar(\Hbar\Gamma_c)
=
-\bigl[
(\Hbar\Gamma_c)^{-1}
 \delbar(\Hbar\Gamma_c),
 (\Hbar\Gamma_c)^{-1}
 (\Gamma_c^{-1}\del\Gamma_c)\Hbar\Gamma_c
\bigr]
\end{equation}
Because the off-diagonal part of
$\Hbar\Gamma_c$
and $\delbar(\Hbar\Gamma_c)$
are dominated by 
$C_{6}\exp(-\epsilon_{6}|\eta|^{-m-1})$,
the term (\ref{eq;15.7.13.11})
is dominated by
$C_{7}\exp(-\epsilon_{7}|\eta|^{-m-1})\,
 d\zeta\,d\zetabar$.
Similarly,
\[
 \delbar\Gamma_c^{-1}
 (\Hbar^{-1}\del\Hbar)\Gamma_c
-\Gamma_c^{-1}\Hbar^{-1}\del\Hbar\cdot
 \delbar\Gamma_c
=-\bigl[
 \Gamma_c^{-1}\delbar\Gamma_c,
 \Gamma_c^{-1}(\Hbar^{-1}\del\Hbar)\Gamma_c
 \bigr]
\]
is dominated by
$C_{8}\exp(-\epsilon_{8}|\eta|^{-m-1})\,
 d\zeta\,d\zetabar$.
Then, the claim of the lemma follows.
\hfill\qed

\vspace{.1in}
We have the self-adjoint endomorphisms
$k_c$ of $(\Etilde_{\ell},h^0_{c,\ell})$
determined by 
$h_{c,\ell}(u_1,u_2)=h^0_{c,\ell}(k_cu_1,u_2)$.
Note that $k_c$ are bounded
with respect to $h^0_{c,\ell}$,
although the estimate might depend on $c$
at this stage.
We also remark that
$\Tr k_c(P)\geq 1$ at any $P\in\cnum$.
The claim of Proposition \ref{prop;15.7.13.41}
is deduced from the following proposition.

\begin{prop}
\label{prop;15.7.13.40}
We have the convergence
$k_c\lrarr\id$ $(c\to c^0)$
uniformly on $\cnum$.
\end{prop}
\pf
Take a large $p>1$.
Let $\|k_c\|_{h^0_{c,\ell},g_{\proj^1},L^p}$
be the $L^p$-norm of $k_c$
with respect to $h^0_{c,\ell}$ and $g_{\proj^1}$.
We set $s_c:=k_c/\|k_c\|_{h^0_{c,\ell},g_{\proj^1},L^p}$.
According to \cite[Lemma 3.1]{Simpson88},
we have the following inequality on $\cnum$:
\begin{equation}
\label{eq;15.7.13.20}
 \sqrt{-1}\Lambda_{g_{\proj^1}}
 \delbar\del\Tr(s_c)\leq
 \Bigl|
\Lambda_{g_{\proj^1}}
 \Tr\Bigl(
 \bigl(R(h^0_{c,\ell})
 +[\thetatilde,\thetatilde^{\dagger}_{h^0_{c,\ell}}]\bigr)
 s_c\Bigr)
 \Bigr|
\end{equation}

We recall the following general lemma,
which is a variant of \cite[Lemma 2.2]{Simpson90}.
\begin{lem}
Let $\varphi$  and $g$ be bounded $\real$-valued
$C^{\infty}$-functions on 
a punctured disc $\{x\in\cnum\,|\,0<|x|<1\}$.
Suppose that 
$-\del_x\del_{\xbar}\varphi\leq g$
holds on $\{x\in\cnum\,|\,0<|x|<1\}$.
Then, the inequality holds on 
$\{x\in\cnum\,|\,|x|<1\}$
in the sense of distributions.
\end{lem}
\pf
We give only an outline of the proof.
We take a $C^{\infty}$-function
$\rho:\real\lrarr\real_{\geq 0}$
satisfying
$\rho(t)=1$ $(t\leq 1)$
and 
$\rho(t)=0$ $(t\geq 2)$.
For any $N>0$,
we set
$\chi_N(x)=\rho\bigl(-N^{-1}\log|x|\bigr)$.
Note that
$\del\chi_N$,
$\delbar\chi_N$
and $\del\delbar\chi_N$
are bounded with respect to
the Poincar\'e metric 
$(\log|x|^2)^{-1}|x|^{-2}dx\,d\xbar$,
which are uniformly for $N$.

Let $f$ be any $\real_{\geq 0}$-valued test function on $\{|x|<1\}$.
The claim of the lemma is the following inequality:
\[
 \int_{|x|<1}
 -\varphi\del_x\del_{\xbar}f\,
 |dx\,d\xbar|
\leq
 \int_{|x|<1}
 g\cdot f\,
  |dx\,d\xbar|
\]
By the assumption,
we have
$\int_{|x|<1}
 -\varphi\del_x\del_{\xbar}(\chi_Nf)\,|dx\,d\xbar|
\leq
 \int_{|x|<1}
 g\cdot \chi_N f\,
 |dx\,d\xbar|$.
It is enough to prove that
$\lim_{N\to\infty}
 \int_{|x|<1}
 \varphi\,
 \bigl(
 \del\delbar(\chi_Nf)
-\chi_N\del\delbar(f)
 \bigr)
=0$.
It follows from the uniform boundedness of
$\del\chi_N$, $\delbar\chi_N$
and $\del\delbar\chi_N$
with respect to the Poincar\'e metric.
\hfill\qed

\vspace{.1in}
In particular, the inequality (\ref{eq;15.7.13.20})
holds on $\proj^1$.
The right hand side of (\ref{eq;15.7.13.20})
is uniformly bounded in $L^p$.
We can find $L_2^p$-functions 
$M_c$ on $\proj^1$
and constants $C_i>0$ $(i=10,11)$
such that the following holds
for any $c\in U$:
\[
 \Lambda_{g_{\proj^1}}
 \delbar\del\bigl(
 \Tr(s_c)
-M_c
 \bigr)
\leq C_{10},
\quad\quad
 \sup|M_c|\leq C_{11}
\]
By \cite[Proposition 2.1]{Simpson88},
we can take constants
$C_i>0$ $(i=12,13)$
such that the following holds
for any $c\in U$:
\[
 \sup_{\proj^1}\bigl(
 \Tr(s_c)-M_c
 \bigr)
\leq
 C_{12}\int_{\proj^1}\bigl|
 \Tr(s_c)-M_c\bigr|
 \vol_{\proj^1}
\leq
 C_{13}
\]
Hence, we can take a constant $C_{14}>0$
such that the following holds for any $c\in U$:
\begin{equation}
\label{eq;15.7.13.30}
 \sup_{\proj^1}|s_c|
\leq
 C_{14}
\end{equation}

Again, according to 
\cite[Proposition 3.1]{Simpson88},
we have
\[
 \sqrt{-1}\Lambda_{g_{\proj^1}}\delbar\del\Tr(s_c)
=\sqrt{-1}\Lambda_{g_{\proj^1}}
 \Tr\Bigl(s_c
 \bigl(R(h^0_{c,\ell}
 +[\thetatilde,\thetatilde^{\dagger}_{h^0_{c,\ell}}])\bigr)
 \Bigr)
-\bigl|s_c^{-1/2}(\delbar+\thetatilde)s_c
 \bigr|^2_{h^0_{c,\ell},g_{\proj^1}}.
\]
Hence, we obtain
the boundedness of 
$\bigl\|s_c^{-1/2}(\delbar+\thetatilde)s_c
 \bigr\|_{h^0_{c,\ell},g_{\proj^1},L^2}$
$(c\in U)$.
Moreover,
we have 
\[
\bigl\|
 s_c^{-1/2}(\delbar+\thetatilde)s_c
 \bigr\|_{h^0_{c,\ell},g_{\proj^1},L^2}\to 0
\quad (c\to c^0).
\]
Because $s_c^{1/2}$ is uniformly bounded
with respect to $h^0_{c,\ell}$,
we obtain the boundedness of
$\bigl\|(\delbar+\thetatilde)s_c
 \bigr\|_{h^0_{c,\ell},g_{\proj^1},L^2}$
$(c\in U)$.
We also have
$\bigl\|(\delbar+\thetatilde)s_c
 \bigr\|_{h^0_{c,\ell},g_{\proj^1},L^2}
\to 0$ $(c\to c^0)$.
In particular, we obtain
the uniform boundedness of $s_c$
in $L_1^2$
with respect to $h^0_{c,\ell}$.

We take any subsequence $s_{c^i}$
which is weakly convergent in $L_1^2$
locally on $\proj^1\setminus\{\infty\}$.
Let $s_{\infty}$ denote the limit.
The sequence $s_{c^i}$ converges to $s_{\infty}$
almost everywhere.
By the uniform boundedness (\ref{eq;15.7.13.30}),
we have the boundedness of $s_{\infty}$.
It also implies
$\|s_{\infty}\|_{L^p}=\lim\|s_{c^i}\|_{L^p}=1$.
In particular, $s_{\infty}\neq 0$.
By the construction,
we have $(\delbar+\thetatilde)s_{\infty}=0$.
Hence, 
$s_{\infty}$ gives a non-zero endomorphism of
the stable filtered Higgs bundle
$(\nbigp^{c^0}_{\ast}\Etilde_{\ell},\thetatilde)$.
It implies that $s_{\infty}$ is 
the multiplication of a non-zero complex number.
In particular,
$\det(s_{\infty})\neq 0$.
It implies that  
$\bigl\|k_{c^i}
 \bigr\|_{L^p,h^0_{c,\ell},g_{\proj^1}}$
are bounded.
We obtain that 
$\bigl\|k_{c}
 \bigr\|_{L^p,h^0_{c,\ell},g_{\proj^1}}$
$(c\in U)$
are bounded.

It implies that
the sequence
$k_{c}$ $(c\in U)$ are bounded in $L_1^2$.
We also have $C_{15}>0$
such that
$\sup_{\proj^1}|k_c|<C_{15}$ $(c\in U)$.
We take any subsequence $k_{c^i}$
which is weakly convergent in $L_1^2$
locally on $\proj^1\setminus\{\infty\}$.
Then, the limit $k_{\infty}$ is the multiplication of
a non-zero positive number.
Because $\det(k_{c})=1$ for any $c$,
we obtain that $k_{\infty}=\id$.
It implies that
$k_c$ $(c\in U)$ is weakly convergent to $\id$
in $L_1^2$ locally on $\proj^1\setminus\{\infty\}$.

By using \cite[Proposition 3.1]{Simpson88},
we obtain
$\sqrt{-1}\Lambda_{g_{\proj^1}}
 \delbar\del\bigl(
 \Tr(k_c)-2
 \bigr)\leq C_{16}$ $(c\in U)$
for a constant $C_{16}>0$.
By \cite[Proposition 2.1]{Simpson88},
we obtain
$\sup_{\proj^1}\bigl(
 \Tr(k_c)-2
 \bigr)
\leq
 C_{17}
 \int_{\proj^1}\bigl|\Tr(k_c)-2\bigr|\dvol_{\proj^1}$
 $(c\in U)$
for a constant $C_{17}>0$.
Note that 
we always have $\Tr(k_c)-2\geq 0$,
and $\Tr(k_c)-2=0$ implies that $k_c=\id$.
Then, we obtain the uniform convergence
$\Tr(k_c)-2\to 0$ $(c\to c_0)$.
It implies the uniform convergence $k_c\to\id$ $(c\to c_0)$.
Thus, the claims of Proposition \ref{prop;15.7.13.40}
and Proposition \ref{prop;15.7.13.41}
are proved.
\hfill\qed

\subsubsection{Behaviour of $b_c$ when $c\to 0$}

Recall that we have the following description
on $\{|\zeta|>1\}$:
\[
 \log\bigl|
 \vtilde_{1}
 \bigr|_{h_{c,\ell}}
=c\log|\zeta|+\log b_c+\rho_{1,c},
\quad\quad
 \log\bigl|\vtilde_{2}\bigr|_{h_{c,\ell}}
=-(\ell+c)\log|\zeta|-\log b_c+\rho_{2,c}.
\]
According to Proposition \ref{prop;15.8.6.20},
we have
$|\rho_{i,c}|\leq
 C_{30}\exp(-\epsilon_{30}|\zeta|^{m+1})$ $(i=1,2)$
for positive constants $C_{30}$  and $\epsilon_{30}$
which are independent of $c$.

\begin{prop}
\label{prop;15.7.13.110}
When $c\to 0$,
we have $b_{c}\to\infty$.
We also have the uniform convergences
$\rho_{i,c}\lrarr 0$ $(i=1,2)$.
\end{prop}
\pf
First, let us study the convergence
of $\rho_{i,c}$.
We begin with a preliminary.
We set $Y_0:=\cnum_w\times\proj^1$
and $D_0:=\cnum_w\times\{\infty\}$.
Let $p:Y_0\lrarr\proj^1$ be the projection.
We have the locally free $\nbigo_{Y_0}(\ast D_0)$-module
$p^{\ast}(L_1\oplus L_2)$.
Let $\overline{v}_i$ be the pull back of $\vtilde_i$.
We have the $\nbigo_Y(\ast D_Y)$-submodule
$\overline{E}_{\ell}\subset
 p^{\ast}(L_1\oplus L_2)$
generated by
$\overline{e}_1=
 \overline{v}_1+w\cdot\overline{v}_2$
and 
$\overline{e}_2=\zeta^{\ell}\overline{v}_2$.

The restriction of
$\overline{E}_{\ell}$ to $\{w\}\times\proj^1$
is denoted by
$\overline{E}_{\ell,w}$.
For each $w$,
we may naturally regard
$\overline{E}_{\ell,w}$
as a subsheaf of 
$\Vtilde=L_1\oplus L_2$.
We have
$\thetatilde(\overline{E}_{\ell,w})
\subset
\overline{E}_{\ell,w}\otimes\Omega^1_{\proj^1}$.
So, we have the family of Higgs bundles
$(\overline{E}_{\ell,w},\thetatilde)$.
The restrictions of $\overline{v}_i$
and $\overline{e}_i$ to $\{w\}\times \proj^1$
are denoted by
$\overline{v}_{i,w}$
and $\overline{e}_{i,w}$,
respectively.
We have
\[
 \thetatilde
 \bigl(\overline{v}_{1,w},\overline{v}_{2,w}\bigr)
=
 \bigl(\overline{v}_{1,w},\overline{v}_{2,w}\bigr)
\left(
 \begin{array}{cc}
 \alpha \zeta^{m}d\zeta & 0 \\
 0 & -\alpha\zeta^m d\zeta
 \end{array}
\right),
\]
\[
 \thetatilde
 \bigl(\overline{e}_{1,w},\overline{e}_{2,w}\bigr)
=
 \bigl(\overline{e}_{1,w},\overline{e}_{2,w}\bigr)
\left(
 \begin{array}{cc}
 \alpha \zeta^{m}d\zeta & 0 \\
 -2w\alpha\zeta^{m-\ell}d\zeta
 & -\alpha\zeta^m d\zeta
 \end{array}
\right).
\]

We have a natural isomorphism
$\overline{E}_{\ell,0}
\simeq
 L_1\oplus \zeta^{\ell}L_2$.
For $w\neq 0$, we pick $w^{1/2}$,
then we have the isomorphism
$\Phi_{w}:(\overline{E}_{\ell,w},\thetatilde)
\simeq
 (\Etilde_{\ell},\thetatilde)$
given by
$w^{-1/2}\overline{e}_{1,w}\longleftrightarrow e_1$
and 
$w^{1/2}\overline{e}_{2,w}\longleftrightarrow e_2$,
i.e.,
$w^{-1/2}\overline{v}_{1,w}\longleftrightarrow v_1$
and 
$w^{1/2}\overline{v}_{2,w}\longleftrightarrow v_2$.

We take a $C^{\infty}$-metric $\overline{h}$
of $\overline{E}_{\ell|\cnum_w\times\cnum_{\zeta}}$
satisfying the following conditions:
\begin{itemize}
\item
On  $\cnum_w\times\{|\zeta|\geq 1\}$,
we have
 $\overline{h}(\overline{v}_1,\overline{v}_1)=|\zeta|^{-2|w|}$,
 $\overline{h}(\overline{v}_2,\overline{v}_2)=|\zeta|^{-2(\ell-|w|)}$,
 and $\overline{h}(\overline{v}_1,\overline{v}_2)=0$.
\item
On $\{w=0\}\times\cnum_{\zeta}$,
we have
 $\overline{h}(\overline{e}_{1,0},\overline{e}_{1,0})=1$,
 $\overline{h}(\overline{e}_{2,0},\overline{e}_{2,0})=1$,
 and
 $\overline{h}(\overline{e}_{1,0},\overline{e}_{2,0})=0$.
\end{itemize}
The restriction of $\overline{h}$ to $\{w\}\times\proj^1$
is denoted by $\overline{h}_w$.
We take a small $\delta>0$,
and consider $U_w:=\{|w|\leq \delta\}$.
We have the following uniform boundedness
on $\cnum_{\zeta}$:
\[
\Bigl|
 R(\overline{h}_w)+
 \bigl[\thetatilde,\thetatilde^{\dagger}_{\overline{h}_{w}}
 \bigr]
\Bigr|_{\overline{h}_w,g_{\proj^1}}
\leq C_{31}
\quad\quad
 (w\in U_{w}).
\]
Moreover,
we have the uniform convergence
$\lim_{w\to 0}
\Bigl(
  R(\overline{h}_w)+
 \bigl[\thetatilde,\thetatilde^{\dagger}_{\overline{h}_{w}}
 \bigr]
\Bigr)=0$.

For any $c$ satisfying $-\delta<c<0$,
we have the self-adjoint endomorphism $k_c$
of $(\overline{E}_{\ell,-c},\overline{h}_{-c})$
determined by 
$\Phi_{-c}^{\ast}h_{c,\ell}(u_1,u_2)=
\overline{h}_{-c}(k_cu_1,u_2)$.
Take a large $p>1$.
Let $\|k_c\|_{\overline{h}_{-c},g_{\proj^1},L^p}$
be the $L^p$-norm of $k_c$
with respect to $\overline{h}_{-c}$ and $g_{\proj^1}$.
We set $s_c:=k_c/\|k_c\|_{\overline{h}_{-c},g_{\proj^1},L^p}$.

Suppose that 
$|\rho_{1,c}|+|\rho_{2,c}|$ is not uniformly convergent to $0$
when $c\to 0$,
and we shall deduce a contradiction.
Under the assumption,
we have a positive number $\delta>0$
and a subsequence $c^j\to 0$
such that 
\begin{equation}
\label{eq;15.8.9.10}
\sup_{|\zeta|>1}\bigl(
 |\rho_{1,c^j}|+|\rho_{2,c^j}|
 \bigr)\geq \delta.
\end{equation}
By the argument in the proof of Proposition \ref{prop;15.7.13.40},
we can assume that $s_{c^j}$ weakly converges to 
a non-zero endomorphism $s_{\infty}$
of $\overline{E}_{\ell,0}$ 
in $L_1^2$
locally on $\proj^1\setminus\{\infty\}$,
such that
(i) $(\delbar+\thetatilde)s_{\infty}=0$,
(ii) $s_{\infty}$ is bounded with respect to $\overline{h}_0$.
So, it gives an endomorphism of
the poly-stable parabolic Higgs bundle
$(\nbigp_{\ast}^{\overline{h}_0}\overline{E}_0,\thetatilde)$.
We obtain that
$s_{\infty}=\alpha_1\id_{L_1}\oplus \alpha_2\id_{\zeta^{\ell}L_2}$.
Here, $\alpha_i$ are non-negative real numbers,
and $(\alpha_1,\alpha_2)\neq (0,0)$.
Suppose that $\alpha_1\neq 0$.
We have the following uniform convergence on 
any compact subset in $\cnum$.
\[
\lim_{j\to\infty}
 \Bigl(
 h_{c^j,\ell}(\vtilde_1,\vtilde_1)
 |\zeta|^{-2c^j}
 \|k_{c^j}\|_{\overline{h}_{-c^j},g_{\proj^1},L^p}^{-1}
 \alpha_1^{-1}
 \Bigr)
=1
\]
It implies that
$\rho_{1,c^j}$ is convergent to $0$
on any compact subset in 
$\{|\zeta|\geq 1\}$.
Together with the uniform estimate
$\bigl|
 \rho_{1,c^j}
 \bigr|\leq C_{30}\exp(-\epsilon_{30}|\zeta|^{m+1})$,
we obtain the uniform convergence
$\lim_{j\to\infty}
 \sup_{|\zeta|>1}|\rho_{1,c^j}|=0$.
We have 
$\bigl|
 \vtilde_1\wedge
 \vtilde_2
 \bigr|_{h_{c,\ell}}=|\zeta|^{-\ell}$.
By (\ref{eq;15.7.13.1}),
we have the following estimate:
\[
\Bigl|
 \bigl|\vtilde_1\bigr|_{h_{c,\ell}}
 \cdot
 \bigl|\vtilde_2\bigr|_{h_{c,\ell}}
-|\zeta|^{-\ell}
\Bigr|
\leq
 C_{32}\exp(-\epsilon_{32}|\zeta|^{m+1})
\]
Here, the constants
$C_{32}$ and $\epsilon_{32}$ are independent of $c$.
Then, we can deduce 
the uniform convergence $\rho_{2,c^j}\to 0$.
Hence, we obtain 
$\lim_{j\to\infty}\sup_{|\zeta|>1}|\rho_{i,c^j}|=0$ $(i=1,2)$
in the case $\alpha_1\neq 0$.
Similarly,
we can deduce the uniform convergences
 $\lim_{j\to\infty}
 \sup_{|\zeta|>1}|\rho_{i,c^j}|=0$ $(i=1,2)$
in the case $\alpha_2\neq 0$.
But, it contradicts with 
(\ref{eq;15.8.9.10}).
Thus, we can conclude the uniform convergences
$\lim_{c\to 0}\rho_{i,c}=0$ $(i=1,2)$.

\vspace{.1in}
Let us study the divergence of $b_c$.
By the convergence of $\rho_{i,c}$,
we have the following on $\{|\zeta|\geq 1\}$:
\begin{equation}
\label{eq;15.8.18.3}
 \lim_{c\to 0}
 \Bigl(
 b_c^{-1}
 \bigl|\vtilde_1\bigr|_{h_{c,\ell}}
 \Bigr)
=1,
\quad\quad
 \lim_{c\to 0}
 \Bigl(
 b_c
 \bigl|\vtilde_2\bigr|_{h_{c,\ell}}
 \Bigr)
=|\zeta|^{-\ell}.
\end{equation}

Suppose that there exists a subsequence
$c^i\to 0$ such that
$b_{c^i}$ are bounded,
and we shall deduce a contradiction.
We may assume the convergence
$b_{c^i}\to \check{b}$.
We shall give a detailed argument in the case $\check{b}=0$.
Later, we shall sketch the argument for
the simpler case $\check{b}\neq 0$.

\vspace{.1in}

We set $Y_1:=\cnum_x\times \proj^1$
and $D_1:=\cnum_x\times\{\infty\}$.
Let $p_1:Y_1\lrarr \proj^1$ be the projection.
We have the locally free $\nbigo_{Y_1}(\ast D_1)$-module
$p_1^{\ast}(L_1\oplus L_2)$.
The pull back of $v_i$
are denoted by 
$\widehat{v}_i$.
Let $\widehat{E}_{\ell}\subset p_1^{\ast}(L_1\oplus L_2)$
generated by 
$\widehat{f}_1=x\widehat{v}_1+\widehat{v}_2$
and 
$\widehat{f}_2=\zeta^{\ell}\widehat{v}_1$.
The restriction of
$\widehat{E}_{\ell}$ to $\{x\}\times\proj^1$
is denoted by
$\widehat{E}_{\ell,x}$.
The restriction of
$\widehat{v}_i$ and $\widehat{f}_i$
to $\{x\}\times\proj^1$
are denoted by
$\widehat{v}_{i,x}$ and $\widehat{f}_{i,x}$.
We have
\[
 \thetatilde(\widehat{v}_{1,x},\widehat{v}_{2,x})
=
(\widehat{v}_{1,x},\widehat{v}_{2,x})
 \left(
 \begin{array}{cc}
 \alpha \zeta^{m}d\zeta & 0 \\
 0 & -\alpha\zeta^md\zeta
 \end{array}
 \right)
\]
\[
  \thetatilde(\widehat{f}_{1,x},\widehat{f}_{2,x})
=(\widehat{f}_{1,x},\widehat{f}_{2,x})
 \left(
 \begin{array}{cc}
 \alpha \zeta^{m}d\zeta & 0\\
 2x\alpha\zeta^{m-\ell}d\zeta
  & -\alpha\zeta^md\zeta
 \end{array}
 \right)
\]
We have
$\widehat{E}_{\ell,0}\simeq
 \zeta^{\ell}L_1\oplus L_2$.
For $x\neq 0$,
taking $x^{1/2}$,
we have the isomorphism
$\Psi_x:(\widehat{E}_{\ell,x},\thetatilde)
\simeq 
 (\Etilde_{\ell},\thetatilde)$
given by the correspondences
$x^{-1/2}\widehat{f}_{1,x}\longleftrightarrow
 v_1+v_2$
and 
$x^{1/2}\widehat{f}_{2,x}\longleftrightarrow
 \zeta^{\ell}v_1$,
i.e.,
$x^{1/2}\widehat{v}_1\longleftrightarrow v_1$
and 
$x^{-1/2}\widehat{v}_2\longleftrightarrow v_2$.

\vspace{.1in}
We have the isomorphisms
of holomorphic vector bundles
$\Upsilon_{x_i}:
 \widehat{E}_{\ell,0}\simeq
 \widehat{E}_{\ell,x_i}$
given by
$\Upsilon_{x_i}(\widehat{f}_{j,0})
=\widehat{f}_{j,x}$.
We shall implicitly identify the vector bundles
$\widehat{E}_{\ell,0}$
and 
$\widehat{E}_{\ell,x_i}$
by $\Upsilon_{x_i}$ in the following argument.

\vspace{.1in}

Let $x_i:=b_{c^i}^2$.
By the assumption,
we have $\lim_{i\to\infty}x_i=0$.
We obtain the harmonic metrics
$\Psi_{x_i}^{\ast}h_{c^i,\ell}$
on $(\widehat{E}_{\ell,x_i},\theta)$.
We take a family of $C^{\infty}$-Hermitian metrics
$h^0_{x_i}$ of $\widehat{E}_{\ell,x_i|\cnum_{\zeta}}$
satisfying the following conditions.
\begin{itemize}
\item
 $h^0_0(\widehat{v}_{1,0},\widehat{v}_{1,0})=1$,
 $h^0_0(\widehat{v}_{2,0},\widehat{v}_{2,0})=|\zeta|^{-2\ell}$,
 and 
 $h^0_0(\widehat{v}_{1,0},\widehat{v}_{2,0})=0$.
\item
 $\Upsilon_{x_i}^{\ast}h^0_{x_i}\to h^0_0$
 in the $C^{\infty}$-sense
 on any compact subset in $\cnum_{\zeta}$.
\item
We have
$h^0_{x_i}(\widehat{v}_{1,x},\widehat{v}_{1,x})=
 |\zeta|^{2c^i}$,
$h^0_{x_i}(\widehat{v}_{2,x},\widehat{v}_{2,x})=
 |\zeta|^{-2c^i-2\ell}$
and 
$h^0_{x_i}(\widehat{v}_{1,x},\widehat{v}_{2,x})=0$
on $\{|\zeta|>1\}$.
\end{itemize}

Let $k_i$ be the self-adjoint endomorphism
of $(\widehat{E}_{\ell,x_i},h^0_{x_i})$
determined by
$\Psi_{x_i}^{\ast}h_{c^i,\ell}=h^0_{x_i}\cdot k_i$.
By (\ref{eq;15.8.18.3}),
we have the convergence of 
$\Upsilon_{x_i}^{\ast}k_i$ to the identity
on any compact subset in $\{|\zeta|>1\}$.
We have the convergence of
$\Upsilon_{x_i}^{\ast}k_i^{-1}$ to the identity
on any compact subset in $\{|\zeta|>1\}$.

\vspace{.1in}

Let us study the convergence 
of $\Upsilon_{x_i}^{\ast}k_i$ and $\Upsilon_{x_i}^{\ast}k_i^{-1}$
on $\{|\zeta|<2\}$.
We have the following uniform estimate:
\[
 \bigl|
 R(h_{c^i,\ell})
\bigr|_{h_{c^i,\ell},g_{\proj^1}}
\leq C_{42}.
\]
Hence, we have a constant $C_{43}>0$
which is independent of $c^i$
such that 
$-\del_{\zeta}\del_{\zetabar}
 \log|\widehat{f}_1|_{h_{c^i,\ell}}
\leq
 C_{43}$
on $\{|\zeta|<2\}$.
Namely,
we have
\[
 -\del_{\zeta}\del_{\zetabar}
 \Bigl(
 \log|\widehat{f}_1|_{h_{c^i,\ell}}
-C_{43}|\zeta|^2
 \Bigr)
\leq 0
\]
We have already known the uniform boundedness of
$|\widehat{f}_1|_{h_{c^i,\ell}}$
on $\{1<|\zeta|<2\}$.
Then, we obtain the uniform boundedness of
$|\widehat{f}_1|_{h_{c^i,\ell}}$
on $\{|\zeta|<2\}$.
Similarly, we obtain the uniform boundedness
of $|\widehat{f}_2|_{h_{c^i,\ell}}$
on $\{|\zeta|<2\}$.
Hence, $k_i$ and $k_i^{-1}$ are uniformly bounded
on $\{|\zeta|<2\}$.

Then, we may assume that
the sequence 
$\Upsilon_{x_i}^{\ast}k_i$ 
is 
weakly convergent to
$k_{\infty}$ in $L_2^p$
on any compact subset in $\cnum_{\zeta}$.
We have the convergence of 
$\Upsilon_{x_i}^{\ast}k_i^{-1}$
to $k_{\infty}^{-1}$.
We obtain a harmonic metric $
h_{\infty}=h^0_0\cdot k_{\infty}$.
By the construction,
we have
$\nbigp^{0}_{\ast}\widehat{E}_{\ell,0}
\subset
 \nbigp^{h_{\infty}}_{\ast}\widehat{E}_{\ell,0}$.
Because 
$\deg(\nbigp^0_{\ast}\widehat{E}_{\ell,0})
=\deg(\nbigp^{h_{\infty}}_{\ast}\widehat{E}_{\ell,0})
=0$,
we obtain 
$\nbigp^{0}_{\ast}\widehat{E}_{\ell,0}
=\nbigp^{h_{\infty}}_{\ast}\widehat{E}_{\ell,0}$
by Lemma \ref{lem;15.8.7.11}.
But,
$(\nbigp^0_{\ast}\widehat{E}_{\ell,0},\thetatilde)$
is not poly-stable
because 
$\nbigp^0_{0}L_2\simeq
 \nbigo_{\proj^1}(\ell)$
and 
$\deg\nbigp^0_{0}L_2=\ell>0$.
Hence, we have deduced a contradiction
from the assumption
$b_{c^i}\to 0$.

\vspace{.1in}

Let us give a sketch of the argument
in the case $\lim b_{c^i}=\check{b}\neq 0$.
Under the assumption,
the sequence $h_{c^i,\ell}$ is convergent on
$\{|\zeta|\geq 1\}$.
As in the previous case,
by taking a subsequence,
we may assume that
the sequence $h_{c^i,\ell}$
is weakly convergent 
in $L_2^p$
on any compact subset in $\cnum_{\zeta}$,
and the limit $h_{\infty}$ is a harmonic metric
of $(\Etilde_{\ell},\theta)$ 
which is adapted to
the filtered bundle $\nbigp^0_{\ast}\Etilde_{\ell}$.
But, $(\nbigp^0_{\ast}\Etilde_{\ell},\theta)$
is not poly-stable.
Hence, we obtain a contradiction 
even in the case $\check{b}\neq 0$.
\hfill\qed

\subsubsection{Convergence of some sequences}

Let $t_i$ be a sequence of positive numbers
such that $t_i\to\infty$.
According to Proposition \ref{prop;15.7.13.41}
and Proposition \ref{prop;15.7.13.110},
we can take a sequence of negative numbers $c_i$
such that 
$c_i\to 0$
and 
$b_{c_i}t_i^{c_i/(m+1)}=1$.
We set $\tau_i:=t_i^{1/2(m+1)}$.

We use the notation in \S\ref{subsection;15.8.6.12}.
We have the isomorphisms
$\Phi_{\tau_i}:
 \varphi_{\tau_i}^{\ast}
 \Etilde_{\ell}
\simeq
 \Etilde_{\ell}$.
Let $h_i$ denote the Hermitian metric of
$\Etilde_{\ell}$
induced by
$\Phi_{\tau_i}$
and $\varphi_{\tau_i}^{\ast}h_{c_i,\ell}$.
Take any $T>0$.
We have the following estimates 
on $\{|\zeta|\geq T\}$
for any $i$ such that $t_i^{1/(m+1)}T>1$,
where $C$ and $\epsilon$ are positive constants
independent of $i$ and $T$:
\[
\Bigl|
\log
 \bigl|
 \vtilde_1
 \bigr|_{h_i}
-\log\bigl(
 |\zeta|^{c_i}t_i^{\ell/2(m+1)}
 \bigr)
\Bigr|
\leq
 C \exp(-\epsilon|\zeta|^{m+1}t_i)
\]
\[
\Bigl|
\log
 \bigl|
 \vtilde_2
 \bigr|_{h_i}
-\log\bigl(
 |\zeta|^{-\ell-c_i}t_i^{-\ell/2(m+1)}
 \bigr)
\Bigr|
\leq
 C \exp(-\epsilon|\zeta|^{m+1}t_i)
\]
\[
 \bigl|
 h_i(\vtilde_1,\vtilde_2)
 \bigr|
\leq 
 C\exp(-\epsilon|\zeta|^{m+1}t_i)
\cdot
 |\vtilde_1|_{h_i} 
\cdot
 |\vtilde_2|_{h_i}
\]

Let $h^{\lim}$ be the Hermitian metric
of $\Vtilde_{|\cnum^{\ast}}$ given by
\[
 h^{\lim}(\vtilde_1,\vtilde_1)=1,
\quad
 h^{\lim}(\vtilde_2,\vtilde_2)=|\zeta|^{-2\ell},
\quad
 h^{\lim}(\vtilde_1,\vtilde_2)=0.
\]

For any $\gamma>0$,
the automorphism $\Psi_{\gamma}$ on $\Vtilde_{|\cnum^{\ast}}$
is given by
$\Psi_{\gamma}=\gamma\id_{L_1}\oplus \gamma^{-1}\id_{L_2}$.
We define
$\Psi_{\gamma}^{\ast}h_{i}(u_1,u_2)
:=h_{i}(\Psi_{\gamma}u_1,\Psi_{\gamma}u_2)$.
The following is clear.
\begin{prop}
\label{prop;15.8.6.101}
Set $\gamma_i:=t_i^{-\ell/2(m+1)}$.
Then, 
we have the convergence
$\lim_{i\to\infty}
 \Psi_{\gamma_i}^{\ast}h_{i}=h^{\lim}$
on $\cnum^{\ast}$.
\hfill\qed
\end{prop}

\subsection{Complement}
\label{subsection;15.7.15.2}

We use the notation in \S\ref{subsection;15.7.28.1}.
Take $0<\kappa<1$.
We take a $C^{\infty}$-function
$\rho:\cnum\lrarr\real_{\geq 0}$
satisfying
$\rho(\zeta)=1$ $(|\zeta|\leq 1/2)$
and 
$\rho(\zeta)=0$ $(|\zeta|\geq 1)$.
We set 
\[
 u_1:=-(\kappa\rho(\zeta)+|\zeta|^{2\ell})^{-1}\zetabar^{\ell}e_2+e_1,
\quad\quad
u_2:=e_2.
\]
They give a $C^{\infty}$-frame of $\Etilde_{\ell}$.
Note that 
we have
$u_1=v_1$ 
and $u_2=\zeta^{\ell}v_2$
on $\{|\zeta|\geq 1\}$.

We take a large integer $L$.
Let $h_{\kappa}$ be the Hermitian metric of 
$\Etilde_{\ell}$
determined by the following conditions:
\[
 |u_2|_{h_{\kappa}}=\kappa^{L},
\quad
 |u_1|_{h_{\kappa}}=\kappa^{-L},
\quad
 h_{\kappa}(u_1,u_2)=0
\]
Let $\nabla_{\kappa}$ be the Chern connection
of $(\Etilde_{\ell},h_{\kappa})$,
and let $R(h_{\kappa})$ denote 
the curvature of $\nabla_{\kappa}$.
Let $\thetatilde^{\dagger}_{\kappa}$
denote the adjoint of $\thetatilde$
with respect to $h_{\kappa}$.
Let $\thetatilde^{\dagger}_{\kappa}$
denote the adjoint of $\thetatilde$ 
with respect to $h_{\kappa}$.

\begin{lem}
On $\bigl\{|\zeta|\leq 1\bigr\}$,
we have 
$\bigl|
 R(h_{\kappa})
 \bigr|_{h_{\kappa},g_{\cnum}}
\leq
 C_{50}\kappa^L$
and 
$\bigl|
 [\thetatilde,\thetatilde^{\dagger}_{\kappa}]
 \bigr|_{h_{\kappa},g_{\cnum}}
\leq
 C_{50}\kappa^L$
for a constant $C_{50}>0$.
On $\bigl\{|\zeta|\geq 1\bigr\}$,
we have 
$R(h_{\kappa})=
 [\thetatilde,\thetatilde^{\dagger}_{\kappa}]=0$.
\end{lem}
\pf
In the proof,
$O(\kappa^L)$
denotes functions which are dominated by
$C_{51}\kappa^L$
for a constant $C_{51}>0$.
The equalities on $\{|\zeta|\geq 1\}$ are clear.
Let us argue the estimates on $\bigl\{|\zeta|\leq 1\bigr\}$.
We have
$\delbar\bigl(\kappa^{-L}u_2\bigr)=0$,
and
\[
 \delbar(\kappa^Lu_1)=
 -\frac{\delbar(\zetabar^{\ell}\rho)\cdot
 \kappa^{2L+1}}{(\kappa\rho(\zeta)+|\zeta|^{2\ell})^2}
 (\kappa^{-L}u_2).
\]
We have
\[
 \nabla_{\kappa}
 \bigl(\kappa^Lu_1,\kappa^{-L}u_2\bigr)
=
 \bigl(\kappa^Lu_1,\kappa^{-L}u_2\bigr)A,
\quad\quad
A=
 \left(
 \begin{array}{cc}
 0 & \frac{\del(\zeta^{\ell}\rho(\zeta))
 \kappa^{2L+1}}{(\kappa\rho(\zeta)+|\zeta|^{2\ell})^2}
 \\
\frac{-\delbar(\zetabar^{\ell}\rho(\zeta))
 \kappa^{2L+1}}{(\kappa\rho(\zeta)+|\zeta|^{2\ell})^2}
 & 0
 \end{array}
 \right)
\]
It is easy to obtain the following estimate
on $\{|\zeta|\leq 1\}$:
\[
 \frac{\delbar(\zetabar^{\ell}\rho(\zeta))
 \kappa^{2L+1}}{(\kappa\rho(\zeta)+|\zeta|^{2\ell})^2}
=O(\kappa^L)d\zetabar,
\quad\quad
  \frac{\del(\zeta^{\ell}\rho(\zeta))
 \kappa^{2L+1}}{(\kappa\rho(\zeta)+|\zeta|^{2\ell})^2}
=O(\kappa^L)d\zeta
\]
We also have the following
on $\{|\zeta|\leq 1\}$:
\[
 d\Bigl( 
 \frac{\delbar(\zetabar^{\ell}\rho(\zeta))
 \kappa^{2L+1}}{(\kappa\rho(\zeta)+|\zeta|^{2\ell})^2} 
 \Bigr)
=O(\kappa^L)\,d\zeta\,d\zetabar
\]
Hence, we have the following
on $\{|\zeta|\leq 1\}$:
\[
 dA+A\wedge A
=O(\kappa^L)\,d\zeta\,d\zetabar
\]
It implies the estimate for $R(h_{\kappa})$
on $\{|\zeta|\leq 1\}$.

\vspace{.1in}

The Higgs field $\thetatilde$ is represented as follows:
\[
 \thetatilde(\kappa^L u_1,\kappa^{-L}u_2)
=(\kappa^Lu_1,\kappa^{-L}u_2)
 \left(
 \begin{array}{cc}
 \alpha \zeta^md\zeta
 & 0 \\
 \mbox{{}}\\
 \frac{-2\kappa^{2L+1} \alpha
 \rho(\zeta)\zeta^{m-\ell}}{\kappa\rho(\zeta)+|\zeta|^{2\ell}}
 d\zeta
 &
 -\alpha \zeta^md\zeta
 \end{array}
 \right)
\]
The adjoint $\thetatilde^{\dagger}_{\kappa}$
is represented as follows:
\[
 \thetatilde_{\kappa}^{\dagger}
 (\kappa^Lu_1,\kappa^{-L}u_2)
=(\kappa^Lu_1,\kappa^{-L}u_2)
 \left(
 \begin{array}{cc}
 \alphabar \zetabar^md\zetabar
 &  \frac{-2\kappa^{2L+1}\alphabar
 \rho(\zeta)\zetabar^{m-\ell}}{\kappa\rho(\zeta)+|\zeta|^{2\ell}}
 d\zetabar
 \\
 \mbox{{}}\\
 0 
 &
 -\alphabar \zetabar^md\zetabar
 \end{array}
 \right)
\]
Hence,
$[\thetatilde,\thetatilde^{\dagger}_{\kappa}]$
is represented as follows:
\[
 \bigl[
 \thetatilde,\thetatilde^{\dagger}_{\kappa}
 \bigr]
 (\kappa^Lu_1,\kappa^{-L}u_2)
=
 (\kappa^Lu_1,\kappa^{-L}u_2)
 A_2
\]
\[
 A_2:=\left(
 \begin{array}{cc}
 4\kappa^{4L+2}|\alpha|^2\rho(\zeta)^2
 \frac{|\zeta|^{2(m-\ell)}}{(\kappa\rho(\zeta)+|\zeta|^{2\ell})^2}
 d\zetabar d\zeta
 & 
-4\kappa^{2L+1}|\alpha|^2\rho(\zeta)
 \frac{|\zeta|^{2m}\zeta^{-\ell}}
 {\kappa\rho(\zeta)+|\zeta|^{2\ell}}
 d\zetabar\,d\zeta
 \\
 \mbox{{}}
 \\
 4\kappa^{2L+1}|\alpha|^2\rho(\zeta)
 \frac{ |\zeta|^{2m}\zeta^{-\ell}}
 {\kappa\rho(\zeta)+|\zeta|^{2\ell}}
 d\zetabar\,d\zeta
 &
 -4\kappa^{4L+2}|\alpha|^2\rho(\zeta)^2
 \frac{|\zeta|^{2(m-\ell)}}{(\kappa\rho(\zeta)+|\zeta|^{2\ell})^2}
 d\zetabar d\zeta
 \end{array}
 \right)
\]
In particular,
we have the following on $\{|\zeta|\leq 1\}$:
\[
 A_2=O\Bigl(
 \frac{\kappa^{4L+2}\rho(\zeta)|\zeta|^{2(m-\ell)}}
 {(\kappa\rho(\zeta)+|\zeta|^{2\ell})^2}
+
 \frac{\kappa^{2L+1}\rho(\zeta)|\zeta|^{2m-\ell}}
 {\kappa\rho(\zeta)+|\zeta|^{2\ell}}
 \Bigr)
 d\zeta\,d\zetabar
=O(\kappa^L)d\zeta\,d\zetabar
\]
Thus, we are done.
\hfill\qed

\subsubsection{Convergence of some sequences}

Suppose that we are give 
a number $a>\ell/2(m+1)$.
Set $\nu:=a-\ell/2(m+1)>0$.
Let $t_i$ be a sequence of positive numbers
such that $t_i\to\infty$.
We set $\kappa_i:=t_i^{-\nu/L}$,
for which we have
$\kappa_i\to 0$.
We have the Hermitian metric $h_{\kappa_i}$
of $\Etilde_{\ell}$ as above.

We use the notation in \S\ref{subsection;15.8.6.12}.
We set $\tau_i:=t_i^{1/2(m+1)}$.
We have the isomorphisms
$\Phi_{\tau_i}:
 \varphi_{\tau_i}^{\ast}
 \Etilde_{\ell}
\simeq
 \Etilde_{\ell}$.
Let $h_i$ be the Hermitian metric of $\Etilde_{\ell}$
induced by
$\varphi_{\tau_i}^{\ast}h_{\kappa_i}$ and $\Phi_{\tau_i}$.

On $\{|\tau|\geq \kappa_i\}$,
we have the following:
\[
|\vtilde_1|_{h_i}
=t_i^{a},
\quad
|\vtilde_2|_{h_i}
=t_i^{-a}|\zeta|^{-\ell},
\quad
 h_i(\vtilde_1,\vtilde_2)=0
\]

Let $h^{\lim}$ be the Hermitian metric
of $\Vtilde_{|\cnum^{\ast}}$ given by
\[
 h^{\lim}(\vtilde_1,\vtilde_1)=1,
\quad
 h^{\lim}(\vtilde_2,\vtilde_2)=|\zeta|^{-2\ell},
\quad
 h^{\lim}(\vtilde_1,\vtilde_2)=0.
\]

For any $\gamma>0$,
the automorphism $\Psi_{\gamma}$ on $\Vtilde_{|\cnum^{\ast}}$
is given by
$\Psi_{\gamma}=\gamma\id_{L_1}\oplus \gamma^{-1}\id_{L_2}$.
We define
$\Psi_{\gamma}^{\ast}h_{i}(u_1,u_2)
:=h_{i}(\Psi_{\gamma}u_1,\Psi_{\gamma}u_2)$.
The following is clear.
\begin{prop}
\label{prop;15.8.6.102}
Set $\gamma_i:=t_i^{-a}$.
Then, 
we have the convergence
$\lim_{i\to\infty}
 \Psi_{\gamma_i}^{\ast}h_{i}=h^{\lim}$
on $\cnum^{\ast}$.
\hfill\qed
\end{prop}

\section{Limiting configurations}

\subsection{A description of 
generically regular semisimple Higgs bundles of rank $2$}
\label{subsection;15.7.12.1}

Let $X$ be a connected complex curve.
Let $(E,\delbar_E,\theta)$ be a Higgs bundle on $X$ 
of rank $2$
which is generically regular semisimple.
For simplicity, we assume that $\tr(\theta)=0$.
Let us consider the case that 
the spectral curve $\Sigma(E,\theta)$ is reducible,
i.e.,
we have a holomorphic $1$-form $\omega\neq 0$ on $X$
such that
$\Sigma(E,\theta)$ is the union of
the images of $\omega$ and $-\omega$.
Such $\omega$ is determined up to the multiplication of
$\pm 1$.

\begin{rem}
Suppose that $\tr(\theta)\neq 0$.
Then, we set
$\theta':=\theta-(\tr\theta/2)\id_E$.
Then, the Higgs bundle $(E,\delbar_{E},\theta')$
satisfies $\tr(\theta')=0$.
Hence, it is enough to study Higgs bundles
such that the trace of the Higgs field is $0$.

Suppose that $\Sigma(E,\theta)$ is irreducible.
We take a normalization 
$\varphi:\Xtilde\lrarr \Sigma(E,\theta)$.
We have the induced morphism
$\varphi_1:\Xtilde\lrarr X$.
The Higgs bundle
$\varphi_1^{\ast}(E,\delbar_E,\theta)$
satisfies the above condition.
If $X$ is compact and $(E,\delbar_E,\theta)$ is stable,
then $\varphi_1^{\ast}(E,\delbar_E,\theta)$
is stable or poly-stable.
So, the study on the irreducible case
can be reduced to the reducible case.
(See also {\rm\S\ref{subsection;16.5.21.20}}.)
\hfill\qed
\end{rem}

Let $Z(\omega)$ denote the zero set of $\omega$.
We have the decomposition of 
the $\nbigo_X\bigl(\ast Z(\omega)\bigr)$-module
\[
 E\otimes\nbigo_X\bigl(\ast Z(\omega)\bigr)
=L'_{\omega}\oplus L'_{-\omega}
\]
corresponding to the decomposition of the spectral curve,
i.e.,
$\theta
=\omega\cdot \id_{L_{\omega}'}
\oplus
 (-\omega)\cdot \id_{L_{-\omega}'}$.
Let $L_{\omega}$ (resp. $L_{-\omega}$)
denote the $\nbigo_X$-module
obtained as the image of $E$ by the induced morphism
$E\lrarr L_{\omega}'$
(resp. $E\lrarr L_{-\omega}$).
We can regard $E$
as an $\nbigo_X$-submodule 
of $L_{\omega}\oplus L_{-\omega}$.

We obtain the $\nbigo_X$-module
$\det(L_{\omega}\oplus L_{-\omega})\big/\det(E)$
whose supports are contained in $Z(\omega)$.
For each $P\in Z(\omega)$,
let $\ell_P$ denote the length of
the stalk of $\det(L_{\omega}\oplus L_{-\omega})/\det(E)$ at $P$.
Because we have the exact sequence
\[
 0\lrarr L_{\omega}\cap E \lrarr
 E\lrarr L_{-\omega}\lrarr 0,
\]
we have 
$L_{\omega}\cap E=
 L_{\omega}\bigl(-\sum_{P\in Z(\omega)}\ell_P\cdot P\bigr)$
in $L_{\omega}\oplus L_{-\omega}$.
Similarly, we have
$L_{-\omega}\cap E=
 L_{-\omega}\bigl(-\sum_{P\in Z(\omega)}\ell_P\cdot P\bigr)$
in $L_{\omega}\oplus L_{-\omega}$.

\paragraph{Local description}

Let $(U_P,z)$ be a small holomorphic coordinate neighbourhood
around $P$ satisfying $z(P)=0$.
Let $s_P$ be any frame of $\det(E)_{|U_P}$.

\begin{lem}
We have local frames $v_{\pm}$ of $L_{\pm \omega}$ 
around $P$
such that 
(i) $e_1=v_++v_-$
and $e_2=z^{\ell_P}v_-$
is a frame of $E_{|U_P}$,
(ii) $e_1\wedge e_2=s_P$.
\end{lem}
\pf
If $\ell_P=0$,
we have $E=L_{\omega}\oplus L_{-\omega}$ on $U_P$,
and hence the claim is clear.
We shall consider the case $\ell_P>0$.
We omit to distinguish the restriction to $U_P$.
We take a frame $e_2'$ of $z^{\ell}L_{-\omega}$.
We take a section $e_1'$ of $E$ such that
$e'_1$ and $e'_2$ give a frame of $E$.
We have the unique decomposition
$e'_1=v_+'+v_-'$,
where $v_{\pm}'$ are sections of $L_{\pm\omega}$.
By the construction,
$L_{\pm\omega}$ is generated by 
the images of $e_1'$ and $e_2'$.
Because $\ell_P>0$,
we can observe that
the image of $v_{\pm}'$ generates $L_{\pm\omega}$,
i.e., $v_{\pm}'$ are frames of $L_{\pm\omega}$.
We may assume that $e_2'=z^{\ell_P}v_-'$.

Let $g$ be the holomorphic function 
determined by $e'_1\wedge e'_2=g\cdot s$.
Because $g(P)\neq 0$,
we can take a holomorphic function $g_1$
such that $g_1^2=g$.
We set $e_i:=g_1^{-1}\cdot e'_i$
and $v_{\pm}:=g_1^{-1}v_{\pm}'$.
Then, they satisfy the desired condition.
\hfill\qed

\vspace{.1in}
Let $m_P$ denote the order of zero of $\omega$ at $P$,
i.e., we have $\omega=z^{m_P}g_P(z)\,dz$ on $U_P$
for a holomorphic function $g_P(z)$ with $g_P(0)\neq 0$.
For the frame $(e_1,e_2)$,
the Higgs field $\theta_{|U_P}$ is described as 
\[
 \theta(e_1,e_2)
=(e_1,e_2)
 \left(
 \begin{array}{cc}
 \omega & 0 \\
 -2z^{-\ell_P}\omega
 & -\omega
 \end{array}
\right).
\]
In particular,
we have $\ell_P\leq m_P$.

\subsection{The limiting configurations of 
stable Higgs bundles}
\label{subsection;15.7.24.120}

Let $X$ be a compact connected Riemann surface.
Let $(E,\delbar_E,\theta)$ be a Higgs bundle of rank $2$ on $X$
such that 
(i) $(E,\delbar_E,\theta)$ is stable,
(ii) $(E,\delbar_E,\theta)$ is generically regular semisimple,
(iii) $\tr(\theta)=0$.
We assume that the spectral curve of $\theta$ is reducible.
We obtain a holomorphic $1$-form $\omega\neq 0$
and the line bundles $L_{\pm\omega}$ 
as in \S\ref{subsection;15.7.12.1}.
We impose that 
$\deg(L_{\omega})\leq \deg(L_{-\omega})$.
We set $L_1:=L_{\omega}$ and $L_2:=L_{-\omega}$.

We set $d_i:=\deg(L_i)$ $(i=1,2)$.
We have
\[
 d_1+d_2-\sum_{P\in Z(\omega)}\ell_P=\deg(E).
\]
The stability condition for $(E,\theta)$
is equivalent to the inequalities
$\deg(L_j\cap E)=d_j-\sum_{P\in Z(\omega)}\ell_P<\deg(E)/2$ 
$(j=1,2)$,
i.e.,
\[
 d_j-\deg(E)/2>0 \quad (j=1,2).
\]

The local scaling factor $a_P$ of $(E,\delbar,\theta)$
at $P\in Z(\omega)$
is defined as follows:
\[
 a_P:=\frac{\ell_P}{2(m_P+1)}.
\]
We have the functions 
$\chi_P:\real_{\geq 0}\lrarr\real_{\leq 0}$ $(P\in Z(\omega))$
given as follows:
\[
 \chi_P(a):=
 \left\{
 \begin{array}{ll}
 (m_P+1)(a-a_P)
 & (0\leq a\leq a_P)
 \\
 0 &
 (a\geq a_P)
 \end{array}
 \right.
\]
We set $\chi_{E,\theta}(a):=\sum_{P\in Z(\omega)}\chi_i(a)$.

\begin{lem}
We have the unique number $a_{E,\theta}$ 
satisfying the conditions
\[
d_1-\frac{\deg(E)}{2}+\chi_{E,\theta}(a_{E,\theta})=0,
\quad
 0\leq a_{E,\theta}<
 \max\bigl\{a_P\,\big|\,P\in Z(\omega)\bigr\}.
\]
\end{lem}
\pf
The function $\chi_{E,\theta}$ is strictly increasing
for $0\leq a\leq a_1:=\max\{a_P\,|\,P\in Z(\omega)\}$,
and we have
$\chi_{E,\theta}(a)=0$ for $a\geq a_1$.
Hence, we have
$d_1-\frac{\deg(E)}{2}+\chi_{E,\theta}(a_1)=d_1-\frac{\deg(E)}{2}>0$.
We have 
$\chi_{E,\theta}(0)=-\sum_{P\in Z(\omega)}\ell_P/2$.
Because of
$d_1+d_2-\deg(E)-\sum_{P\in Z(\omega)}\ell_P=0$
and $d_1\leq d_2$,
we have
$d_1-\frac{\deg(E)}{2}+\chi_{E,\theta}(0)\leq 0$.
Hence, the claim of the lemma follows.
\hfill\qed

\vspace{.1in}

As in \S\ref{subsection;15.7.21.10},
the numbers $-\chi_P(a_{E,\theta})$ $(P\in Z(\omega))$
give a parabolic structure on the line bundle $L_1$.
For any $\vecc=(c_P\,|\,P\in Z(\omega))\in \real^{Z(\omega)}$,
let $\vecn(\vecc)=(n_P(\vecc)\,|\,P\in Z(\omega))
 \in\seisuu^{Z(\omega)}$ 
be determined by the condition
\[
c_P-1<n_P(\vecc)-\chi_P(a_{E,\theta})\leq c_P.
\]
Then, we set 
$\nbigp^{\lim}_{\vecc}(L_1)
=L_1\Bigl(
 \sum_{P\in Z(\omega)}n_P(\vecc)P
 \Bigr)$.
We obtain a filtered bundle
$\nbigp_{\ast}^{\lim}(L_1)=
 \bigl(
 \nbigp^{\lim}_{\vecc}(L_1)\,\big|\,
 \vecc\in\real^{Z(\omega)}
 \bigr)$
over the meromorphic line bundle $L_1\bigl(\ast Z(\omega)\bigr)$.
The parabolic degree of 
$\nbigp^{\lim}_{\ast}(L_1)$ is 
\[
d_1-\sum_{P\in Z(\omega)}(-\chi_P(a_{E,\theta}))
=d_1+\chi_{E,\theta}(a_{E,\theta})=\frac{\deg(E)}{2}.
\]

Similarly,
the numbers 
$\chi_P(a_{E,\theta})+\ell_P$ $(P\in Z(\omega))$
determine a filtered bundle over $L_2(\ast Z(\omega))$.
The parabolic degree of 
$\nbigp^{\lim}_{\ast}(L_2)$ is
\[
d_2-\sum_{P\in Z(\omega)}(\chi_P(a_{E,\theta})+\ell_P)
=d_1-\frac{\deg(E)}{2}+d_2-\sum_{P\in Z(\omega)}\ell_P=
\frac{\deg(E)}{2}.
\]
The direct sum
$\nbigp^{\lim}_{\ast}L_1\oplus\nbigp^{\lim}_{\ast}L_2$ 
is called the limiting configuration of $(E,\delbar_E,\theta)$.

\subsubsection{Hermitian metrics of the limiting configuration}

We fix a Hermitian metric $h_{\det(E)}$
on the determinant line bundle $\det(E)$.
We have the $2$-form $R(h_{\det(E)})$
obtained as the curvature of the Chern connection
of $(\det(E),h_{\det(E)})$.
\begin{lem}
\label{lem;16.5.19.2}
We have Hermitian metrics
$h^{\lim}_{L_j}$ of the holomorphic line bundles
$L_{j|X\setminus Z(\omega)}$
satisfying the following conditions:
\begin{itemize}
\item
The curvature of the Chern connection $\nabla^{\lim}_j$ of
$(L_{j|X\setminus Z(\omega)},h^{\lim}_{L_j})$
is equal to $R(h_{\det(E)})/2$.
\item
For each $P\in Z(\omega)$,
let $(U_P,z_P)$ be any holomorphic coordinate neighbourhood of $P$
with $z_P(P)=0$.
Then,
$|z_P|^{-2\chi_P(a_{E,\theta})}h^{\lim}_{L_1|U_P\setminus\{P\}}$
and
$|z_P|^{2\chi_P(a_{E,\theta})+2\ell_P}h^{\lim}_{L_2|U_P\setminus\{P\}}$
induce Hermitian metrics of $C^{\infty}$-class
on $L_{1|U_P}$ and $L_{2|U_P}$, respectively.
\item
Under the isomorphism
$\det(E)_{|X\setminus Z(\omega)}\simeq
 (L_1\otimes L_2)_{|X\setminus Z(\omega)}$,
we have
$h_{L_1}^{\lim}\otimes h_{L_2}^{\lim}
=h_{\det(E)}$
on $X\setminus Z(\omega)$.
\end{itemize}
\end{lem}
\pf
Because this is standard, we give only a sketch of the proof.
We can take Hermitian metrics $h_{L_j}$ $(j=1,2)$
of $L_{j|X\setminus Z(\omega)}$
such that 
$|z_P|^{-2\chi_P(a_{E,\theta})}h_{L_1|U_P\setminus \{P\}}$
and
$|z_P|^{2\chi_P(a_{E,\theta})+2\ell_P}h_{L_2|U_P\setminus\{P\}}$
induce Hermitian metrics of $C^{\infty}$-class
on $L_{1|U_P}$ and $L_{2|U_P}$, respectively.
Let $R(h_{L_j})$ denote the curvature form
of $(L_{j|X\setminus Z(\omega)},h_{L_j})$.
They naturally induce $2$-forms of $C^{\infty}$-class on $X$,
which are also denoted by $R(h_{L_j})$.
Because 
$\frac{\sqrt{-1}}{2\pi}\int R(h_{L_j})$
is equal to the parabolic degree of $\nbigp^{\lim}_{\ast}L_j$,
we have
$\int\Bigl(R(h_{L_j})-R(h_{\det(E)})/2\Bigr)=0$.
We have $C^{\infty}$-functions $\rho_j$ $(j=1,2)$ on $X$
such that 
$\delbar\del\rho_j=R(h_{L_j})-R(h_{\det(E)})/2$.
We set 
$h^{(1)}_{L_j}:=e^{-\rho_j}h_{L_j}$.
Then, we have
$R(h^{(1)}_{L_j})=R(h_{\det(E)})/2$.
Because
$\det(E)=
 L_1\otimes L_2\otimes
 \nbigo_X\bigl(-\sum \ell_PP\bigr)$,
the tensor product
$h^{(1)}_{L_1}\otimes h^{(2)}_{L_2}$
induces a $C^{\infty}$-Hermitian metric
of $\det(E)$.
By comparison of the curvature,
we have 
$h^{(1)}_{L_1}\otimes h^{(2)}_{L_2}
=\alpha\cdot h_{\det(E)}$
for a positive constant $\alpha$.
Hence, we obtain Hermitian metrics $h^{\lim}_{L_j}$
with the desired property
by adjusting $h^{(1)}_{L_j}$,
for example
by setting
$h^{\lim}_{L_1}=\alpha^{-1}h^{(1)}_{L_1}$
and 
$h^{\lim}_{L_2}=h^{(1)}_{L_2}$.
\hfill\qed

\vspace{.1in}

Such a metric
$h^{\lim}_{E,\theta}:=h^{\lim}_{L_1}\oplus h^{\lim}_{L_2}$
is also called the limiting configuration of
$(E,\delbar_E,\theta)$.
Note that we have the ambiguity
of the actions of automorphisms
$\alpha\id_{L_1}\oplus\alpha^{-1}\id_{L_2}$ $(\alpha>0)$,
i.e.,
the pair of
$\alpha\cdot h^{\lim}_{L_1}$
and
$\alpha^{-1}\cdot h^{\lim}_{L_2}$
satisfies the conditions in Lemma \ref{lem;16.5.19.2}.
But, the induced Chern connection 
$\nabla^{\lim}_{E,\theta}:=
 \nabla^{\lim}_1\oplus \nabla^{\lim}_2$
on 
$E_{|X\setminus Z(\omega)}$ is well defined.
Note that 
$\nabla^{\lim}_{E,\theta}$ is projectively flat,
whose curvature is given by
the multiplication of 
$R(h_{\det(E)})/2$.

\begin{rem}
The metric $h^{\lim}_{E,\theta}$
is also characterized as a Hermitian-Einstein metric
for the Higgs bundle
$\Bigl(
 (L_1,t\omega)\oplus(L_2,-t\omega)
 \Bigr)_{|X\setminus Z(\omega)}$
adapted to the filtered bundle
$\nbigp^{\lim}_{\ast}L_1\oplus\nbigp^{\lim}_{\ast}L_2$
for any $t\neq 0$
such that 
$\det(h^{\lim}_{E,\theta})=h_{\det(E)}$.
Because 
 $(\nbigp^{\lim}_{\ast}L_1,t\omega)\oplus
 (\nbigp^{\lim}_{\ast}L_2,-t\omega)$
is polystable,
we have the ambiguity of the metric $h^{\lim}_{E,\theta}$
by the automorphisms
$\alpha\id_{L_1}\oplus \alpha^{-1}\id_{L_2}$,
as usual.
\hfill\qed
\end{rem}

\subsection{The limiting configuration in complementary cases}

Let $X$ be a compact connected Riemann surface.
Let $(E,\delbar_E,\theta)$ be any Higgs bundle of rank $2$ on $X$
such that 
(i) $(E,\delbar_E,\theta)$ is generically regular semisimple,
(ii) $\tr(\theta)=0$.
We give limiting configurations
in some complementary cases.

\subsubsection{Polystable Higgs bundles}
\label{subsection;16.5.21.1}

Suppose that $(E,\delbar_E,\theta)$ is polystable.
Then, we have the decomposition
$(E,\theta)
=(L_{\omega},\omega)
\oplus
 (L_{-\omega},-\omega)$.
In particular, the spectral curve is reducible,
and we have $\ell_P=0$ in the description 
in \S\ref{subsection;15.7.12.1}.
In this case,
we set $(E,\theta)$
as the limiting configuration.

We take a Hermitian metric $h_{\det(E)}$
on the line bundle $\det(E)$.
We have Hermitian metrics 
$h^{\lim}_{L_{\pm\omega}}$
on $L_{\pm\omega}$ such that
$R(h^{\lim}_{L_{\pm\omega}})
=R(h_{\det(E)})/2$.
We set 
$h^{\lim}_E:=
 h^{\lim}_{L_{\omega}}
 \oplus
 h^{\lim}_{L_{-\omega}}$
such that
$h^{\lim}_{L_{\omega}}
\otimes
 h^{\lim}_{L_{-\omega}}
=h_{\det(E)}$.
The metric $h^{\lim}_{E,\theta}$
is also called limiting configuration.
It is characterized as
a Hermitian-Einstein metric 
on the Higgs bundle $(E,\delbar_E,\theta)$
such that $\det(h^{\lim}_{E,\theta})=h_{\det(E)}$.
We have the ambiguity 
of the metric $h^{\lim}_{E,\theta}$
caused by automorphisms of 
$(E,\theta)$.

The limiting configuration can also be given as
a filtered bundle on $(X,Z(\omega))$
as in the case of stable Higgs bundles
\S\ref{subsection;15.7.24.120}.
We consider the trivial parabolic structures
for $L_{\pm\omega}$ at $P\in Z(\omega)$.
Namely, for any 
$\vecc=(c_P\,|\,P\in Z(\omega))\in\real^{Z(\omega)}$,
let $\vecn(\vecc)=(n_P(\vecc))\in\seisuu^{Z(\omega)}$
be determined by the condition
$c_P-1<n_P(\vecc)\leq c_P$.
Then, we set
$\nbigp^{\lim}_{\vecc}(L_{\pm\omega})
=L_{\pm\omega}\Bigl(
 \sum_{P\in Z(\omega)}n_P(\vecc)P
 \Bigr)$.
The parabolic degrees of
$\nbigp^{\lim}_{\ast}(L_{\pm\omega})$
is $\deg(E)/2$.
The filtered bundle
$\nbigp^{\lim}_{\ast}L_{\omega}\oplus
 \nbigp^{\lim}_{\ast}L_{-\omega}$
is called the limiting configuration of 
$(E,\delbar_E,\theta)$.
\begin{rem}
The metric $h^{\lim}_{E,\theta}$
can be characterized as 
a Hermitian-Einstein metric 
for $(E,\theta)_{|X\setminus Z(\omega)}$
adapted to the filtered bundle
such that
$\det(h^{\lim}_{E,\theta})=h_{\det(E)}$.
\hfill\qed
\end{rem}

\subsubsection{The case where the spectral curve is irreducible}
\label{subsection;16.5.21.20}

Suppose that the spectral curve $\Sigma(E,\theta)$ is irreducible.
It implies that the Higgs bundle  $(E,\delbar_E,\theta)$ is stable.
We take a normalization $\Xtilde\lrarr \Sigma(E,\theta)$.
We have the induced morphism $p:\Xtilde\lrarr X$,
which is a ramified covering of degree $2$.
We have the involution of $\Sigma(E,\theta)$
induced by the multiplication of $-1$ on the cotangent bundle
$T^{\ast}X$.
It induces an involution $\rho$ of $\Xtilde$ over $X$.
We can regard $X$ as the quotient space of $\Xtilde$
by the action of the group $\{1,\rho\}$.

We set 
$(\Etilde,\delbar_{\Etilde},\thetatilde)
:=p^{\ast}(E,\delbar_E,\theta)$.
The spectral curve of
$(\Etilde,\thetatilde)$ is reducible,
i.e.,
it is the union $\Image(\omegatilde)\cup\Image(-\omegatilde)$
for a holomorphic one form $\omegatilde$ on $\Xtilde$.
We have a natural isomorphism
$\rho^{\ast}(\Etilde,\delbar_{\Etilde},\thetatilde)
\simeq
 (\Etilde,\delbar_{\Etilde},\thetatilde)$.
By the construction, we have
$\rho^{\ast}\omegatilde=-\omegatilde$.
We have $Z(\omegatilde)=p^{-1}D(E,\theta)$.
We have the line bundles
$\Ltilde_{\omegatilde}$
and
$\Ltilde_{-\omegatilde}$
on $\Xtilde$
with an inclusion
$\Etilde\lrarr 
 \Ltilde_{\omegatilde}
 \oplus
 \Ltilde_{-\omegatilde}$
as in \S\ref{subsection;15.7.12.1}.
Because $\rho^{\ast}\omegatilde=-\omegatilde$,
we have natural isomorphisms
$\rho^{\ast}\Ltilde_{\pm\omegatilde}
\simeq
 \Ltilde_{\mp\omegatilde}$
such that the following is commutative:
\[
 \begin{CD}
 \rho^{\ast}\Etilde
 @>>>
 \rho^{\ast}\Ltilde_{\omegatilde}
\oplus
 \rho^{\ast}\Ltilde_{-\omegatilde}
 \\
 @V{\simeq}VV @V{\simeq}VV \\
 \Etilde
 @>>>
 \Ltilde_{-\omegatilde}
\oplus
 \Ltilde_{\omegatilde}
 \end{CD}
\]
Note that we have already known that
$(\Etilde,\delbar_{\Etilde},\thetatilde)$
is poly-stable.
Indeed, 
because $(E,\delbar_E,\theta)$ is stable,
we have a Hermitian-Einstein metric $h^{HE}$
for $(E,\delbar_E,\theta)$.
The pull back $p^{\ast}h^{HE}$
is a Hermitian-Einstein metric for
$(\Etilde,\delbar_{\Etilde},\thetatilde)$,
which implies the poly-stability of 
the Higgs bundle.

\paragraph{Stable case}
When $(\Etilde,\delbar_{\Etilde},\thetatilde)$
is stable,
we obtain the limiting configuration
for $(\Etilde,\delbar_{\Etilde},\thetatilde)$
as a filtered bundle on $(\Xtilde,Z(\omegatilde))$
by the procedure in \S\ref{subsection;15.7.24.120}.
Namely, we obtain the filtered line bundles
$\nbigp^{\lim}_{\ast}\Ltilde_{\pm\omegatilde}$
as in \S\ref{subsection;15.7.24.120}.
In this case, we have 
$\deg(\Ltilde_{\omegatilde})=\deg(\Ltilde_{-\omegatilde})$
which implies that $a_{\Etilde,\thetatilde}=0$
and that the parabolic weights are
$-\chi_{P}(a_{\Etilde,\thetatilde})
=\chi_P(a_{\Etilde,\thetatilde})+\ell_P
=\ell_P/2$.
Because $\ell_P=\ell_{\rho(P)}$,
the isomorphisms
$\rho^{\ast}\Ltilde_{\pm\omegatilde}
\simeq
\Ltilde_{\mp\omegatilde}$
induce the isomorphisms
of filtered line bundles
$\rho^{\ast}\nbigp^{\lim}_{\ast}\Ltilde_{\pm\omegatilde}
\simeq
 \nbigp^{\lim}_{\ast}\Ltilde_{\mp\omegatilde}$.
We have the Hermitian metrics 
$h^{\lim}_{\Ltilde_{\pm\omegatilde}}$
of $\Ltilde_{\pm\omegatilde|\Xtilde\setminus Z(\omegatilde)}$
as in Lemma \ref{lem;16.5.19.2}.
We may also impose the condition
$\rho^{\ast}h^{\lim}_{\Ltilde_{\omegatilde}}
=h^{\lim}_{\Ltilde_{-\omegatilde}}$
with which the metrics $h^{\lim}_{\Ltilde_{\pm\omegatilde}}$
are uniquely determined.

We set 
$h^{\lim}_{\Etilde,\thetatilde}:=
 h^{\lim}_{\Ltilde_{\omegatilde}}
\oplus
 h^{\lim}_{\Ltilde_{-\omegatilde}}$
on $\Etilde_{|\Xtilde\setminus Z(\omegatilde)}$.
Because 
$\rho^{\ast}h^{\lim}_{\Etilde,\thetatilde}
=h^{\lim}_{\Etilde,\thetatilde}$,
we have the Hermitian metric
$h^{\lim}_{E,\theta}$ of 
$E_{|X\setminus D(E,\theta)}$
such that
$p^{\ast}h^{\lim}_{E,\theta}
=h^{\lim}_{\Etilde,\thetatilde}$.
The metric $h^{\lim}_{E,\theta}$ 
and the associated Chern connection
$\nabla^{\lim}_{E,\theta}$ of $E_{|X\setminus D(E,\theta)}$
are uniquely determined.

\begin{rem}
The metric $h^{\lim}_{E,\theta}$ is also characterized as follows.
Because the filtered bundle
$\nbigp^{\lim}_{\ast}\Ltilde_{\omega}
\oplus
 \nbigp^{\lim}_{\ast}\Ltilde_{-\omega}$
is equivariant with respect to the action of
$\{1,\rho\}$,
we have the filtered  bundle
$\nbigp^{\lim}_{\ast}E$
on $(X,D(E,\theta))$
obtained as the descent of
$\nbigp^{\lim}_{\ast}\Ltilde_{\omega}
\oplus
 \nbigp^{\lim}_{\ast}\Ltilde_{-\omega}$.
(See {\rm \S\ref{subsection;16.5.21.10}} below
for the descent of filtered bundles in this situation.)
We can easily observe that
the filtered Higgs bundles
$(\nbigp^{\lim}_{\ast}E,t\theta)$ $(t\neq 0)$
are stable,
and that the parabolic degree is $\deg(E)$.
The metric $h^{\lim}_{E,\theta}$ is 
a unique Hermitian-Einstein metric 
for the Higgs bundle 
$(E,\delbar_E,t\theta)_{|X\setminus D(E,\theta)}$
for any $t$,
such that
$\det(h^{\lim}_{E,\theta})=h_{\det(E)}$
adapted to the filtered bundle
$\nbigp^{\lim}_{\ast}E$.
\hfill\qed
\end{rem}

\paragraph{Polystable case}

We have 
$\Etilde=
 \Ltilde_{\omegatilde}
\oplus
 \Ltilde_{-\omegatilde}$.
As in \S\ref{subsection;16.5.21.1},
we take Hermitian metrics
$h^{\lim}_{\Ltilde_{\pm\omegatilde}}$
of $\Ltilde_{\pm\omegatilde}$
satisfying 
$R(h^{\lim}_{\Ltilde_{\pm\omegatilde}})
=p^{\ast}R(h_{\det(E)})/2$
and 
$h^{\lim}_{\Ltilde_{\omegatilde}}
\otimes
 h^{\lim}_{\Ltilde_{-\omegatilde}}
=p^{\ast}h_{\det(E)}$,
and we set
$h^{\lim}_{\Etilde,\thetatilde}:=
 h^{\lim}_{\Ltilde_{\omegatilde}}
\oplus
 h^{\lim}_{\Ltilde_{-\omegatilde}}$.
We also impose the condition
$\rho^{\ast}h^{\lim}_{\Ltilde_{\omegatilde}}
=h^{\lim}_{\Ltilde_{-\omegatilde}}$,
with which the metrics
$h^{\lim}_{\Ltilde_{\pm\omegatilde}}$
and $h^{\lim}_{\Etilde,\thetatilde}$
are uniquely determined.
Because $\rho^{\ast}h^{\lim}_{\Etilde,\thetatilde}$,
we have a unique Hermitian metric
$h^{\lim}_{E,\theta}$ of $E$
such that 
$p^{\ast}h^{\lim}_{E,\theta}
=h^{\lim}_{\Etilde,\thetatilde}$.

The metric $h^{\lim}_{E,\theta}$
is also characterized as follows.
\begin{lem}
\label{lem;16.5.21.22}
For any $t\neq 0$,
$h^{\lim}_{E,\theta}$
is a unique Hermitian-Einstein metric
for the stable Higgs bundle $(E,\delbar_E,t\theta)$
such that $\det(h^{\lim}_{E,\theta})=h_{\det(E)}$.
In particular,
the Hermitian-Einstein metrics
for the Higgs bundles $(E,\delbar_E,t\theta)$
are independent of $t$.
\end{lem}
\pf
The claim is clear by the construction of
$h^{\lim}_{E,\theta}$.
\hfill\qed

\begin{rem}
We have the filtered bundle
$\nbigp^{\lim}_{\ast}\Etilde$
as in {\rm\S\ref{subsection;16.5.21.1}}.
Because it is equivariant with respect to $\rho$,
we obtain a filtered bundle
$\nbigp^{\lim}_{\ast}E$
obtained as the descent of $\nbigp^{\lim}_{\ast}\Etilde$.
The metric $h^{\lim}_{E,\theta}$
is also characterized as a unique Hermitian metric
for the stable Higgs bundle
$(\nbigp^{\lim}_{\ast}E,t\theta)$ $(t\neq 0)$.
\hfill\qed
\end{rem}

\subsubsection{Descent (Appendix)}
\label{subsection;16.5.21.10}

Let $\Xtilde,X,p,\rho$, $D(E,\theta)$ and $Z(\omegatilde)$
be as in \S\ref{subsection;16.5.21.20}.
Let $\Vtilde$ be a locally free 
$\nbigo_{\Xtilde}(\ast Z(\omegatilde))$-module
of finite rank
which is equivariant with respect to $\{1,\rho\}$,
i.e.,
we are given an isomorphism
$\Phi:\rho^{\ast}\Vtilde\simeq \Vtilde$
such that $\Phi\circ\rho^{\ast}\Phi=\id$.
Let $\nbigp_{\ast}\Vtilde$ be 
a filtered bundle over $\Vtilde$,
which is equivariant with respect to $\{1,\rho\}$,
i.e.,
$\rho^{\ast}\nbigp_{\veca}\Vtilde
=\nbigp_{\veca}\Vtilde$
for any $\veca\in\real^{Z(\omegatilde)}$
under the above isomorphism.
In this case,
the decent of $\nbigp_{\ast}\Vtilde$ 
is given as follows.
By the equivariance of $\Vtilde$,
we have a locally free $\nbigo_X(\ast D(E,\theta))$-module
$V$ with an isomorphism
$p^{\ast}V\simeq \Vtilde$.
It is also described as follows.
We have the locally free $\nbigo_X(\ast D(E,\theta))$-module
$p_{\ast}\Vtilde$.
It is equivariant with respect to $\{1,\rho\}$,
where the action of $\{1,\rho\}$ on $X$ is trivial.
Then, $V$ is the invariant part of $p_{\ast}\Vtilde$
with respect to the action.
For any $P\in D(E,\theta)$,
let $q(P)$  denote the number of the set $p^{-1}(P)$,
which are $1$ or $2$.
Let $\veca=(a_P\,|\,P\in D(E,\theta))\in\real^{D(E,\theta)}$.
For $Q\in p^{-1}(P)$,
we set $\atilde_Q:=q(P)\cdot a_P$.
We obtain the locally free $\nbigo_{\Xtilde}$-module
$\nbigp_{\vecatilde}(\Vtilde)$.
It is equivariant with respect to $\{1,\rho\}$.
We obtain the locally free $\nbigo_X$-module
$p_{\ast}\nbigp_{\vecatilde}\Vtilde$.
It is equivariant with respect to $\{1,\rho\}$.
The invariant part is denoted by
$\nbigp_{\veca}V$.
Thus, we obtain a filtered bundle
$\nbigp_{\ast}V$ over $V$,
which is the decent of $\nbigp_{\ast}\Vtilde$.

\section{Convergence to the limiting configurations}

\subsection{Statements}

\subsubsection{General case}
\label{subsection;15.10.13.1}

Let $(E,\delbar_E,\theta)$ be a stable Higgs bundle
of rank $2$
on a compact connected Riemann surface $X$,
such that 
(i) $(E,\delbar_E,\theta)$ is generically regular semisimple,
(ii) the spectral curve is reducible,
(iii) $\tr(\theta)=0$.
We fix a Hermitian metric $h_{\det(E)}$ of $\det(E)$.
We use the notation in \S\ref{subsection;15.7.24.120}.
We have the limiting configuration
$\nbigp_{\ast}^{\lim}L_1
\oplus
 \nbigp_{\ast}^{\lim}L_2$.
We take Hermitian metrics
$h^{\lim}_{L_j}$ for the parabolic line bundle
$\nbigp_{\ast}^{\lim}L_j$
satisfying the condition in Lemma \ref{lem;16.5.19.2}.
We set $h^{\lim}_{E,\theta}:=h^{\lim}_{L_1}\oplus h^{\lim}_{L_2}$.
We have the associated Chern connection 
$\nabla^{\lim}_{E,\theta}$
of $E_{|X\setminus Z(\omega)}$,
which is projectively flat.

For any $t>0$,
the Higgs bundle 
$(E,\delbar_E,t\theta)$ is also stable.
We have the Hermitian-Einstein metrics $h_t$
of the Higgs bundles $(E,\delbar_E,t\theta)$,
i.e.,
$R(h_t)+[t\theta,t\theta^{\dagger}_{h_t}]$
is equal to the multiplication of $R(h_{\det(E)})/2$
according to 
Hitchin \cite{Hitchin-self-duality}
and Simpson \cite{Simpson88}.
We impose that $\det(h_t)=h_{\det(E)}$.
We have the Chern connection
$\nabla_{h_t}$ of $(E,\delbar_E,h_t)$.

For any $\gamma>0$,
let $\Psi_{\gamma}$ denote the automorphism
of $L_1\oplus L_2$
given by $\Psi_{\gamma}=
 \gamma\id_{L_1}\oplus\gamma^{-1}\id_{L_2}$.
We define the metric 
$\Psi_{\gamma}^{\ast}h_{t}$
of $E_{|X\setminus Z(\omega)}$
by 
$\Psi_{\gamma}^{\ast}h_{t}(u_1,u_2)
=h_t(\Psi_{\gamma}u_1,\Psi_{\gamma}u_2)$
for local sections $u_i$ of $E_{|X\setminus Z(\omega)}$.

Take any $Q\in X\setminus Z(\omega)$.
Let $v_Q$ be any frame of $L_{1|Q}$.
We set
\[
 \gamma(t,Q):=
 \left(
 \frac{h^{\lim}_{L_1}(v_Q,v_Q)}{h_{t}(v_Q,v_Q)}
\right)^{1/2}
\]
We shall prove the following theorem
in \S\ref{subsection;16.5.19.22}--\S\ref{subsection;15.7.24.141}.

\begin{thm}
\label{thm;15.7.15.10}
When $t$ goes to $\infty$,
the sequence 
$\Psi_{\gamma(t,Q)}^{\ast}h_t$
converges to 
$h^{\lim}_{E,\theta}$
in the $C^{\infty}$-sense
on any compact subset 
in $X\setminus Z(\omega)$.
\end{thm}

In particular,
we obtain the following convergence of 
unitary connections $\nabla_{h_t}$.

\begin{cor}
\label{cor;16.5.19.10}
The sequence $\nabla_{h_t}$ $(t\to\infty)$
converges to 
$\nabla^{\lim}_{E,\theta}$
on any compact subset 
in $X\setminus Z(\omega)$.
\end{cor}
\pf
If $\deg(E)=0$ and $R(h_{\det(E)})=0$,
according to the estimates in \S\ref{subsection;15.7.17.30},
it is enough to prove that 
the sequence of the Chern connections of
$(L_{j|X\setminus Z(\omega)},h_{t|L_j})$ converges to
the Chern connection of 
$(L_{j|X\setminus Z(\omega)}, h^{\lim}_{L_j})$.
It follows from Theorem \ref{thm;15.7.15.10}.

Let us consider the case where $\deg(E)$ is even.
We have a holomorphic line bundle $L_0$ on $X$
with an isomorphism $\det(E)\simeq L_0\otimes L_0$.
We have the Hermitian metric $h_{L_0}$ of $L_0$
such that $h_{\det(E)}=h_{L_0}\otimes h_{L_0}$
under the isomorphism.
Then, $h_t\otimes h_{L_0}^{-1}$
on $E\otimes L_0^{-1}$
is a harmonic metric of
$(E\otimes L_0^{-1},
 \delbar_{E\otimes L_0^{-1}},
 \theta,h_t\otimes h_{L_0}^{-1})$.
We have the convergence of
$\nabla_{h_t\otimes h_{L_0}^{-1}}$
to the Chern connection of
$\bigl(
 (L_1\otimes L_0^{-1})_{|X\setminus Z(\omega)},
 h^{\lim}_{L_1}\otimes h_{L_0}^{-1}\bigr)
\oplus
 \bigl
 (L_2\otimes L_0^{-1})_{|X\setminus Z(\omega)},
 h^{\lim}_{L_2}\otimes h_{L_0}^{-1}\bigr)
 \bigr)$
by the consideration in the case
$\deg(E)=0$ and $h_{\det(E)}=0$.
Hence, we obtain the convergence of
$\nabla_{h_t}$ to $\nabla^{\lim}_{E,\theta}$
in the case where $\deg(E)$ is even.

Let us consider the case where $\deg(E)$ is odd.
We take a covering map $p:\Xtilde\lrarr X$ of degree $2$
such that $\Xtilde$ is connected.
Note that 
$p^{\ast}(E,\delbar_E,\theta)$ is stable
because
$\deg(p^{\ast}L_i)-\deg(p^{\ast}E)/2>0$.
Because $\deg(p^{\ast}E)$ is even,
we have the convergence of 
$p^{\ast}\nabla_{h_t}$
to $p^{\ast}\nabla^{\lim}_{E,\theta}$.
Hence, we obtain the convergence of
$\nabla_{h_t}$ to $\nabla^{\lim}_{E,\theta}$.
\hfill\qed

\subsubsection{Symmetric case}
\label{subsection;16.5.19.3}

We can deduce a stronger result 
if $X$ and $(E,\delbar_E,\theta)$
are equipped with an extra symmetry.
Suppose that $X$ is equipped with 
a holomorphic non-trivial involution $\rho$,
i.e.,
$\rho$ is an automorphism of $X$
such that $\rho\circ\rho=\id_X$
and $\rho\neq\id_X$.
Let $(E,\delbar_E,\theta)$ be 
as in \S\ref{subsection;15.10.13.1}.
We impose the following additional conditions.
\begin{itemize}
\item
 $(E,\delbar_E,\theta)$ is equivariant 
 with respect to the action of $\{\id_X,\rho\}$.
 Namely,
 we have an isomorphism 
$\upsilon_{\rho}:\rho^{\ast}(E,\delbar_E,\theta)
\simeq
 (E,\delbar_E,\theta)$
 such that
 $\rho^{\ast}\upsilon_{\rho}\circ\upsilon_{\rho}=\id$.
\item
 We have
 $\rho^{\ast}\omega=-\omega$.
\end{itemize}
We impose the condition
$\rho^{\ast}h_{\det(E)}=h_{\det(E)}$
to the metric $h_{\det(E)}$
under the induced isomorphism
$\rho^{\ast}\det(E)\simeq \det(E)$.

The conditions imply that
we have natural isomorphisms
$\rho^{\ast}L_1\simeq L_2$
and $\rho^{\ast}L_2\simeq L_1$
which are compatible with
$\rho^{\ast}E\simeq E$.
Because $\deg(L_1)=\deg(L_2)$,
we have $a_{E,\theta}=0$ and 
$-\chi_P(a_{E,\theta})
=\chi_P(a_{E,\theta})+\ell_P=\ell_P/2$
for any $P\in Z(\omega)$.
We also have
$\ell_P=\ell_{\rho(P)}$.
So, we have natural isomorphisms
$\rho^{\ast}\nbigp_{\ast}^{\lim}L_1
\simeq
 \nbigp^{\lim}_{\ast}L_2$
and 
$\rho^{\ast}\nbigp_{\ast}^{\lim}L_2
\simeq
 \nbigp^{\lim}_{\ast}L_1$
compatible with the isomorphism
$\rho^{\ast}E\simeq E$.
We can impose the additional condition
$\rho^{\ast}h^{\lim}_{L_1}=h^{\lim}_{L_2}$
to the conditions in Lemma \ref{lem;16.5.19.2},
with which the metrics
$h^{\lim}_{L_i}$ are uniquely determined.
We shall prove the following theorem
in \S\ref{subsection;15.7.24.142}.
\begin{thm}
\label{thm;15.7.24.143}
Suppose the symmetric property of $(E,\delbar_E,\theta)$
as above.
For any $t>0$, let $h_t$ be the Hermitian-Einstein metric
for $(E,\delbar_E,t\theta)$
satisfying $\det(h_t)=h_{\det(E)}$.
Then, when $t$ goes to $\infty$,
the sequence $h_t$ is convergent
to $h^{\lim}_{E,\theta}$
in the $C^{\infty}$-sense
on any compact subset in $X\setminus Z(\omega)$.
\end{thm}

\subsubsection{The case where the spectral curve is irreducible}
\label{subsection;16.5.21.30}

As a complement,
we explain how to deduce the results
in the case where the spectral curve is irreducible,
from Theorem \ref{thm;15.7.24.143}.

Let $(E,\delbar_E,\theta)$ be a Higgs bundle of rank $2$
such that 
(i) $(E,\delbar_E,\theta)$ is generically regular semisimple,
(ii) the spectral curve $\Sigma(E,\theta)$ is irreducible,
(iii) $\tr(\theta)=0$.
Note that the Higgs bundle is stable.
We fix a Hermitian metric $h_{\det(E)}$ on $\det(E)$.
We have the Hermitian-Einstein metrics
$h_t$ of $(E,\delbar_E,t\theta)$
such that $\det(h_t)=h_{\det(E)}$.
Recall that we have constructed
a Hermitian metric $h^{\lim}_{E,\theta}$
of $E_{|X\setminus D(E,\theta)}$
in \S\ref{subsection;16.5.21.20}.

\begin{cor}
When $t$ goes to $\infty$,
the sequence $h_t$ converges to 
$h^{\lim}_{E,\theta}$
in the $C^{\infty}$-sense 
on any compact subset in $X\setminus D(E,\theta)$.
\end{cor}
\pf
We use the notation in \S\ref{subsection;16.5.21.20}.
If the Higgs bundle
$(\Etilde,\delbar_{\Etilde},\thetatilde)$ is polystable,
then the metric $h_t$ are independent of $t$
and equal to $h^{\lim}_{E,\theta}$
as remarked in Lemma \ref{lem;16.5.21.22}.
Let us consider the case where
$(\Etilde,\delbar_{\Etilde},\thetatilde)$
is stable.
Then, 
by Theorem \ref{thm;15.7.24.143},
$p^{\ast}h_t$ is convergent to
$h^{\lim}_{\Etilde,\thetatilde}
=p^{\ast}h^{\lim}_{E,\theta}$
in the $C^{\infty}$-sense
on any compact subset in $\Xtilde\setminus Z(\omegatilde)$.
Hence, we obtain the convergence 
of $h_t$ to $h^{\lim}_{E,\theta}$.
\hfill\qed

\subsection{A reduction for the proof of Theorem \ref{thm;15.7.15.10}}
\label{subsection;16.5.19.22}

Let us observe that 
for the proof of Theorem \ref{thm;15.7.15.10}
it is enough to consider the case where
$\deg(E)=0$ and $R(h_{\det(E)})=0$.
The argument already appeared in the proof of
Corollary \ref{cor;16.5.19.10}.

\begin{lem}
\label{lem;16.5.19.21}
Suppose that we have already proved
the claim of Theorem {\rm\ref{thm;15.7.15.10}}
in the case where $\deg(E)=0$ and $R(h_{\det(E)})=0$.
Then, we obtain the claim of Theorem {\rm \ref{thm;15.7.15.10}}
in the general case.
\end{lem}
\pf
Let us consider the case where $\deg(E)$ is even.
We have a holomorphic line bundle $L_0$
with an isomorphism
$L_0^{\otimes 2}\simeq \det(E)$.
We have a Hermitian metric $h_{L_0}$ on $L_0$
such that $h_{L_0}\otimes h_{L_0}=h_{\det(E)}$.
By the assumption and the construction of the limiting configuration,
we obtain that
the sequence
$\Psi_{\gamma(t,Q)}^{\ast}(h_{t}\otimes h^{-1}_{L_0})$
converges to
$h^{\lim}_{E,\theta}\otimes h_{L_0}^{-1}$
in the $C^{\infty}$-sense 
on any compact subset in $X\setminus Z(\omega)$.
Hence,
we obtain the convergence of 
$\Psi_{\gamma(t,Q)}^{\ast}(h_{t})$
to $h^{\lim}_{E,\theta}$.

Let us consider the case where $\deg(E)$ is odd.
We take a covering $p:\Xtilde\lrarr X$ of degree $2$
such that $\Xtilde$ is connected.
Note that 
$p^{\ast}(E,\delbar_E,\theta)$ is stable
because
$\deg(p^{\ast}L_i)-\deg(p^{\ast}E)/2>0$.
We take $\Qtilde$ such that $p(\Qtilde)=Q$.
Because $\deg(p^{\ast}E)$ is even,
we obtain that 
the sequence
$\Psi_{\gamma(t,\Qtilde)}^{\ast}(p^{\ast}h_{t})$
is convergent 
to $p^{\ast}h^{\lim}_{E,\theta}$
in the $C^{\infty}$-sense
on any compact subset in 
$\Xtilde\setminus p^{-1}D(E,\theta)$.
Hence, we obtain the convergence of
$\Psi_{\gamma(t,Q)}^{\ast}h_t$
to $h^{\lim}_{E,\theta}$.
Thus, the proof of Lemma \ref{lem;16.5.19.21}
is finished.
\hfill\qed

\vspace{.1in}
It remains to prove Theorem \ref{thm;15.7.15.10}
in the case where $\deg(E)=0$ and $R(h_{\det(E)})=0$,
which will be established in 
\S\ref{subsection;15.7.24.140}--\ref{subsection;15.7.24.141}.

\subsection{Construction of approximate solutions}
\label{subsection;15.7.24.140}

Let $X$ and $(E,\delbar_E,\theta)$ be the Higgs bundle
as in \S\ref{subsection;15.10.13.1}.
We impose $\deg(E)=0$.
We fix a flat metric $h_{\det(E)}$ on $\det(E)$.

\subsubsection{Rescaling around the zeroes}
\label{subsection;15.7.23.101}

Let $P\in Z(\omega)$.
We take a holomorphic coordinate system $(U_P,z)$
such that
the eigenvalues of $\theta$ are
$\pm d(z^{m_P+1})=\pm(m_P+1)z^{m_{P}}dz$.
Such $z$ is determined up to the multiplication
of a $(m_P+1)$-th square root of $1$.
We take frames $v_i$ of $L_{i|U_P}$ $(i=1,2)$
such that 
(i) $e_1=v_1+v_2$
and $e_2=z^{\ell_P}v_2$ give a frame of 
$E_{|U_P}$,
(ii) $|e_1\wedge e_2|_{h_{\det(E)}}=1$.

We use the Higgs bundle
$(\Etilde_{\ell},\thetatilde)$
in \S\ref{section;15.7.15.1}
by setting $\alpha=m_P+1$.
Let $\varphi_t:U_P\lrarr\cnum$
be given by
$\varphi_t(z)=t^{1/(m_P+1)}z=\zeta$.
We have 
\[
 \varphi_t^{\ast}\thetatilde
 \bigl(
 \varphi_t^{\ast}\etilde_1,
 \varphi_t^{\ast}\etilde_2
 \bigr)
=(\varphi_t^{\ast}\etilde_1,\varphi_t^{\ast}\etilde_2)
 \left(
 \begin{array}{cc}
 t\alpha z^{m_P}dz
 & 0 \\
 -2\alpha t\cdot t^{-\ell_P/(m_P+1)}
 z^{m_P-\ell_P}dz
 & -t\alpha z^{m_P}dz
 \end{array}
 \right).
\]
Hence, 
we have the isomorphism
$\varphi_t^{\ast}(\Etilde_{\ell},\thetatilde)
\simeq
(E,t\theta)_{|U_P}$
given by the following correspondence:
\[
 t^{\ell_P/2(m_P+1)}
 \varphi_t^{\ast}\etilde_1
\longleftrightarrow
 e_1,
\quad\quad
 t^{-\ell_P/2(m_P+1)}\varphi_t^{\ast}\etilde_2
\longleftrightarrow
 e_2
\]
Moreover, we have
$t^{\ell_P/2(m_P+1)}\varphi_t^{\ast}\vtilde_j
\longleftrightarrow v_j$.

\subsubsection{Local constructions around the zeroes}
\label{subsection;15.7.23.110}

In the following,
for a given positive function $\nu$,
let $O(\nu)$ denote a function $f$
such that $|f|\leq C\nu$,
where $C$ is a positive constant
independent of $t$.
Let $\epsilon$ denote small positive numbers
which are independent of $t$.

Let $P\in Z(\omega)$.
Suppose that 
$a_{E,\theta}-\ell_P/2(m_P+1)<0$.
We have the harmonic metric
$h_{\chi_P(a_{E,\theta}),\ell_P}$ 
of $(\Etilde_{\ell},\thetatilde)$
as in \S\ref{subsection;15.7.28.2}.
We obtain harmonic metrics
$h^0_{t,P}:=\varphi_t^{\ast}h_{\chi_P(a_{E,\theta}),\ell_P}$
of $(E,\delbar_E,t\theta)_{|U_P}$.
By construction, we obtain the following
from Proposition \ref{prop;15.8.6.20}:
\begin{lem}
\label{lem;15.8.6.40}
Take $R_{1,P}>0$ such that
$\{|z|\leq R_{1,P}\}\subset U_{P}$.
Take $0<R_{2,P}<R_{1,P}$.
On $\{R_{2,P}\leq |z|<R_{1,P}\}\subset U_P$,
we have
\[
 \bigl|
 v_1
 \bigr|_{h^0_{t,P}}
=t^{a_{E,\theta}}
 |z|^{\chi_P(a_{E,\theta})}
 \cdot
 b_{\chi_P(a_{E,\theta})}
 \Bigl(
1+O\bigl(\exp(-\epsilon |z|^{m_P+1}t)\bigr)
 \Bigr)
\]
\[
 \bigl|
 v_2
 \bigr|_{h^0_{t,P}}
=
 t^{-a_{E,\theta}}
 |z|^{-\ell_P-\chi_P(a_{E,\theta})}
 b_{\chi_P(a_{E,\theta})}^{-1}
 \Bigl(
1+O\bigl(\exp(-\epsilon |z|^{m_P+1}t)\bigr)
 \Bigr)
\]
\[
 h^0_{t,P}\bigl(v_1,v_2\bigr)
=O\bigl(
 \exp(-\epsilon|z|^{m_P+1}t)
 \bigr)
\]
\hfill\qed
\end{lem}

We also have the following lemma.
\begin{lem}
\label{lem;15.7.16.2}
On $\bigl\{R_{2,P}\leq |z|\leq R_{1,P}\bigr\}$,
we have the following:
\[
 \del \log|v_1|^2_{h^0_{t,P}}
=\chi_P(a_{E,\theta})dz/z
+O\bigl(
 \exp(-\epsilon t)
 \bigr)\,dz
=O(1)\,dz
\]
\[
  \del \log|v_2|^2_{h^0_{t,P}}
=-(\ell_P+\chi_P(a_{E,\theta}))dz/z
+O\bigl(
 \exp(-\epsilon t)
 \bigr)\,dz
=O(1)\,dz
\]
\[
 \delbar\del\log|v_j|^2_{h^0_{t,P}}=
 O\bigl(\exp(-\epsilon t)\bigr)\,dz\,d\zbar
\quad 
 (j=1,2)
\]
\[
 \del h^0_{t}(v_1,v_2)=O\bigl(\exp(-\epsilon t)\bigr)\,dz,
\quad\quad
 \delbar h^0_{t}(v_1,v_2)=O\bigl(\exp(-\epsilon t)\bigr)\,d\zbar,
\]
\[
  \del\delbar 
 h^0_{t}(v_1,v_2)=O\bigl(\exp(-\epsilon t)\bigr)\,dz\,d\zbar
\]
\end{lem}
\pf
According to Proposition \ref{prop;15.8.6.20},
we have
$z\del_z\log|v_1|^2_{h^0_{t,P}}
-\chi_{P}(a_{E,\theta})
=O\bigl(
 \exp(-\epsilon |z|^{m+1}t)
 \bigr)$.
Hence, we obtain the estimate for 
$\del\log|v_1|^2_{h^0_{t,P}}$
on the domain.
We obtain the estimate for
$\del\log|v_2|^2_{h^0_{t,P}}$
in a similar way.
We obtain the estimate for
$\delbar\del\log|v_j|^2_{h^0_{t,P}}$
from Lemma \ref{lem;15.7.16.1}.
We obtain the estimate for 
$\del h_t^0(v_1,v_2)$
and $\del\delbar h_t^0(v_1,v_2)$
from Lemma \ref{lem;15.8.18.40}.
\hfill\qed

\vspace{.1in}

Suppose that
$a_{E,\theta}-\ell_P/2(m_P+1)=0$.
If $\ell_P=a_{E,\theta}=0$,
we set
$h_{t_i,P}=\varphi_{t_i}^{\ast}h_{\Etilde_0}$,
where $h_{\Etilde_0}$ is the harmonic metric
given in \S\ref{subsection;15.7.28.10}.
Suppose $\ell_P>0$.
According to Proposition \ref{prop;15.7.13.41}
and Proposition \ref{prop;15.7.13.110},
for a given sequence $t_i\to\infty$,
we can take a sequence of negative numbers
$c_i\to 0$
such that 
\[
 \frac{\log b_{c_i}}{-c_i}
=\frac{\log t_i}{m_P+1},
\quad
\mbox{\rm i.e.,}\quad
b_{c_i}t^{c_i/(m_P+1)}=1
\]
We obtain the sequence of harmonic metrics
$h_{t_i,P}:=\varphi_{t_i}^{\ast}h_{c_i,\ell_P}$
of $(E,\delbar_E,t_i\theta)_{|U_P}$,
where $h_{c_i,\ell_P}$ are given as in 
\S\ref{subsection;15.7.28.2}.
By Proposition \ref{prop;15.8.6.20},
we have the following.
\begin{lem}
Take $0<R_{2,P}<R_{1,P}$ 
as in Lemma {\rm\ref{lem;15.8.6.40}}.
On $\{R_{2,P}\leq|z|\leq R_{1,P}\}$,
we have the following estimates:
\[
 \bigl|
 v_1
 \bigr|_{h^0_{t_i,P}}
=t_i^{a_{E,\theta}}|z|^{c_i}
 \Bigl(
1+O\bigl(\exp(-\epsilon |z|^{m_P+1}t_i)
 \bigr)
 \Bigr)
\]
\[
 \bigl|
 v_2
 \bigr|_{h^0_{t_i,P}}
=
 t_i^{-a_{E,\theta}}
 |z|^{-c_i-\ell_P}
 \Bigl(
1+O\bigl(\exp(-\epsilon |z|^{m_P+1}t_i)\bigr)
 \Bigr)
\]
\[
 h^0_{t_i,P}\bigl(v_1,v_2\bigr)
=O\bigl(
 \exp(-\epsilon|z|^{m_P+1}t_i)
 \bigr)
\]
\hfill\qed
\end{lem}

As in the case of Lemma \ref{lem;15.7.16.2},
we have the following.
\begin{lem}
Take $0<R_{2,P}<R_{1,P}$.
On $\bigl\{R_{2,P}\leq |z|\leq R_{1,P}\bigr\}$,
we have the following:
\[
 \del \log|v_1|^2_{h^0_{t_i,P}}
=c_idz/z
+O\bigl(
 \exp(-\epsilon t_i)
 \bigr)\,dz
=O(1)\,dz
\]
\[
 \del \log|v_2|^2_{h^0_{t_i,P}}
=-(c_i+\ell_P)dz/z
+O\bigl(
 \exp(-\epsilon t_i)
 \bigr)\,dz
=O(1)\,dz
\]
\[
 \delbar\del\log|v_j|^2_{h^0_{t_i,P}}
=O\bigl(
 \exp(-\epsilon t_i)
 \bigr)\,dz\,d\zbar
\]
\[
 \del h^0_{t_i}(v_1,v_2)=O\bigl(\exp(-\epsilon t_i)\bigr)\,dz,
\quad\quad
 \delbar h^0_{t_i}(v_1,v_2)=O\bigl(\exp(-\epsilon t_i)\bigr)\,d\zbar
\]
\[
  \del\delbar 
 h^0_{t_i}(v_1,v_2)=O\bigl(\exp(-\epsilon t_i)\bigr)\,dz\,d\zbar
\]
\hfill\qed
\end{lem}

Suppose that 
$j_P:=a_{E,\theta}-\ell_P/2(m_P+1)>0$.
We use the notation in \S\ref{subsection;15.7.15.2}.
For a given sequence $t_i\to\infty$,
we set
$\kappa_i:=t_i^{-j_P/L}$.
We have the Hermitian metrics $h_{\kappa_i}$
of $\Etilde_{\ell}$ as in \S\ref{subsection;15.7.15.2}.
We obtain the Hermitian metrics
$h^0_{t_i,P}:=\varphi_{t_i}^{\ast}h_{\kappa_i}$
of $E_{|U_P}$.
By construction,
we have the following
on $\{t_i^{-1/(m_P+1)}\leq |z|<1\}$:
\[
 \bigl|v_1\bigr|_{h^0_{t_i,P}}=t_i^{a_{E,\theta}},
 \quad
 \bigl|v_2\bigr|_{h^0_{t_i,P}}=t_i^{-a_{E,\theta}}|z|^{-\ell_P},
\quad
 h^0_{t_i,P}(v_1,v_2)=0
\]

\begin{lem}
\label{lem;15.7.23.100}
Take 
$0<\epsilon_P<\!<(m_P+1)j_P$.
We have the following estimate:
\[
 R(h^0_{t_i,P})
=O\Bigl(t_i^{-j_P+\epsilon_P/(m_P+1)}\Bigr)
 |z|^{\epsilon_P-2}dz\,d\zbar
\]
\[
 \bigl[t_i\theta,(t_i\theta)^{\dagger}_{h^0_{t_i,P}}\bigr]
=O\Bigl(t_i^{-j_P+\epsilon_P/(m_P+1)}\Bigr)
 |z|^{\epsilon_P-2}dz\,d\zbar
\]
\end{lem}
\pf
Because we have
$R(h^0_{t_i,P})=
 [t\theta,(t\theta)^{\dagger}_{h^0_{t_i,P}}]=0$
on $\{|z|\geq t^{-1/(m_P+1)}\}$,
it is enough to argue the estimates
on $|z|<t_i^{-1/(m_P+1)}$.
We have the following:
\[
 R(h^0_{t_i,P})=
 O\Bigl(
 t_i^{-j_P}t_i^{2/(m_P+1)}
 \Bigr)dz\,d\zbar,
\quad
 \bigl[
 t_i\theta,(t_i\theta)^{\dagger}_{h^0_{t_i,P}}
 \bigr]
=O\Bigl(
 t_i^{-j_P}t_i^{2/(m_P+1)}
 \Bigr) dz\,d\zbar
\]
Both of them are dominated by
\[
 O\Bigl(t_i^{-j_P+2/(m_P+1)}\Bigr)|z|^{2-\epsilon_P}
 \bigl(|z|^{\epsilon_P-2}dzd\zbar\bigr)
=O\Bigl(t_i^{-j_P+\epsilon/(m_P+1)}\Bigr)
 |z|^{\epsilon_P-2}dz\,d\zbar
\]
Thus, we obtain the claim of the lemma.
\hfill\qed

\subsubsection{Global construction}
\label{subsection;15.7.24.1}

We take a K\"ahler metric $g_X$ of $X$.
Let $t_i$ be any sequence of positive numbers
going to $\infty$.
We set $\beta_i:=t_i^{a_{E,\theta}}$.
We shall construct 
a family of Hermitian metrics 
$h^0_{t_i}$ of $(E,\delbar_E,\theta)$
with the following property:
\begin{itemize}
\item
 There exists $p>1$ such that 
 the $L^p$-norms of 
 $R(h^0_{t_i})+\bigl[
 t_i\theta,(t_i\theta)^{\dagger}_{h^0_{t_i}}
 \bigr]$ with respect to $g_X$ and $h^0_{t_i}$
 are bounded.
\item
 There exists $C>0$ such that 
 $C^{-1}h^0_{t_i}\leq h^0_{t_i,P}
 \leq C h^0_{t_i}$
 on the neighbourhood $U_P$
 for each $P\in Z(\omega)$,
and that
 $C^{-1}
 h^0_{t_i}
 \leq
 \Psi_{\beta_i}^{\ast}
 h^{\lim}_{E,\theta}
 \leq
 Ch^0_{t_i}$
 on $X\setminus \bigcup_{P\in Z(\omega)} U_P$.
\end{itemize}

Let $P\in Z(\omega)$.
We take $0<R_{2,P}<R_{1,P}$
such that
$\{|z|\leq R_{1,P}\}\subset U_P$.
We take a function 
$\rho_P:\real\lrarr \real_{\geq 0}$
such that
$\rho_P(s)=1$ $(s\leq R_{2,P})$
and
$\rho_P(s)=0$ $(s\geq R_{1,P})$.
On $\{R_{2,P}\leq |z|\leq R_{1,P}\}\subset U_P$,
we define $h_{t_i}^0$ by the following conditions:
\[
 \log h_{t_i}^0(v_j,v_j)
=\rho_{P}(|z|)
 \log h^0_{t_i,P}(v_j,v_j)
+(1-\rho_P(|z|))
 \log h_{E,\theta}^{\lim}(\Psi_{\beta_i}v_j,\Psi_{\beta_i}v_j)
\quad
 (j=1,2)
\]
\[
 h^0_{t_i}(v_1,v_2)
=\rho_P(|z|)h^0_{t_i,P}(v_1,v_2)
\]
Note that
$\log h^0_{t_i,P}(v_j,v_j)
-\log h^{\lim}_{E,\theta}(\Psi_{\beta_i}v_j,\Psi_{\beta_i}v_j)$
are uniformly bounded on $\{R_{2,P}\leq |z|\leq R_{1,P}\}$.
On $\{|z|\leq R_{2,P}\}$,
we set 
$h^0_{t_i}:=h^0_{t_i,P}$.
On $X\setminus\bigcup_{P\in Z(\omega)}\{|z|\leq R_{1,P}\}$,
we set
$h^0_{t_i}:=
 \Psi_{\beta_i}^{\ast}h^{\lim}_{E,\theta}$.
Then, we can check that
the family of the Hermitian metrics
$h^0_{t_i}$ has the desired property
by using the estimates in \S\ref{subsection;15.7.23.110}.

The following lemma is clear by the construction
and Proposition \ref{prop;15.8.6.100},
Proposition \ref{prop;15.8.6.101}
and Proposition \ref{prop;15.8.6.102}
\begin{lem}
\label{lem;15.7.24.10}
The sequence of Hermitian metrics
$\Psi_{\beta_i^{-1}}^{\ast}h^0_{t_i|X\setminus Z(\omega)}$
is convergent in the $C^{\infty}$-sense
on any compact subset in $X\setminus Z(\omega)$.
For the limit
$\htilde^0_{\infty}$,
the decomposition $L_1\oplus L_2$ is orthogonal.
There exists $M_1>0$
such that
$M_1^{-1}h^{\lim}_{E,\theta}
\leq
 \htilde^0_{\infty}
\leq
 M_1h^{\lim}_{E,\theta}$.
In particular, there exists $M_2>0$
with the following property.
\begin{itemize}
\item
 For any neighbourhood $N$ of $Z(\omega)$,
 there exists $i_0(N)$ such that 
 $M_2^{-1} h^{\lim}_{E,\theta}
\leq
\Psi_{\beta_i^{-1}}^{\ast}h_{t_i}^0
\leq 
M_2h^{\lim}_{E,\theta}$
on $X\setminus N$
for any $i\geq i_0(N)$.
\hfill\qed
\end{itemize}
\end{lem}

Let $\rho_i$ be the self-adjoint endomorphisms of 
$\bigl(
 E_{|X\setminus Z(\omega)},h^{\lim}_{E,\theta}
 \bigr)$
determined by
$\Psi_{\beta_i^{-1}}^{\ast}h^0_{t_i}(u,v)=
 h^{\lim}_{E,\theta}(\rho_iu,v)$
for any local sections $u$ and $v$.
We also have the following.
\begin{lem}
\label{lem;15.7.23.200}
The sequence $\rho_i$ are convergent in the $C^{\infty}$-sense
with respect to $h^{\lim}_{E,\theta}$
on any compact subset in $X\setminus Z(\omega)$.
The limit $\rho_{\infty}$ preserves
the decomposition $L_1\oplus L_2$.
We have the boundedness of
$\rho_{\infty}$ and $\rho_{\infty}^{-1}$
with respect to $h^{\lim}_{E,\theta}$.
\hfill\qed
\end{lem}

\subsection{Proof of Theorem \ref{thm;15.7.15.10}}
\label{subsection;15.7.24.141}

We continue to use the notation in 
\S\ref{subsection;15.7.24.140}.

\subsubsection{Boundedness of a modified sequence}

Let $t_i\to\infty$ be any sequence.
It is enough to prove that
we can take a subsequence $t_i'$
such that 
the sequence 
$\Psi^{\ast}_{\gamma(t_i',Q)}h_{t_i'}$
converges to $h^{\lim}_{E,\theta}$.

Let $\Delta_{X}$ be the Laplacian with respect to 
the K\"ahler metric $g_X$ of $X$.
We construct the family of Hermitian metrics of $h_{t_i}^0$ 
on $E$ as in \S\ref{subsection;15.7.24.140}.
Let $k_i$ be the self adjoint endomorphism
of $(E,h_{t_i}^0)$
determined by $h_{t_i}=h_{t_i}^{0}k_i$,
i.e.,
$h_{t_i}(u,v)=h_{t_i}^{0}(k_iu,v)$
for local sections $u$ and $v$.
According to \cite[Proposition 3.1]{Simpson88},
we have the following on $X$:
\[
 \Delta_{X}
 \Tr(k_i)\leq 
 \Bigl|
 \Lambda_{g_{X}}
 \Tr\Bigl(k_i
 \cdot
\bigl(
 R(h^{0}_{t_i})
+[t_i\theta,(t_i\theta)^{\dagger}_{h^0_{t_i}}]
\bigr)
 \Bigr)
 \Bigr|
\]
Let $p>1$  be as in 
\S\ref{subsection;15.7.24.1},
i.e.,
the $L^p$-norms of 
$R(h^0_{t_i})+\bigl[
 t_i\theta,(t_i\theta)^{\dagger}_{h^0_{t_i}}
 \bigr]$ with respect to $g_X$ and $h^0_{t_i}$
 are bounded.
We take $q>1$ such that 
$p^{-1}+q^{-1}<1$.
Set $r:=(p^{-1}+q^{-1})^{-1}$.
Let $\nu_i:=\|k_i\|_{L^q,h^0_{t_i},g_X}$
be the $L^q$-norm of $k_i$
with respect to $h^0_{t_i}$ and $g_X$.
Set $s_i:=\nu_i^{-1}k_i$.
We have
\begin{equation}
 \label{eq;15.7.24.2}
 \Delta_{X}
 \Tr(s_i)\leq 
 \Bigl|
 \Lambda_{g_{X}}
 \Tr\Bigl(s_i
 \cdot
\bigl(
 R(h^{0}_{t_i})
+[t_i\theta,(t_i\theta)^{\dagger}_{h^0_{t_i}}]
\bigr)
 \Bigr)
 \Bigr|
\end{equation}
The $L^r$-norms of the right hand side in (\ref{eq;15.7.24.2})
are bounded for $i$.
So we have a constant $C_1>0$ and 
$L_2^r$-functions $G_i$ such that 
\[
 \Delta_X\bigl(
  \Tr(s_i)-G_i
 \bigr)
\leq C_1,
\quad
 \|G_i\|_{L_2^r}\leq C_1.
\]
Hence, we have $C_2>0$
such that 
$\sup_X\bigl|s_i\bigr|_{h^0_{t_i}}
 \leq C_2$
holds for any $i$.
Again, according to \cite[Proposition 3.1]{Simpson88},
we have the following:
\[
 \Delta_{X}\Tr(s_i)
=\sqrt{-1}\Lambda_{X}
 \Tr\Bigl(
 s_i\cdot
\bigl(
 R(h^{0}_{t_i})
+[t_i\theta,(t_i\theta)^{\dagger}_{h^0_{t_i}}]
\bigr)
 \Bigr)
-\bigl|
 s_i^{-1/2}(\delbar+t_i\theta)s_i
 \bigr|^2_{h^0_{t_i},g_X}
\]
We have a constant $C_3>0$
such that 
$\bigl|
 (\delbar_E+t_i\theta)s_i
 \bigr|^2_{h^0_{t_i},g_X}
\leq
 C_3
 \bigl|
 s_i^{-1/2}(\delbar_E+t_i\theta)s_i
 \bigr|^2_{h^0_{t_i},g_X}$.
Hence, we obtain the following for a constant $C_4>0$:
\[
 \int_X
  \bigl|\delbar_Es_i\bigr|^2_{h^0_{t_i},g_X}
+ \int_X
  \bigl|t_i[\theta, s_i]\bigr|^2_{h^0_{t_i},g_X}
\leq
 C_4
\]

\subsubsection{Weak convergence of a subsequence}
\label{subsection;16.5.19.30}

We consider the sequences of metrics
$\overline{h}_{t_i}:=\Psi_{\beta_i^{-1}}^{\ast}h_{t_i}$
and 
$\overline{h}^0_{t_i}:=\Psi_{\beta_i^{-1}}^{\ast}h^0_{t_i}$
on $E_{|X\setminus Z(\omega)}$.
Let $\overline{k}_i$ be the self-adjoint endomorphism
of $\bigl(E_{|X\setminus Z(\omega)},
 \overline{h}^0_{t_i}\bigr)$
determined by
$\overline{h}_{t_i}
=\overline{h}_{t_i}^{0}\overline{k}_i$.
We have
$\|\overline{k}_i\|_{L^q,\overline{h}^0_{t_i},g_X}=\nu_i$.
Set 
$\overline{s}_i:=
 \nu_i^{-1}\overline{k}_i$.
We have 
\begin{equation}
\label{eq;15.7.24.12}
 \sup_X\bigl|\overline{s}_i\bigr|_{\overline{h}^0_{t_i}}
 \leq C_2.
\end{equation}
We also have
\begin{equation}
 \label{eq;15.7.24.11}
 \int_X
 \bigl|
 \delbar_E\overline{s}_i
 \bigr|^2_{\overline{h}^0_{t_i},g_X}
+\int_X
 \bigl|
 t_i[\theta,\overline{s}_i]
 \bigr|^2_{\overline{h}^0_{t_i},g_X}
\leq
 C_4
\end{equation}

By Lemma \ref{lem;15.7.24.10}
and (\ref{eq;15.7.24.11}),
we may assume that
the sequence $\sbar_i$
is weakly convergent in $L_1^2$
on any compact subset in
$X\setminus Z(\omega)$
with respect to
$g_X$ and $h^{\lim}_{E,\theta}$.
Let $\sbar_{\infty}$ denote the weak limit.

\begin{lem}
\label{lem;15.7.24.21}
$\sbar_{\infty}$ is bounded with respect to
$h^{\lim}_{E,\theta}$.
We have $\sbar_{\infty}\neq 0$.
\end{lem}
\pf
By Lemma \ref{lem;15.7.23.200} and (\ref{eq;15.7.24.12}),
there exists $M_3>0$
with the following property.
\begin{itemize}
\item
 For any neighbourhood $N$ of $Z(\omega)$,
 there exists $i_3(N)$
 such that 
 $|\sbar_i|_{h^{\lim}_{E,\theta}}\leq M_3$
 on $X\setminus N$
 for any $i\geq i_3(N)$.
\end{itemize}
Hence, we have the boundedness
$\bigl|\sbar_{\infty}
 \bigr|_{h^{\lim}_{E,\theta}}\leq M_3$.

Take a small $\delta>0$.
By (\ref{eq;15.7.24.12}),
we have a small neighbourhood 
$N_1$ of $Z(\omega)$
such that
\[
 \int_{X\setminus N_1}
 \bigl|\overline{s}_i\bigr|^q_{\overline{h}^0_{t_i}}
 \dvol_{g_X}
\geq 1-\delta>0
\]
Hence, we obtain
$\int_{X\setminus N_1}
 \bigl|\overline{s}_{\infty}\bigr|^q_{h^{\lim}_{E,\theta}}
 \dvol_{g_X}>0$.
In particular,
we have $\overline{s}_{\infty}\neq 0$.
\hfill\qed

\begin{lem}
\label{lem;15.7.24.20}
We have $[\sbar_{\infty},\theta]=0$.
In particular,
$\sbar_{\infty}$ preserves the decomposition
$L_1\oplus L_2$.
\end{lem}
\pf
We have
$\lim_{i\to\infty}\int_X\bigl|
 [\sbar_i,\theta]
 \bigr|^2_{\overline{h}^0_{t_i},g_X}
=0$
from (\ref{eq;15.7.24.12}),
which implies the claim of the lemma.
\hfill\qed

\subsubsection{Modification}
\label{subsection;16.5.19.31}

Let $\rho_i$ be as in Lemma \ref{lem;15.7.23.200}.
We set $\sbar_i^1:=\rho_i\circ \sbar_i$.
It is self-adjoint with respect to 
$h^{\lim}_{E,\theta}$,
and we have
$\nu_i^{-1}\Psi_{\beta_i^{-1}}^{\ast}h_{t_i}
=h^{\lim}_{E,\theta}\cdot \sbar_i^1$.
The sequence $\sbar_i^1$ is weakly convergent
in $L^2_1$ on any compact subset in $X\setminus Z(\omega)$.
Let $\sbar_{\infty}^1$ denote the weak limit.
We have $\sbar_{\infty}^1=\rho_{\infty}\circ\sbar_{\infty}$.
We obtain the following from
Lemma \ref{lem;15.7.23.200},
Lemma \ref{lem;15.7.24.21}
and Lemma \ref{lem;15.7.24.20}.
\begin{lem}
\label{lem;15.7.24.30}
$\sbar^1_{\infty}$ is bounded
with respect to $h^{\lim}_{E,\theta}$.
We have
$\sbar^1_{\infty}\neq 0$
and 
$[\sbar^1_{\infty},\theta]=0$.
\hfill\qed
\end{lem}

\begin{lem}
\label{lem;15.7.24.31}
We have
$\delbar_E\sbar^1_{\infty}=0$.
\end{lem}
\pf
Applying \cite[Proposition 3.1]{Simpson88}
to $h^{\lim}_{\infty}$ and 
$\nu_i^{-1}\Psi_{\beta_i^{-1}}^{\ast}h_{t_i}$,
we obtain the following 
on $X\setminus Z(\omega)$:
\[
 \bigl|
 (\sbar_i^1)^{-1/2}\delbar_E\sbar^1_i
 \bigr|^2_{h^{\lim}_{E,\theta},g_X}
\leq
\bigl|(\sbar^1_i)^{-1/2}(\delbar+t_i\theta)\sbar^1_i
 \bigr|^2_{h^{\lim}_{E,\theta},g_X}
=- \sqrt{-1}\Lambda_{g_X}\delbar\del\Tr(\sbar^1_i)
\]
We take a $C^{\infty}$-function
$\mu:\real\lrarr \real_{\geq 0}$
such that
$\mu(s)=1$ $(s\leq 1)$
and $\mu(s)=0$ $(s\geq 2)$.
For any sufficiently large real number $T$,
let $\chi_T:X\lrarr \real_{\geq 0}$
be the $C^{\infty}$-function
such that
(i) $\chi_T\equiv 1$  on 
$X\setminus\bigcup_{P\in Z(\omega)}U_P$,
(ii) $\chi_T(z)=\mu\bigl(-T^{-1}\log|z|\bigr)$ 
on the coordinate neighbourhoods $(U_P,z)$
for $P\in Z(\omega)$.

We have a constant $M_4>0$
with the following property.
\begin{itemize}
\item
 For any neighbourhood $N$ of $Z(\omega)$,
 there exists $i_4(N)$ such that 
 $|\sbar^1_i|_{h^{\lim}_{E,\theta}}\leq M_4$
 on $X\setminus N$
 for any $i\geq i_4(N)$.
\end{itemize}
Then, 
we have a constant $C_{10}$
such that 
for any fixed $T>0$
the following holds for a large $i$:
\begin{multline}
 \int\chi_T
 \bigl|\delbar_E\sbar^1_i\bigr|^2_{h^{\lim}_{E,\theta},g_X}
 \dvol_{g_X}
\leq
 C_{10}
 \int
 \chi_T
 \bigl|(\sbar^1_i)^{-1/2}\delbar_E\sbar^1_i
 \bigr|^2_{h^{\lim}_{E,\theta},g_X}
 \dvol_{g_X}
 \\
\leq
 C_{10}
 \left|
\int \chi_T
 \delbar\cdot\del\Tr(\sbar^1_i)
 \right|
=C_{10}
\left|
 \int
\bigl(
 \delbar\del\chi_T
\bigr)
\cdot
 \Tr(\sbar^1_i)
 \right|
\end{multline}
By taking the limit for $i\to\infty$, we obtain
\[
  \int
\chi_T
 \bigl|\delbar_E\sbar^1_{\infty}
 \bigr|^2_{h^{\lim}_{E,\theta},g_X}
 \dvol_{g_X}
\leq
 C_{10}
 \left|
\int 
 \bigl(\delbar\del\chi_T\bigr)
\cdot
 \Tr(\sbar^1_{\infty})
\right|
\]
Note that 
we have already known 
the boundedness of
$\Tr(\sbar^1_{\infty})$.
We also have the uniform boundedness 
of $\delbar\del\chi_T$
with respect to the Poincar\'e like metric
on $X\setminus Z(\omega)$.
Hence, by taking the limit for $T\to\infty$,
we obtain
\[
\int
 \bigl|\delbar_E\sbar^1_{\infty}
 \bigr|^2_{h^{\lim}_{E,\theta},g_X}
 \dvol_{g_X}
\leq 0
\]
Hence, we obtain
$\delbar_E\sbar^1_{\infty}=0$.
\hfill\qed

\vspace{.1in}

By Lemma \ref{lem;15.7.24.30}
and Lemma \ref{lem;15.7.24.31},
we have
$\sbar_{\infty}^1=\alpha_1\id_{L_1}\oplus\alpha_2\id_{L_2}$
for non-negative real numbers $\alpha_i$ $(i=1,2)$.
We have $(\alpha_1,\alpha_2)\neq(0,0)$.

\subsubsection{End of the proof of Theorem \ref{thm;15.7.15.10}}

Suppose that 
$\alpha_1\neq 0$.
Let $u_j$ be local frames of $L_j$ $(j=1,2)$
on a relatively compact open subset
in $X\setminus Z(\omega)$.
We set
$\gamma_i:=\beta_i^{-1}\nu_i^{-1/2}\alpha_1^{-1/2}$.
Then, we have the following:
\begin{equation}
\label{eq;15.7.24.40}
 \lim_{i\to\infty}
 \bigl|u_1
 \bigr|_{h_{t_i}}
 \gamma_i
=\bigl| u_1 \bigr|_{h^{\lim}_{L_1}}
\end{equation}
By the asymptotic orthogonality in \S\ref{section;15.7.16.10},
we have the following on $X\setminus N$,
where $N$ is any neighbourhood of $Z(\omega)$:
\[
 \bigl|u_1\wedge u_2\bigr|_{\det(E)}^2
=
 |u_1|_{h_{t_i}}^2|u_2|_{h_{t_i}}^2
-|h_{t_i}(u_1,u_2)|^2
=|u_1|_{h_{t_i}}^2|u_2|_{h_{t_i}}^2
 \cdot
 \Bigl(
 1+O\bigl(\exp(-\epsilon t_i)\bigr)
 \Bigr)
\]
We also have
$\bigl|
 u_1\wedge u_2
 \bigr|_{\det(E)}
=
\bigl|
 u_1
 \bigr|_{h^{\lim}_{E,\theta}}
 \cdot
\bigl|
 u_2
 \bigr|_{h^{\lim}_{E,\theta}}$.
Hence, we obtain the following
from (\ref{eq;15.7.24.40}):
\begin{equation}
\label{eq;15.7.24.41}
 \lim_{i\to\infty}
 \bigl|u_2\bigr|_{h_{t_i}}
 \gamma_i^{-1}
=\bigl| u_2\bigr|_{h^{\lim}_{L_2}}
\end{equation}
We obtain the following for $j=1,2$
from (\ref{eq;15.7.24.40}) and (\ref{eq;15.7.24.41}):
\[
 \lim_{i\to\infty}
 \Psi_{\gamma_i}^{\ast}
 h_{t_i}(u_j,u_j)
=h^{\lim}_{E,\theta}(u_j,u_j).
\]
We also have the following
on $X\setminus N$,
where $N$ is any neighbourhood of $Z(\omega)$:
\[
 \Psi_{\gamma_i}^{\ast}h_{t_i}(u_1,u_2)
=h_{t_i}(u_1,u_2)
=O(\exp(-\epsilon t_i))
 \cdot
 |u_1|_{h_{t_i}}
 \cdot
 |u_2|_{h_{t_i}}
=O(\exp(-\epsilon t_i))
 \cdot
 |u_1|_{\Psi_{\gamma_i}^{\ast}h_{t_i}}
 \cdot
 |u_2|_{\Psi_{\gamma_i}^{\ast}h_{t_i}}
\]
So, we obtain
$\lim_{i\to\infty}
 \Psi_{\gamma_i}^{\ast}h_{t_i}(u_1,u_2)
=0
=h^{\lim}_{E,\theta}(u_1,u_2)$.
Hence, we have
the convergence of 
$\Psi_{\gamma_i}^{\ast}h_{t_i}$
to $h^{\lim}_{E,\theta}$
in $C^0$
on any compact subset in $X\setminus Z(\omega)$.

By the construction of the sequence
$\gamma(t_i,Q)$,
we have
$\lim_{i\to\infty}\gamma(t_i,Q)\cdot \gamma_i^{-1}=1$.
Hence, 
we have the convergence of 
$\Psi_{\gamma(t_i,Q)}^{\ast}h_{t_i}$
to $h^{\lim}_{E,\theta}$
in $C^0$
on any compact subset in $X\setminus Z(\omega)$.
We can obtain the convergence of the higher derivative
from Corollary \ref{cor;15.7.24.100}.
Thus, we are done 
in the case $\alpha_1\neq 0$.

We can argue the case $\alpha_2\neq 0$
in a similar way.
Thus, the proof of Theorem \ref{thm;15.7.15.10}
is finished.
\hfill\qed

\subsection{Proof of Theorem \ref{thm;15.7.24.143}}
\label{subsection;15.7.24.142}

\subsubsection{Preliminary}
\label{subsection;16.5.19.33}

Let $(E,\delbar_E,\theta)$ be any stable Higgs bundle of rank $2$
on $X$ such that 
(i) $(E,\delbar_E,\theta)$ is generically regular semisimple,
(ii) the spectral curve $\Sigma(E,\theta)$ is reducible,
(iii) $\tr(\theta)=0$.
We have holomorphic line bundles $L_i$ $(i=1,2)$
with an inclusion $E\lrarr L_1\oplus L_2$
as in \S\ref{subsection;15.7.24.120}.
We assume that $\deg(L_1)=\deg(L_2)$.
Then, we have $a_{E,\theta}=0$
and 
$-\chi_P(a_{E,\theta})
=\chi_P(a_{E,\theta})+\ell_P=\ell_P/2$.

We fix a Hermitian metric $h_{\det(E)}$ of $\det(E)$.
For any $t>0$,
we have Hermitian-Einstein metrics $h_t$
for $(E,\delbar_E,t\theta)$ such that
$\det(h_t)=h_{\det(E)}$.
We have the metric $h^{\lim}_{E,\theta}$
as in \S\ref{subsection;15.10.13.1}.

We take any sequence $t_i\to\infty$.
Let $k^{(2)}_{i}$ be the self-adjoint endomorphisms
of $(E_{|X\setminus Z(\omega)},h^{\lim}_{E,\theta})$
determined by
$h_{t_i}=h^{\lim}_{E,\theta}\cdot k^{(2)}_i$.

\begin{lem}
\label{lem;16.5.19.34}
After going to a subsequence $\{i(p)\}\subset \{i\}$
there exists a sequence of positive numbers $\nu_{i(p)}$
such that the sequence
$\nu_{i(p)}^{-1}k^{(2)}_{i(p)}$ weakly converges to
a morphism
$\alpha_1\cdot \id_{L_1}
 \oplus
 \alpha_2\cdot\id_{L_2}$
in $L_1^2$ locally on $X\setminus Z(\omega)$
for non-negative real numbers $\alpha_j$ $(j=1,2)$
with $(\alpha_1,\alpha_2)\neq(0,0)$.
\end{lem}
\pf
Let us consider the case 
where $\deg(E)$ is even.
We take a holomorphic line bundle $L_0$
with an isomorphism
$L_0\otimes L_0\simeq \det(E)$.
We have a Hermitian metric $h_{L_0}$
such that $h_{L_0}\otimes h_{L_0}=h_{\det(E)}$.
We have the harmonic metrics
$h_{E\otimes L_0^{-1},t_i}$
for $(E\otimes L_0^{-1},\delbar_{E\otimes L_0^{-1}},t_i\theta)$.
We have
$h_{E\otimes L_0^{-1},t_i}=h_{t_i}\otimes h_{L_0}^{-1}$
and 
$h^{\lim}_{E\otimes L_0^{-1},\theta}
=h^{\lim}_{E,\theta}\otimes h_{L_0}^{-1}$.
Hence, 
 $k^{(2)}_{i}$ is the self-adjoint endomorphisms
of $((E\otimes L_0^{-1})_{|X\setminus Z(\omega)},
 h^{\lim}_{E\otimes L_0^{-1},\theta})$
determined by
$h_{E\otimes L_0^{-1},t_i}=
 h^{\lim}_{E\otimes L_0^{-1},\theta}\cdot k^{(2)}_i$.
Note that 
$a_{E\otimes L_0^{-1},\theta}=0$ and 
$\beta_i=t_i^{a_{E\otimes L_0^{-1},\theta}}=1$ 
in \S\ref{subsection;15.7.24.1}
and \S\ref{subsection;15.7.24.141}.
Then, the claim of the lemma for $k^{(2)}_i$
has been already observed in 
\S\ref{subsection;16.5.19.30}--\S\ref{subsection;16.5.19.31}.
We can easily reduce the case where $\deg(E)$ is odd
to the case where $\deg(E)$ is even,
by taking the pull back by
a covering $\Xtilde\lrarr X$ of degree $2$.
\hfill\qed

\subsubsection{Proof of Theorem \ref{thm;15.7.24.143}}

We take any sequence $t_i\to\infty$.
It is enough to prove that 
we can take a subsequence $t_i'$
such that $h_{t_i'}$ converges to $h^{\lim}_{E,\theta}$
on any compact subsets in $X\setminus Z(\omega)$.

Let $k^{(2)}_{i}$ be the self-adjoint endomorphisms
of $(E_{|X\setminus Z(\omega)},h^{\lim}_{E,\theta})$
determined by
$h_{t_i}=h^{\lim}_{E,\theta}\cdot k^{(2)}_i$,
as in \S\ref{subsection;16.5.19.33}.
Because
$\rho^{\ast}h^{\lim}_{E,\theta}=h^{\lim}_{E,\theta}$
and $\rho^{\ast}h_{t_i}=h_{t_i}$,
we have $\rho^{\ast}k^{(2)}_i=k^{(2)}_i$.
As remarked in Lemma \ref{lem;16.5.19.34},
by going to a subsequence,
we may assume to have a sequence
of positive numbers $\nu_i$
such that
the sequence
$\nu_i^{-1}k_i^{(2)}$
is weakly convergent to
$\alpha_1\id_{L_1}\oplus\alpha_2\id_{L_2}$
in $L_1^2$ locally on $X\setminus Z(\omega)$,
where $\alpha_i$ are non-negative numbers
such that $(\alpha_1,\alpha_2)\neq(0,0)$.
Because
$\rho^{\ast}\bigl(
 \nu_i^{-1}k^{(2)}_i
 \bigr)
=\nu_i^{-1}k^{(2)}_i$,
$\rho^{\ast}(L_1)=L_2$
and $\rho^{\ast}(L_2)=L_1$,
we have $\alpha_1=\alpha_2$.
In particular, $\alpha_1\cdot\alpha_2\neq 0$.
Because $\det(k_i^{(2)})=1$,
the sequence $\nu_i^{-2}$ converges to
$\alpha_1\cdot\alpha_2$.
In particular, the sequences
$\nu_i$ and $\nu_i^{-1}$ are bounded.

We take a subsequence $h_{t_{i(p)}}$ 
for which the sequence $k^{(2)}_{i(p)}$ is convergent 
to $\beta\id_{L_1}\oplus \beta\id_{L_2}$
for a positive number $\beta$.
But, we have $\det(k^{(2)}_{i(p)})=1$
and hence $\beta=1$,
i.e., $k^{(2)}_{i(p)}$ is convergent to
the identity $\id$, indeed.
Hence, we can conclude that
the sequence $k^{(2)}_i$ is convergent to $\id$,
and the proof of Theorem \ref{thm;15.7.24.143}
is finished.
\hfill\qed

\vspace{.1in}

\noindent
{\em Address\\
Research Institute for Mathematical Sciences,
Kyoto University,
Kyoto 606-8502, Japan,\\
takuro@kurims.kyoto-u.ac.jp
}

\end{document}